\documentclass[letterpaper]{article} % DO NOT CHANGE THIS
\usepackage[preprint]{aaai2027}  % DO NOT CHANGE THIS
\usepackage[hyphens]{url}  % DO NOT CHANGE THIS
\usepackage{graphicx} % DO NOT CHANGE THIS
\usepackage{natbib}  % DO NOT CHANGE THIS AND DO NOT ADD ANY OPTIONS TO IT
\usepackage{caption} % DO NOT CHANGE THIS AND DO NOT ADD ANY OPTIONS TO IT
\usepackage{algorithm}
\usepackage{algorithmicx}
\usepackage{algpseudocode}
\algrenewcommand{\algorithmicrequire}{\textbf{Input:}}
\algrenewcommand{\algorithmicensure}{\textbf{Output:}}

\usepackage{newfloat}
\usepackage{listings}
\DeclareCaptionStyle{ruled}{labelfont=normalfont,labelsep=colon,strut=off} % DO NOT CHANGE THIS
\floatstyle{ruled}
\newfloat{listing}{tb}{lst}{}
\floatname{listing}{Listing}

\usepackage{booktabs}

\usepackage{amsmath}
\usepackage{amssymb}
\usepackage{multirow}
\usepackage{subcaption}
\definecolor{arxivlinkcolor}{rgb}{0.10,0.25,0.60}
\definecolor{arxivcitecolor}{rgb}{0.10,0.40,0.20}
\definecolor{arxivurlcolor}{rgb}{0.55,0.10,0.10}
\usepackage[colorlinks=true,
            linkcolor=arxivlinkcolor,
            citecolor=arxivcitecolor,
            urlcolor=arxivurlcolor,
            breaklinks=true,
            bookmarks=true,
            bookmarksnumbered=true,
            pdfstartview=FitH]{hyperref}
\usepackage{cleveref}  % after hyperref, so that \cref and \Cref link
\makeatletter
\providecommand{\theHALG@line}{}
\renewcommand{\theHALG@line}{alg.\thealgorithm.\arabic{ALG@line}}
\DeclareRobustCommand{\eqref}[1]{\hyperref[#1]{\textup{\tagform@{\ref*{#1}}}}}
\expandafter\let\csname ver@hyperref.sty\endcsname\relax
\makeatother
\newenvironment{summarythm}[1]{%
  \par\addvspace{5pt}\noindent\textbf{#1.}\itshape\ }{\par\addvspace{5pt}}

\newlength{\algwidth}
\newcommand{\algnarrow}{%
  \hsize=\algwidth \linewidth=\algwidth \columnwidth=\algwidth}

\newtheorem{theorem}{Theorem}
\newtheorem{proposition}{Proposition}
\newtheorem{lemma}{Lemma}
\newtheorem{corollary}{Corollary}
\newtheorem{assumption}{Assumption}
\newtheorem{remark}{Remark}

\crefname{algorithm}{Algorithm}{Algorithms}
\crefname{equation}{Equation}{Equations}
\crefname{figure}{Figure}{Figures}
\crefname{table}{Table}{Tables}
\crefname{theorem}{Theorem}{Theorems}
\crefname{assumption}{Assumption}{Assumptions}
\crefname{lemma}{Lemma}{Lemmas}
\crefname{corollary}{Corollary}{Corollaries}
\crefname{remark}{Remark}{Remarks}
\crefname{hypothesis}{Hypothesis}{Hypotheses}

\title{Sign compression for Muon: SignMuon, MuonSign, and the Limits of Error Feedback}

\author{
    Maria Smirnova\equalcontrib,
    Alexey Kravatskiy\equalcontrib
}

\affiliations{
    MIRIAI\\
    Moscow, Russia\\
    smirnova.m@miriai.org, kravatskii.a@miriai.org
}

\newcommand{\arxivpdfauthor}{Maria Smirnova, Alexey Kravatskiy}

\newcommand{\arxivrepourl}{https://github.com/intsystems/signmuon}
\newcommand{\arxivrepo}{\expandafter\url\expandafter{\arxivrepourl}}

\hypersetup{%
  pdftitle={Sign compression for Muon: SignMuon, MuonSign, and the Limits of Error Feedback},
  pdfauthor={\arxivpdfauthor},
  pdfcreator={pdfLaTeX}}

\begin{document}

%%% >>> begin signmuon_body.tex
﻿\maketitle

\begin{abstract}\label{sec:abs}
  SignMuon compresses the Muon update to one bit per parameter by taking its elementwise sign, providing the most direct way to run a matrix-aware optimizer under an extremely low communication budget. It outperforms SignSGD in practice, yet it can ascend even on a linear function. Signing the gradient \emph{before} the Linear Minimization Oracle (LMO), rather than after, does not repair this: we construct a small explicit instance on which sign-before (MuonUSign) and sign-on-both-sides (MuonSign) ascend as well, so no placement of the sign around the oracle descends in general. Error feedback, the standard remedy for a biased compressor, does not rescue SignMuon: when applied to Muon’s output, error feedback can fail for every smoothness constant, step size and momentum. Applied to the gradient, error feedback does work, and EF21-MuonUSign and EF21-MuonSign attain the standard $\mathcal{O}(T^{-1/2})$ rate for the squared gradient norm on smooth nonconvex problems, the latter at one bit in each direction. Experiments then reverse the ordering: across centralized CIFAR-10, federated CIFAR-10 and the nanoGPT speedrun, the strongest compressed method is consistently sign-\emph{after}-the-LMO, precisely the placement we prove divergent, with the provably convergent variants trailing it. Compressing after the LMO, a heuristic, matters more at these scales than the guarantee does.
\end{abstract}

% \textbf{Keywords:} federated learning, stochastic optimization, gradient compression, 1-bit quantization, matrix orthogonalization, LMO, SignSGD, Muon.

\section{Introduction}\label{sec:intro}
Training a deep network across clients consumes bandwidth as well as computation: every round each client transmits a full update, and under a limited link those updates dominate the wall-clock cost. Compression is the standard remedy \cite{bernstein2018signsgd, bernstein2019signsgdmajority, beznosikov2020biased}, and SignSGD is its extreme point: one bit per coordinate, at slight cost in accuracy. SignSGD, however, flattens each weight matrix into a vector, discarding structure that other optimizers exploit. Muon exploits precisely that structure, orthonormalizing the momentum matrix before stepping, and in several settings surpasses tuned adaptive methods \cite{jordan2024muon, bernstein2024old, shah2025practical}. What it transmits, however, is a \emph{dense} matrix at full precision: thirty-two bits per parameter where the budget allows one. Matrix geometry and a one-bit budget are therefore difficult to obtain together.

% In this work, we propose SignMuon,
% % \footnote{Code is available at: \url{https://github.com/intsystems/2026-Project-203/tree/main}},
% an optimization algorithm that applies sign compression to the Linear Minimization Oracle (LMO) update of Muon. SignSGD transmits the sign of the coordinates of the stochastic gradient direction. Our algorithm first computes an update direction using LMO, and then transmits only the sign of this direction. This approach allows for the matrix structure of the parameters to be accounted for and preserves the advantages of Muon, while reducing the amount of transmitted information to one bit per coordinate. Consequently, SignMuon combines the fast convergence of structured optimizers with the communication efficiency of sign-based approaches.
We study the natural ways to combine sign compression with the Muon LMO at one bit per parameter: \textbf{SignMuon}, which signs \emph{after} the LMO, $\operatorname{sign}(\operatorname{LMO}(\cdot))$; \textbf{MuonUSign}, which signs \emph{before}, $\operatorname{LMO}(\operatorname{sign}(\cdot))$; and \textbf{MuonSign}, which signs on \emph{both} sides and, like SignMuon, emits a $\pm1$-valued step, so that uplink and downlink alike cost one bit. All three build Muon's matrix-aware geometry into the step without preserving it intact, and they are not interchangeable: on federated CIFAR-10 (\hyperref[tab:exp_3]{Table~\ref*{tab:exp_3}}) they span $2.8$ accuracy points, in the order after, before, both sides, and only sign-after matches full-precision Muon.

None of them, however, descends in general. We prove that each can \emph{ascend} on a linear objective, the simplest smooth function there is, at every step size and every momentum: SignMuon on a $4\times4$ gradient (\hyperref[th:1]{Theorem~\ref*{th:1}}), MuonUSign and MuonSign on a single $5\times5$ one (\hyperref[th:2]{Theorems~\ref*{th:2}}--\ref{th:3}), both shapes minimal. The standard repair for a biased compressor is error feedback, and for sign-after it is unavailable: applied to the oracle's \emph{output} it fails for every triple $(L,\eta,\mu)$: there is then an $L$-smooth objective on which the method diverges (\hyperref[th:ef_div]{Theorem~\ref*{th:ef_div}}), so no step-size rule built from the smoothness and momentum constants can save it.

What error feedback does repair is the placement that compresses the gradient instead. Adapting EF21-Muon~\citep{gruntkowska2025error} to sign compression gives \textbf{EF21-MuonUSign}, which drives $\min_{t\le T}\mathbb{E}\|\nabla f(\mathbf{X}_t)\|_*^2$ to zero at the standard $\mathcal{O}(T^{-1/2})$ rate for smooth nonconvex problems, at a one-bit uplink, and \textbf{EF21-MuonSign}, which adds a second error-feedback loop on the downlink so that both channels carry one bit for little further cost. Both descend on the counterexamples above. To our knowledge the sign-before and both-sides placements are new, as are the two error-feedback methods; SignMuon itself was proposed concurrently by \citet{mishra2026signmuon}, whose guarantee is proved for the gradient-sign oracle and not for SignMuon (\hyperref[app:mishra]{Appendix~\ref*{app:mishra}}).

% The effectiveness of the proposed optimizer is evaluated on synthetic and practical tasks: centralized image classification on CIFAR-10 with a ResNet-18 backbone, federated image classification on CIFAR-10 with a CNN across several client counts, and on a synthetic smooth convex least-squares problem.
We evaluate all six against Muon, SignSGD, SGD and Adam on centralized CIFAR-10 with a ResNet-18, federated CIFAR-10 at $N=11$ clients, the nanoGPT speedrun, and a synthetic convex quadratic on which the quantity the counterexamples attack can be measured directly. Theory and experiment then disagree. On random quadratics the gradient/step alignment that \hyperref[th:1]{Theorems~\ref*{th:1}}--\ref{th:3} drive negative stays positive for all three placements, so the constructions describe a worst case that ordinary data does not produce; and on all three networks the best compressed method is sign-after, exactly the placement we prove unrepairable, the two provably convergent variants trailing it by several seed spreads. At these scales, compressing only after the oracle outweighs the guarantee. A tuned five-seed federated comparison establishes that ordering. Proofs and every result cited but not stated here are in the appendix, which continues this numbering.

% Experimental results show that applying sign compression to Muon directions reduces communication overhead (by $32\times$) while maintaining optimization rates. In the centralized setting, SignMuon consistently outperforms the baseline SignSGD method in both convergence speed and final accuracy. At the same time, the accuracy of the proposed 1-bit method is only $1-1.5\%$ lower than that of the Muon optimizer, while significantly outperforming basic signed algorithms. In federated setting, integrating SignMuon with client momentum and Majority Vote aggregation method provides stable convergence and robustness to anomalous client updates even under severe communication constraints. Furthermore, we present an analytical example that formalizes the limitations of applying sign compression to LMO-based updates in the deterministic setting and for smooth convex optimization problems. The results of our experiments demonstrate that the proposed algorithm successfully combines LMO-based methods with sign compression without any significant loss in training accuracy.

\section{Related Work}\label{sec:relWork}
\paragraph{Sign compression.} The sign is the most widely used compressor in this literature \citep{alistarh2017qsgd,horvath2023stochastic,beznosikov2020biased}: one bit per coordinate, no index set or scale sent beside it, and a majority vote of signs that is again a sign, so both directions stay at one bit \citep{bernstein2019signsgdmajority}. Its theory is correspondingly well developed. SignSGD \citep{bernstein2018signsgd} needs growing batches to converge; its bias is otherwise repaired by error feedback \citep{seide20141bit,karimireddy2019efsign}, sharpened into EF21 \citep{richtarik2021ef21}, by momentum \citep{cutkosky2020momentum,sun2023signmomentum}, or by randomizing the sign operator \citep{chen2020median,safaryan2021ssdm,jin2024stosign}.

\paragraph{LMO optimizers.} Muon \citep{jordan2024muon} is the spectral-norm instance of a norm-constrained LMO step \citep{bernstein2024old,pethick2025training}, analysed by \citet{li2025note} and \citet{kovalev2025understanding}, then generalized to arbitrary layer norms by Gluon \citep{riabinin2025gluon}.

\paragraph{Compressed and federated Muon.} \citet{gruntkowska2025error} give the error-feedback framework for bidirectionally compressed Muon/Gluon that EF21-MuonUSign and EF21-MuonSign instantiate. Around it: \citet{qian2026communicationefficientgluon} compress Gluon with SARAH-type variance reduction; \citet{takezawa2025fedmuonfederatedlearningbiascorrected} and \citet{zhang2025fedmuon} study federated LMO steps without compression; \citet{ahn2025dion} make the orthonormalization low-rank with error feedback, \citet{gupta2025effective} quantize Muon's optimizer states, \citet{therien2025muloco} quantize the delta of DiLoCo's Muon inner loop \citep{douillard2024diloco} to two bits at next to no loss. Concurrently, \citet{mishra2026signmuon} proposed plain SignMuon with an $\mathcal{O}(1/\sqrt{T})$ guarantee, which \hyperref[app:mishra]{Appendix~\ref*{app:mishra}} shows to be a guarantee for a different method: the rate is proved only for the gradient-sign oracle, which is SignSGD's update and invokes no polar factor, while the one bound covering the polar-sign update carries an unestimated residual that exceeds the bounded quantity exactly where the expected step ascends, as it does throughout the instance of \hyperref[th:1]{Theorem~\ref*{th:1}}. Two further methods relate sign descent to Muon differently. \citet{bolatov2026lionmuon} \emph{alternate} spectral and sign steps rather than composing them; \citet{kravatskiy2025kyfannorms} mix the geometries \emph{inside} the oracle, so S-Muon stays an LMO for an explicit norm, with the attendant theory. Our placements put the sign \emph{around} the oracle, where it is a compressor and not a geometry, so the counterexamples below bear on neither (\hyperref[app:smuon]{Appendix~\ref*{app:smuon}}).

\section{Problem Statement}\label{sec:proState}

We consider the stochastic optimization problem
\begin{equation}
\min_{\mathbf{X} \in \mathcal{X}} \{f(\mathbf{X}) := \mathbb{E}_{\xi \sim \mathcal{D}} [f(\mathbf{X};\xi)]\},
\end{equation}
where $\mathcal{X}$ is the parameter space (e.g., $\mathbb{R}^d$ or $\mathbb{R}^{m \times n}$) and $f(\cdot\,;\xi) : \mathcal{X} \to \mathbb{R}$ is continuously differentiable and possibly non-convex: $f(\mathbf{X};\xi)$ is the loss of a model $\mathbf{X}$ at a data point $\xi \sim \mathcal{D}$. In the federated setting the data are split across clients, and the training problem becomes
\begin{equation}
\min_{\mathbf{X} \in \mathcal{X}} \Bigl\{ f(\mathbf{X}) := \tfrac{1}{N}\textstyle\sum_{j=1}^N f_j(\mathbf{X})\Bigr\},
\label{fed_eq}
\end{equation}
where $N \ge 1$ is the number of clients and $f_j(\mathbf{X})=\mathbb{E}_{\xi_j\sim\mathcal{D}_j}[f_j(\mathbf{X};\xi_j)]$ is the loss on the data $\mathcal{D}_j$ held by client $j \in [N] := \{1, \dots, N\}$.

The model is a tuple of layers $\mathbf{X} = [\mathbf{X}_1,\dots,\mathbf{X}_p]$, $\mathbf{X}_i \in \mathbb{R}^{m_i\times n_i}$ with $m_i,n_i\ge1$, with $\nabla_i$ the gradient block of layer $i$; iterates carry a time index, $\mathbf{X}_t$, and $\mathbf{X}_{t,i}$ is layer $i$ of iterate $t$. This is Gluon's parameter space \citep{riabinin2025gluon} under a single geometry: each layer carries the spectral norm $\|\cdot\|_{2\to2}$, whose dual is the nuclear norm $\|\cdot\|_{*}$, and we define $\|\nabla f(\mathbf{X})\|_*^2 := \sum_i\|\nabla_i f(\mathbf{X})\|_*^2$; a single matrix parameter is the case $p=1$. A vector parameter is a block of width one, on which the spectral norm is Euclidean and the oracle only normalizes, so all three methods of \hyperref[sec:theory]{Section~\ref*{sec:theory}} coincide with SignSGD there and the divergence results below need $\min(m_i,n_i)\ge2$ (\hyperref[app:vector]{Appendix~\ref*{app:vector}}).

Two assumptions run through the paper.

\begin{assumption}[\textbf{Lower boundedness}]\label{as:1}
$f(\mathbf{X}) \ge f^*$ for all $\mathbf{X}$; where explicitly invoked, each local $f_j \ge f_j^*$ as well.
\end{assumption}

\begin{assumption}[\textbf{Layer-wise smoothness}]\label{as:2}
For every layer $i$ and all $\mathbf{X},\mathbf{Y}$,
\begin{equation}
\|\nabla_i g(\mathbf{X}) - \nabla_i g(\mathbf{Y})\|_{*} \le L_i^{g}\,\|\mathbf{X}_i - \mathbf{Y}_i\|_{2\to 2}
\end{equation}
for $g = f$ (constant $L_i$) and $g = f_j$ (constant $L_{i,j}$). A weaker layer-wise $(L^0,L^1)$ form, stated in \hyperref[app:prelim]{Appendix~\ref*{app:prelim}}, is used only in \hyperref[cor:l0l1]{Corollary~\ref*{cor:l0l1}} (\hyperref[app:hyp]{Appendix~\ref*{app:hyp}}). Both forms are quantified over arbitrary $\mathbf{X},\mathbf{Y}$, as in \citet{gruntkowska2025error}; \hyperref[app:prelim]{Appendix~\ref*{app:prelim}} records what that strength implies.
\end{assumption}

A Linear Minimization Oracle (LMO) minimizes the first-order model of the objective over the unit ball of a norm, returning the steepest-descent direction in that norm (\hyperref[app:prelim]{Appendix~\ref*{app:prelim}}). Muon is the oracle of the spectral norm fixed above: writing $\mathbf{M}=\mathbf{U}\boldsymbol{\Sigma}\mathbf{V}^\top$ for the rank-$r$ singular value decomposition of $\mathbf{M}\in\mathbb{R}^{m\times n}$, the LMO direction is
\begin{equation}
A(\mathbf{M})=-\mathbf{U}\mathbf{V}^\top \in \mathbb{R}^{m\times n},
\qquad \operatorname{polar}(\mathbf{M}):=-A(\mathbf{M}),
\end{equation}
an orthogonalization, unscaled. Muon's reference implementation rescales it by $\sqrt{\max(1,m/n)}$ \citep{jordan2024muon} and the RMS$\to$RMS reading of the layer norm by $\sqrt{m/n}$: positive per-layer constants, invisible to $\operatorname{sign}(\cdot)$, which we carry in the step size~\eqref{eq:unit_gain_main} and not in the geometry (\hyperref[app:prelim]{Appendix~\ref*{app:prelim}}). In practice the oracle is applied to a momentum estimate rather than to a gradient; \eqref{eq:signa_general} below gives the full iteration. What Muon must transmit is that dense direction, precisely what a bandwidth-limited link cannot accommodate; our object is therefore its matrix-aware geometry at a one-bit budget.

\section{Theory}\label{sec:theory}
\paragraph{From SignSGD to Sign\,A.} SignSGD was introduced as a compressed SGD: a convergent method, of which only the sign of the update is transmitted. Nothing in that construction is specific to SGD, so one may apply it to any optimizer that emits a direction, expecting the compressed method to inherit what the direction contributed. LMO-based optimizers give grounds for that expectation, their theory being uniform in the geometry: one analysis covers every norm, which enters only through its oracle and a pair of norm-equivalence constants \citep{pethick2025training,kovalev2025understanding,riabinin2025gluon}. Were signing a direction harmless, it would be harmless across that family, and the member one seeks to compress is Muon.

The second ground for the expectation is where it fails. The sign step is itself an LMO, steepest descent in $\ell_\infty$ \citep{bernstein2024old}, so signing a Muon direction composes two oracles, each sound alone. The composition is an LMO for \emph{no} norm. An oracle for a norm $\|\cdot\|$ returns $\mathbf{D}=-A(\mathbf{G})$ with $\langle\mathbf{G},\mathbf{D}\rangle=\max_{\|\mathbf{Y}\|\le1}\langle\mathbf{G},\mathbf{Y}\rangle=\|\mathbf{G}\|_{\mathrm{dual}}$, nonnegative since $\mathbf{0}$ is feasible; \hyperref[th:1]{Theorem~\ref*{th:1}} provides a gradient on which the composition drives that inner product strictly negative. The guarantee is forfeited in the composition, not in either factor.

\paragraph{The Sign\,A family.} Let A be any optimizer built on LMO directions. A Sign\,A method initializes $\mathbf{M}_0 = \mathbf{0}$ and at each iteration $t \ge 1$ performs
\begin{equation}
\label{eq:signa_general}
\begin{aligned}
\mathbf{M}_t &= \mu \mathbf{M}_{t-1} + (1-\mu)\,\mathbf{G}_t,\\
\tilde{\mathbf{M}}_t &= \mathbf{M}_t
   \ \ \text{or}\ \ (1-\mu)\mathbf{G}_t + \mu \mathbf{M}_t\ \ \text{(Nesterov)},\\
\mathbf{D}_t &= -A(\tilde{\mathbf{M}}_t),
   \qquad \mathbf{s}_t = \operatorname{sign}(\mathbf{D}_t) \in \{\pm 1\}^{m \times n},\\
\mathbf{X}_{t} &= \mathbf{X}_{t-1} - \eta_t \mathbf{s}_t,
\end{aligned}
\end{equation}
where $\mathbf{X}_t \in \mathbb{R}^{m \times n}$ is the parameter matrix, $\mathbf{G}_t = \nabla f(\mathbf{X}_{t-1};\xi_t)$ the stochastic gradient at the current iterate, $\mu \in [0, 1)$ the momentum coefficient and $\eta_t$ the learning rate; $A(\cdot)$ is the LMO, and $\operatorname{sign}(\cdot)$ acts element-wise on the structured direction $\mathbf{D}_t$, with $\operatorname{sign}(0)$ resolved to an independent random $\pm1$ so that the transmitted alphabet stays binary. Momentum is in exponential-moving-average form throughout; the heavy-ball form differs by the positive factor $1/(1-\mu)$, which $\operatorname{sign}(\cdot)$ and the LMO both absorb, so no result below distinguishes them.

\paragraph{Three placements of the sign.} Of the three, only SignMuon is a Sign\,A method. All three keep the momentum of \eqref{eq:signa_general} and the update $\mathbf{X}_t=\mathbf{X}_{t-1}-\eta_t\,\mathbf{s}_t$, and differ only in where the sign acts on $\tilde{\mathbf{M}}_t$:
\begin{equation}
\label{eq:three_placements}
\mathbf{s}_t=
\begin{cases}
\operatorname{sign}\!\bigl(\operatorname{polar}(\tilde{\mathbf{M}}_t)\bigr) & \text{(SignMuon: after)},\\[2pt]
\operatorname{polar}\!\bigl(\operatorname{sign}(\tilde{\mathbf{M}}_t)\bigr) & \text{(MuonUSign: before)},\\[2pt]
\operatorname{sign}\!\bigl(\operatorname{polar}(\operatorname{sign}(\tilde{\mathbf{M}}_t))\bigr) & \text{(MuonSign: both)}.
\end{cases}
\end{equation}
The three are \hyperref[central_alg]{Algorithms~\ref*{central_alg}},~\ref{alg:muon_usign} and~\ref{alg:muon_sign} (\hyperref[app:alg]{Appendix~\ref*{app:alg}}), all transmitting one bit per matrix parameter on the uplink. On the downlink SignMuon and MuonSign distribute a $\pm1$-valued object and so cost one bit each way, whereas MuonUSign's server-side $\operatorname{polar}(\cdot)$ is dense and goes at full precision.

\subsection{Centralized Learning}
% In the centralized setting, SignMuon is applied to the matrix-valued parameters of a neural network. The practical implementation takes into account the tensor structure of model parameters. The optimization process is divided into two separate streams: matrix parameters (weights of convolutional and fully connected layers, $n_{dim} \ge 2$) are updated using the SignMuon LMO-based procedure, whereas vector parameters (bias, Batch Normalization layer, $n_{dim} = 1$) as well as the final classification layer are optimized using the standard AdamW optimizer. This design was originally introduced in the empirical description of Muon by \citet{jordan2024muon}. The first theoretical convergence result was derived by Li and Hong, who analyzed the smooth non-convex setting \cite{li2025note} and subsequently formalized within the Gluon framework \cite{riabinin2025gluon}.
Centrally, SignMuon is applied to the matrix-valued parameters ($n_{dim}\ge 2$), while vector parameters (biases, BatchNorm) and the final classification layer are trained with AdamW, the design of \citet{jordan2024muon}, whose LMO branch is what \citet{li2025note} and \citet{riabinin2025gluon} analyse. We approximate $\operatorname{polar}(\cdot)$ by a 5th-order Newton--Schulz iteration rather than a full SVD (\hyperref[alg:muon_lmo]{Algorithm~\ref*{alg:muon_lmo}}). One further implementation choice has consequences for every experiment below.

\paragraph{Per-layer step sizes: the unit-gain heuristic.} Write $\mathbf{P}_\ell$ for the matrix a method applies to layer $\ell$, of shape $m\times n$ (output $\times$ input dimension). It belongs to one of two families, each of \emph{exactly} known Frobenius norm: $\mathbf{P}_\ell=\mathbf{U}\mathbf{V}^\top$ with $\|\mathbf{P}_\ell\|_F=\sqrt{r}$, $r=\min(m,n)$, for the \textsc{lmo}-terminated methods, and $\mathbf{P}_\ell=\mathbf{s}_\ell\in\{\pm1\}^{m\times n}$ with $\|\mathbf{P}_\ell\|_F=\sqrt{mn}$ for the \textsc{sign}-terminated ones. The two scale differently with layer shape, so no single global $\eta$ is appropriate for both families, or across layers within one. We fix the shape dependence a priori and tune only a shape-free base rate $\eta_0$, giving layer $\ell$ the step size $\eta_\ell=\eta_0\lambda_\ell$ with
\begin{equation}
\label{eq:unit_gain_main}
\lambda_\ell=\frac{\sqrt{m}}{\|\mathbf{P}_\ell\|_F}
\quad\Longrightarrow\quad
\begin{cases}
\lambda_\ell^{\textsc{lmo}}=\sqrt{\max(1,m/n)},\\[2pt]
\lambda_\ell^{\textsc{sign}}=1/\sqrt{n}.
\end{cases}
\end{equation}
The criterion is that every layer's update have the same root-mean-square input--output gain, the average-case form of the spectral scaling condition \citep{yang2023spectral,large2024modular}, so that $\eta_0$ is the per-step RMS gain. Its \textsc{lmo} branch is not new: it reproduces the aspect-ratio factor that the reference Muon implementation already applies \citep{jordan2024muon}. That agreement is the external check we rely on, and it is what licenses applying the same criterion to the \textsc{sign} family, for which no such convention exists. We treat \eqref{eq:unit_gain_main} as a heuristic and apply it uniformly. \hyperref[app:lrscale]{Appendix~\ref*{app:lrscale}} gives the derivation and the measurement selecting the exponent $\tfrac12$ in $\lambda_\ell^{\textsc{sign}}=n^{-a}$ over $\mu$P's $1$ \citep{yang2021tuning}.

% Declared ahead of its first reference below: a two-column float can only
% move forward on the page, and from the next subsection it would land a
% page after the text that discusses it.
\begin{figure*}[!t]
    \centering
    \includegraphics[width=\textwidth]{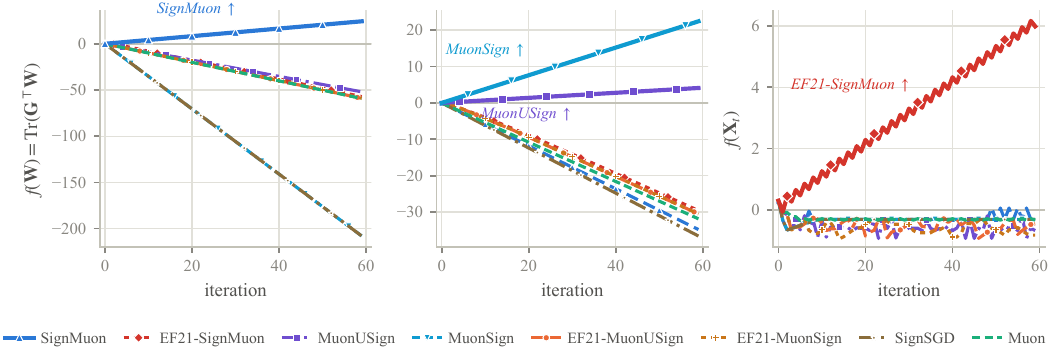}
    \caption{All eight methods on the three counterexample instances; the ascending method is drawn heavy and named in each panel. \emph{Left:} SignMuon, $4\times4$ (\hyperref[th:1]{Theorem~\ref*{th:1}}). \emph{Centre:} MuonUSign and MuonSign on their shared $5\times5$ instance (\hyperref[th:2]{Theorems~\ref*{th:2}}--\ref{th:3}); both panels plot~\eqref{eq:lin_obj}. \emph{Right:} EF21-SignMuon, $2\times2$ (\hyperref[th:ef_div]{Theorem~\ref*{th:ef_div}}), normalized units at $\eta=1$. Trajectories are momentum-free without loss (\hyperref[prop:reduction]{Proposition~\ref*{prop:reduction}}; \hyperref[lem:efsm_realize]{Lemma~\ref*{lem:efsm_realize}} for EF21-SignMuon).}
    \label{fig:divergence_plot}
\end{figure*}

\subsection{Divergence of the Three Sign Placements}
The descent property of the Muon LMO direction is lost under sign compression, before or after the oracle: in each of the three placements the compressed step can become an \emph{ascent} direction on a smooth objective. We refute the descent property on \emph{linear} objectives, the simplest smooth functions:
\begin{equation}
\label{eq:lin_obj}
f(\mathbf{X})=\langle \mathbf{G},\mathbf{X}\rangle=\operatorname{Tr}(\mathbf{G}^\top\mathbf{X}),\qquad \nabla f(\mathbf{X})\equiv \mathbf{G}.
\end{equation}
Gradient descent and full-precision Muon drive $f\to-\infty$ here, so a rule that instead drives $f\to+\infty$ is unambiguously ascending (\hyperref[rem:efsm_bdd]{Remark~\ref*{rem:efsm_bdd}} in \hyperref[app:proof_th4]{Appendix~\ref*{app:proof_th4}} restores \hyperref[as:1]{Assumption~\ref*{as:1}} without moving any trajectory below). On \eqref{eq:lin_obj} the step is the constant matrix $\mathbf{s}(\mathbf{G})$ whatever the momentum, so momentum affects neither convergence nor divergence (\hyperref[prop:reduction]{Proposition~\ref*{prop:reduction}}, \hyperref[app:proof_reduction]{Appendix~\ref*{app:proof_reduction}}) and
\begin{equation}
\label{eq:lin_increment}
f(\mathbf{X}_{t})-f(\mathbf{X}_{t-1})=-\eta_t\,\bigl\langle \mathbf{G},\ \mathbf{s}(\mathbf{G})\bigr\rangle
\end{equation}
for every $\mu\in[0,1)$ and either momentum rule. A single $\mathbf{G}$ with $\langle\mathbf{G},\mathbf{s}(\mathbf{G})\rangle<0$ therefore makes $f$ increase at every step, and such a $\mathbf{G}$ exists for each placement at small size.

\begin{summarythm}{\hyperref[th:1]{Theorems~\ref*{th:1}}--\ref{th:3} (summary)}
There exists $\mathbf{G}\in\mathbb{R}^{4\times4}$ with $\langle\mathbf{G},\mathbf{s}(\mathbf{G})\rangle<0$ for SignMuon, and $\mathbf{G}\in\mathbb{R}^{5\times5}$ with $\langle\mathbf{G},\mathbf{s}(\mathbf{G})\rangle<0$ for MuonUSign and for MuonSign simultaneously. On the corresponding objective~\eqref{eq:lin_obj} each method ascends at every iteration, for every $\eta_t>0$, every $\mu\in[0,1)$ and either momentum rule, and $f(\mathbf{X}_t)\to+\infty$ whenever $\sum_t\eta_t=\infty$.
\end{summarythm}

Both instances are explicit and both inner products exact rationals; \hyperref[app:proof_th1]{Appendices~\ref*{app:proof_th1}}--\ref{app:proof_th3} state the theorems, prove them, and bound the two shapes from below. \hyperref[fig:divergence_plot]{Figure~\ref*{fig:divergence_plot}} (left, centre) demonstrates the ascents. Negative results of this kind exist for uncompressed Muon as well: \citet{parshakova2026muondoesconvergeconvex} show that it need not converge on convex Lipschitz problems. Ours are due to the compressor rather than the geometry, and hold on a smooth objective.

\subsection{Centralized Error Feedback: EF21-SignMuon}

\hyperref[th:1]{Theorems~\ref*{th:1}}--\ref{th:3} rule out all three placements of the sign around the oracle. The classical remedy for a biased compressor is \emph{error feedback}: keep what the compressor discarded and fold it into the next message \citep{seide20141bit,karimireddy2019efsign}. We work throughout in its EF21 form \citep{richtarik2021ef21}, which stores an estimator and compresses the \emph{difference} to it, and so needs no bounded-gradient assumption. Its most direct use for SignMuon keeps the geometry, the sign \emph{after} the LMO, and applies error feedback to the oracle's output. The resulting method, \textbf{EF21-SignMuon} (\hyperref[ef21_signmuon]{Algorithm~\ref*{ef21_signmuon}}), tracks an estimate $\mathbf{d}_t^{\mathrm{est}}$ of the polar factor $\mathbf{D}_t=\operatorname{polar}(\tilde{\mathbf{M}}_t)$, updated by a scaled sign of the residual,
\begin{equation}\label{eq:efsm_est}
\begin{aligned}
\mathbf{d}_t^{\mathrm{est}}&=\mathbf{d}_{t-1}^{\mathrm{est}}+\alpha_t\operatorname{sign}\!\bigl(\mathbf{D}_t-\mathbf{d}_{t-1}^{\mathrm{est}}\bigr),\\
\alpha_t&=\operatorname{mean}\bigl|\mathbf{D}_t-\mathbf{d}_{t-1}^{\mathrm{est}}\bigr|,
\end{aligned}
\end{equation}
and steps $\mathbf{X}_t=\mathbf{X}_{t-1}-\eta_t\,\mathbf{d}_t^{\mathrm{est}}$.

Error feedback does not repair this placement: for \emph{every} step size and momentum setting there is a smooth objective on which EF21-SignMuon diverges.
% The three divergence theorems are now stated in the appendix, next to their
% proofs, but the paper's numbering is fixed by everything that cites it (the
% code, its READMEs, the reproduction guide): Theorems 1-3 are the three sign
% placements, 4 is this one, 5 is the convergence result. The counter is set by
% hand here and reset in the appendix to keep those numbers.
\setcounter{theorem}{3}
\begin{theorem}[Divergence of EF21-SignMuon]\label{th:ef_div}
For every $L>0$, step size $\eta>0$, momentum coefficient $\mu\in[0,1)$, and either momentum variant, there is an $L$-smooth (\hyperref[as:2]{Assumption~\ref*{as:2}}), bounded-below (\hyperref[as:1]{Assumption~\ref*{as:1}}) function $f:\mathbb{R}^{2\times2}\!\to\mathbb{R}$ on which EF21-SignMuon started at $\mathbf{X}_0=\mathbf{0}$ diverges: for an explicit constant $c=c(f)>0$,
\begin{equation}\label{eq:exact_rate}
f(\mathbf{X}_{t+2})-f(\mathbf{X}_t)=c\,L\eta^2>0\qquad(t\ge3),
\end{equation}
so $f(\mathbf{X}_t)\to+\infty$. In particular, no step-size rule $\eta=\eta(L,\mu)$ using only the smoothness and momentum constants can make the method convergent.
\end{theorem}
The construction (\hyperref[app:proof_th4]{Appendix~\ref*{app:proof_th4}}) turns the shared magnitude against the method. Its $2\times2$ LMO target has a large off-diagonal that reverses sign every step, holding the residual, and with it $\alpha_t$, at $\Theta(1)$, and a small constant on the diagonal, which that magnitude overshoots at every step; the diagonal estimate settles into a period-two cycle whose average has the \emph{wrong} sign, and the coordinate it drives diverges. The objective is not convex, unlike the linear ones of \hyperref[th:1]{Theorems~\ref*{th:1}}--\ref{th:3}: sustaining the cycle requires a target sequence that never settles, which here a bounded periodic field supplies. Dimension is not implicated, $2\times2$ being where the scaled sign's worst-case contraction $1/d$ is most favourable, nor is momentum: \hyperref[fig:ef21_momentum]{Figure~\ref*{fig:ef21_momentum}} (\hyperref[app:proof_th4]{Appendix~\ref*{app:proof_th4}}) records the same rate for every $\mu$ and both variants.

\subsection{Centralized Error Feedback: EF21-MuonUSign and EF21-MuonSign}
The shared magnitude $\alpha_t$ is not itself the defect: the two methods below couple all coordinates through the same scalar and converge. What fails is the target. Error feedback needs an estimator whose target varies with the step size, and the polar factor $\operatorname{polar}(\tilde{\mathbf{M}}_t)$ is not one: it is not Lipschitz in its argument, so it can move by $\Theta(1)$ between consecutive rounds however small $\eta$ is. Compressing the \emph{gradient} in its place restores that dependence, and the mechanism is then EF21-Muon~\cite{gruntkowska2025error}, the framework within which we apply sign compression.

That framework compresses the two directions of the link separately, through a pair $(\mathcal{C}^{\uparrow},\mathcal{C}^{\downarrow})$ of contractive compressors. In both directions we take the scaled sign $\mathcal{C}(\mathbf{Y}) = \operatorname{mean}|\mathbf{Y}|\,\operatorname{sign}(\mathbf{Y})$: one bit per parameter, plus one scalar per layer. \textbf{EF21-MuonUSign} takes the pair $(\mathcal{C},I)$, a compressed uplink and a full-precision model back, the appropriate allocation when only the uplink is constrained. It maintains a gradient estimator $\mathbf{g}_t^{\mathrm{est}}$, updated from the residual $\Delta_t = \tilde{\mathbf{M}}_t - \mathbf{g}_{t-1}^{\mathrm{est}}$:
\begin{equation}
\alpha_t = \operatorname{mean}(|\Delta_t|),\quad
\mathbf{g}_{t}^{\mathrm{est}} = \mathbf{g}_{t-1}^{\mathrm{est}} + \alpha_t\,\operatorname{sign}(\Delta_t),
\label{eq:ef21_central}
\end{equation}
and takes the step $\mathbf{X}_{t} = \mathbf{X}_{t-1} - \eta_t\mathbf{D}_t$ with $\mathbf{D}_t = -A(\mathbf{g}_{t}^{\mathrm{est}}) \approx \mathbf{U}_{t}\mathbf{V}^\top_{t}$, the polar factor of the estimator $\mathbf{g}_t^{\mathrm{est}}$ rather than of $\nabla f$. The sign in~\eqref{eq:ef21_central} acts on the internal compression residual alone, so the uplink stays at one bit per parameter (\hyperref[central_alg_ef]{Algorithm~\ref*{central_alg_ef}}).

\textbf{EF21-MuonSign} takes the pair $(\mathcal{C},\mathcal{C})$: one bit in each direction. Its downlink compressor is a second error-feedback loop, on the \emph{model} rather than the gradient, and it splits the method into two iterates:
\begin{equation}\label{eq:two_models}
\begin{aligned}
  \mathbf{X}_{t}&=\mathbf{X}_{t-1}-\eta_t\mathbf{D}_t &&\text{(server model)},\\
  \mathbf{W}_{t}&=\mathbf{W}_{t-1}+\alpha_t^{\downarrow}\operatorname{sign}(\mathbf{X}_{t}-\mathbf{W}_{t-1}) &&\text{(broadcast model)},
\end{aligned}
\end{equation}
with $\mathbf{D}_t=-A(\mathbf{g}_{t}^{\mathrm{est}})$ as before, and $\alpha_t^{\downarrow}=\operatorname{mean}|\mathbf{X}_t-\mathbf{W}_{t-1}|$. The server model $\mathbf{X}$ takes the exact LMO step and never leaves the server; what crosses the downlink is the one-bit update of $\mathbf{W}$, the only model the rest of the method observes, since gradients, momentum and the uplink residual are all computed there, whereas the guarantee bounds the gradient at $\mathbf{X}$. The rate is preserved under the second loop at a step size smaller by a factor of $\sqrt{r}$ in the layer rank, a penalty the analysis cannot avoid because the scaled sign contracts in the Euclidean norm and not in the layer norm the spectral geometry requires (\hyperref[rem:normmismatch]{Remark~\ref*{rem:normmismatch}}, \hyperref[app:hyp]{Appendix~\ref*{app:hyp}}); \hyperref[sec:exp_nanogpt]{Section~\ref*{sec:exp_nanogpt}} measures what it costs in practice.

Both methods descend where the placements they repair ascend: \hyperref[fig:divergence_plot]{Figure~\ref*{fig:divergence_plot}} (centre) has them on the $5\times5$ instance of \hyperref[th:2]{Theorems~\ref*{th:2}}--\ref{th:3}. \hyperref[app:hyp]{Appendix~\ref*{app:hyp}} shows them to be \emph{exact instances} of the EF21-Muon framework, the one substantive check being that the scaled sign is a contractive compressor (\hyperref[lem:signcontr]{Lemma~\ref*{lem:signcontr}}), so they inherit its guarantees (\hyperref[thm:conv]{Theorem~\ref*{thm:conv}}). Under \hyperref[as:1]{Assumptions~\ref*{as:1}}--\ref{as:2} both reach $\min_{t\le T}\mathbb{E}\|\nabla f(\mathbf{X}_t)\|_*^2=\mathcal{O}(T^{-1/2})$, the standard smooth nonconvex rate, which uncompressed Muon attains as well: one bit changes the constant and not the rate. Under $(L^0,L^1)$-smoothness EF21-MuonUSign reaches $\min_t\sum_i\mathbb{E}\|\nabla_i f(\mathbf{X}_t)\|_*=\mathcal{O}(T^{-1/4})$ at a constant step size (\hyperref[cor:smooth]{Corollaries~\ref*{cor:smooth}}--\ref{cor:l0l1}). Memory is the other cost, one model-sized buffer per compressed channel: EF21-MuonUSign holds the gradient estimator on each client, and EF21-MuonSign holds that estimator and, on the server, the broadcast model $\mathbf{W}$. Of these, only the client-side buffer is a practical constraint, since the server is the better-provisioned side.

\subsection{Federated Learning}
In the federated setting~\eqref{fed_eq} the placement decides where the oracle runs. With the sign \emph{after} it, each client must orthogonalize its own momentum: it computes a stochastic gradient at the broadcast model, maintains a momentum buffer, applies the Muon LMO (\hyperref[alg:muon_lmo]{Algorithm~\ref*{alg:muon_lmo}}), then uploads the elementwise sign of the result. The server takes a majority vote, $\mathbf{s}_t^{\mathrm{agg}}=\operatorname{sign}\bigl(\sum_{j} \mathbf{s}_t^{(j)}\bigr)$, and steps $\mathbf{X}_t=\mathbf{X}_{t-1}-\eta_t\mathbf{s}_t^{\mathrm{agg}}$, as SignSGD does \citep{bernstein2019signsgdmajority}. With the sign \emph{before}, the client uploads $\operatorname{sign}(\tilde{\mathbf{M}}_t^{(j)})$; the server votes, then applies a single LMO to the outcome. The vote rather than the average is what keeps the oracle's argument a sign matrix, so that the server-side method is exactly the MuonUSign of \eqref{eq:three_placements}; the error-feedback methods below instead average, as their framework prescribes. The direction that returns is dense, so MuonUSign broadcasts the model at full precision, whereas MuonSign signs that direction once more and broadcasts $\operatorname{sign}(\mathbf{D}_t)$. SignMuon and MuonSign therefore send a $\pm1$-valued object down as well as up, the vote itself in the first case, free of ties at an odd client count, so both directions cost one bit; and since every client applies the same update, copies started from a common $\mathbf{X}_0$ never drift.

Error feedback changes the uplink message, not where the oracle runs: each client sign-compresses the residual between its estimator and the quantity it would otherwise have sent, the polar factor for EF21-SignMuon and the momentum for the other two, and sends $(\mathbf{s}_t^{(j)},\alpha_t^{(j)})$; the server averages these into a global estimator. The extra per-layer scalar leaves the uplink at ${\approx}1$ bit per parameter, but that estimator is a scaled average of signs and so is dense, and the vote argument lapses: the full model must be broadcast unless a second error-feedback loop compresses it, as EF21-MuonSign does. Error feedback meets the same obstruction for SignSGD, where it is analysed for a single worker only \citep{karimireddy2019efsign} and its distributed forms compress the return path with a second loop of their own \citep{tang2019doublesqueeze}. The complete procedures are \hyperref[alg:fed_workerlmo]{Algorithms~\ref*{alg:fed_workerlmo}}--\ref{alg:fed_serverlmo}.

\section{Experiments}\label{sec:exp}

% Declared at the head of the section rather than beside its own
% subsection: a two-column float can only move forward, and from there it
% drifted onto the language-modelling page and shared it with two more.
\begin{figure*}[!t]
    \centering
    \includegraphics[width=\textwidth]{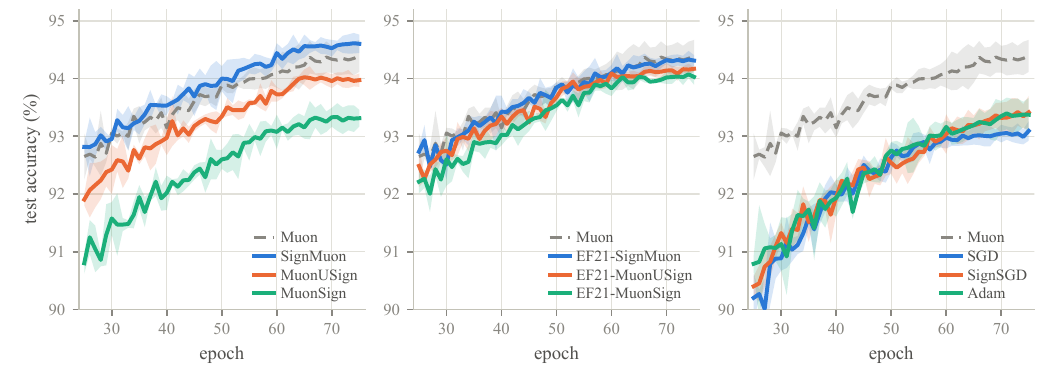}
    \caption{Centralized ResNet-18 on CIFAR-10: test accuracy from epoch $25$ at each method's selected $\eta_0$ (\hyperref[tab:cifar_main]{Table~\ref*{tab:cifar_main}}). Panels group the sign placements, the EF21 variants and the baselines; Muon is the gray dashed reference, bands $\pm1$ s.d.\ over three seeds.}
    \label{fig:cifar_results}
\end{figure*}

\paragraph{Descent in practice.} On a deterministic convex quadratic (\hyperref[app:images_task]{Appendix~\ref*{app:images_task}}) the first-order term of the descent lemma can be measured directly, through the alignment $\rho_t=\langle\nabla f,\mathbf{D}_t\rangle/(\|\nabla f\|_F\|\mathbf{D}_t\|_F)$ between the gradient and the direction actually taken. The direction descends when $\rho_t>0$, and \hyperref[th:1]{Theorems~\ref*{th:1}}--\ref{th:3} provide instances on which $\rho_t<0$; on random instances this does not occur, all three placements keeping $\rho_t$ bounded away from zero at every step. The counterexamples describe a worst case rather than a typical one, which is what permits the network results below to run contrary to them.

\subsection{Centralized Learning}
We compare all six sign-based methods against Muon, SignSGD, SGD and Adam on CIFAR-10 \cite{krizhevsky2009cifar} with a ResNet-18 adapted to low-resolution images, over $75$ epochs with $\eta_0$ cosine-annealed to zero. The only tuned hyperparameter is $\eta_0$, selected per method on a held-out split and read identically across methods through the unit-gain rule~\eqref{eq:unit_gain_main}. Numbers average three seeds; differences below the seed spread are not claimed (\hyperref[app:repro]{Appendix~\ref*{app:repro}}).

\begin{table}[!tb]
\centering
\small
\renewcommand{\arraystretch}{1.15}
\begin{tabular}{@{}lccc@{}}
\toprule
\textbf{Method} & $\boldsymbol{\eta_0}$ & \textbf{Test acc (\%)} & \textbf{Ep.\ to $\mathbf{90\%}$} \\
\midrule
SignMuon & $0.02$ & $94.60 \pm 0.15$ & $7.7$ \\
Muon & $0.1$ & $94.35 \pm 0.27$ & $7.7$ \\
EF21-SignMuon & $0.02$ & $94.31 \pm 0.11$ & $9.0$ \\
EF21-MuonUSign & $0.05$ & $94.14 \pm 0.07$ & $10.3$ \\
EF21-MuonSign & $0.005$ & $94.04 \pm 0.10$ & $11.0$ \\
MuonUSign & $0.02$ & $93.98 \pm 0.12$ & $10.3$ \\
Adam & $0.001$ & $93.37 \pm 0.27$ & $20.0$ \\
SignSGD & $0.002$ & $93.37 \pm 0.26$ & $19.0$ \\
MuonSign & $0.1$ & $93.31 \pm 0.21$ & $17.7$ \\
SGD & $0.02$ & $93.04 \pm 0.14$ & $21.3$ \\
\bottomrule
\end{tabular}
\caption{Centralized ResNet-18 / CIFAR-10, $75$ epochs. Test accuracy is the mean $\pm$ s.d.\ over three seeds of the last five epochs; ``Ep.\ to $90\%$'' is the mean epoch at which it first reaches $90\%$.}
\label{tab:cifar_main}
\end{table}

SignMuon ranks first, but by less than one seed spread over Muon: what \hyperref[tab:cifar_main]{Table~\ref*{tab:cifar_main}} supports is that it \emph{matches} Muon at one bit per parameter, not that it exceeds it. The separation that does resolve is the $1.2$ points over SignSGD, which spends the same budget without the \textsc{lmo}; the geometry accounts for it, not the compressor. Error feedback is not free here, EF21-SignMuon lying $0.29$ points below SignMuon against standard deviations of $0.11$ and $0.15$. Both families order the placements after, before, both sides, the last no better than SignSGD. The threshold column separates the methods more sharply than accuracy does: $7.7$ epochs to $90\%$ for SignMuon against $19.0$--$21.3$ for SGD, SignSGD and Adam.

\subsection{Federated Learning}
\begin{table}[!tb]
\centering
\footnotesize
\setlength{\tabcolsep}{3pt}
\renewcommand{\arraystretch}{1.1}
\begin{tabular}{@{}lccccc@{}}
\toprule
\textbf{Optimizer} & $\boldsymbol{\eta_0}$ & \textbf{Up} & \textbf{Down} &
\begin{tabular}[b]{@{}c@{}}\textbf{Rds.\ to}\\$\mathbf{80\%}$\end{tabular} &
\begin{tabular}[b]{@{}c@{}}\textbf{Test acc}\\\textbf{(\%)}\end{tabular} \\
\midrule
Muon & $0.1$ & $32$ & $32$ & $260$ & $85.98 \pm 0.26$ \\
SignMuon & $0.1$ & $\mathbf{1.09}$ & $\mathbf{1.09}$ & $540$ & $85.72 \pm 0.24$ \\
EF21-SignMuon & $0.02$ & $\mathbf{1.09}$ & $32$ & $640$ & $84.71 \pm 0.15$ \\
MuonUSign & $0.05$ & $\mathbf{1.09}$ & $32$ & $500$ & $84.56 \pm 0.35$ \\
EF21-MuonSign & $0.05$ & $\mathbf{1.09}$ & $\mathbf{1.09}$ & $980$ & $83.99 \pm 0.28$ \\
EF21-MuonUSign & $0.01$ & $\mathbf{1.09}$ & $32$ & $900$ & $83.56 \pm 0.41$ \\
MuonSign & $0.02$ & $\mathbf{1.09}$ & $\mathbf{1.09}$ & $860$ & $82.94 \pm 0.19$ \\
SGD & $0.05$ & $32$ & $32$ & $1240$ & $81.69 \pm 0.37$ \\
SignSGD & $0.01$ & $\mathbf{1.09}$ & $\mathbf{1.09}$ & $1140$ & $81.44 \pm 0.15$ \\
Adam & $0.001$ & $32$ & $32$ & --- & $77.12 \pm 1.04$ \\
\bottomrule
\end{tabular}
\caption{CIFAR-10 federated learning on CNN2, ordered by accuracy: $N=11$ clients, homogeneous split, $2000$ rounds at batch $192$, momentum $0.9$. Test accuracy is the mean of the final five evaluations over five seeds, $\pm$ one s.d.\ \emph{across} seeds; Up and Down are bits per parameter per round each way (\hyperref[app:commacct]{Appendix~\ref*{app:commacct}}). Every method but Adam reached $80\%$ on all five seeds, Adam on none.}
\label{tab:exp_3}
\end{table}

We compare the same ten methods on CNN2 (two convolutional blocks with BatchNorm and an MLP classifier), with CIFAR-10 split homogeneously across $N=11$ clients and no local steps; $N$ is odd, so the majority vote cannot tie. Rates are tuned per method at the reporting horizon on a split held out \emph{before} the client partition; two schemes cover all six, worker- and server-side LMO (\hyperref[alg:fed_workerlmo]{Algorithms~\ref*{alg:fed_workerlmo}}--\ref{alg:fed_serverlmo}, \hyperref[tab:fed_master]{Table~\ref*{tab:fed_master}}).

Federation separates the methods far more. At one bit per parameter in \emph{both} directions SignMuon reaches $85.72\%$ against Muon's $85.98\%$, about one seed standard deviation apart, and exceeds SignSGD by $4.3$ points, far more than centrally. Five seeds resolve the three placements: they span $2.8$ points in the order after, before, both sides, all three clearing SignSGD, the last by $1.5$ points where centrally it was level. The second error-feedback loop is not the source of the cost: EF21-MuonSign, compressed in both directions, stands marginally ahead of the uplink-only EF21-MuonUSign.

Less favourably for the theory, the two methods carrying an unconditional guarantee trail SignMuon by several seed spreads, as they do centrally, so what the guarantee costs here is charged to the placement. EF21-MuonSign is scored at its server model $\mathbf{X}_t$, not the model its clients hold~\eqref{eq:two_models}; the two track each other here, and \hyperref[sec:exp_nanogpt]{Section~\ref*{sec:exp_nanogpt}} shows that they need not (\hyperref[rem:dimension]{Remark~\ref*{rem:dimension}}, \hyperref[app:hyp]{Appendix~\ref*{app:hyp}}).

\subsection{Language Modelling}\label{sec:exp_nanogpt}

\begin{figure}[!tb]
\centering
\includegraphics[width=\columnwidth]{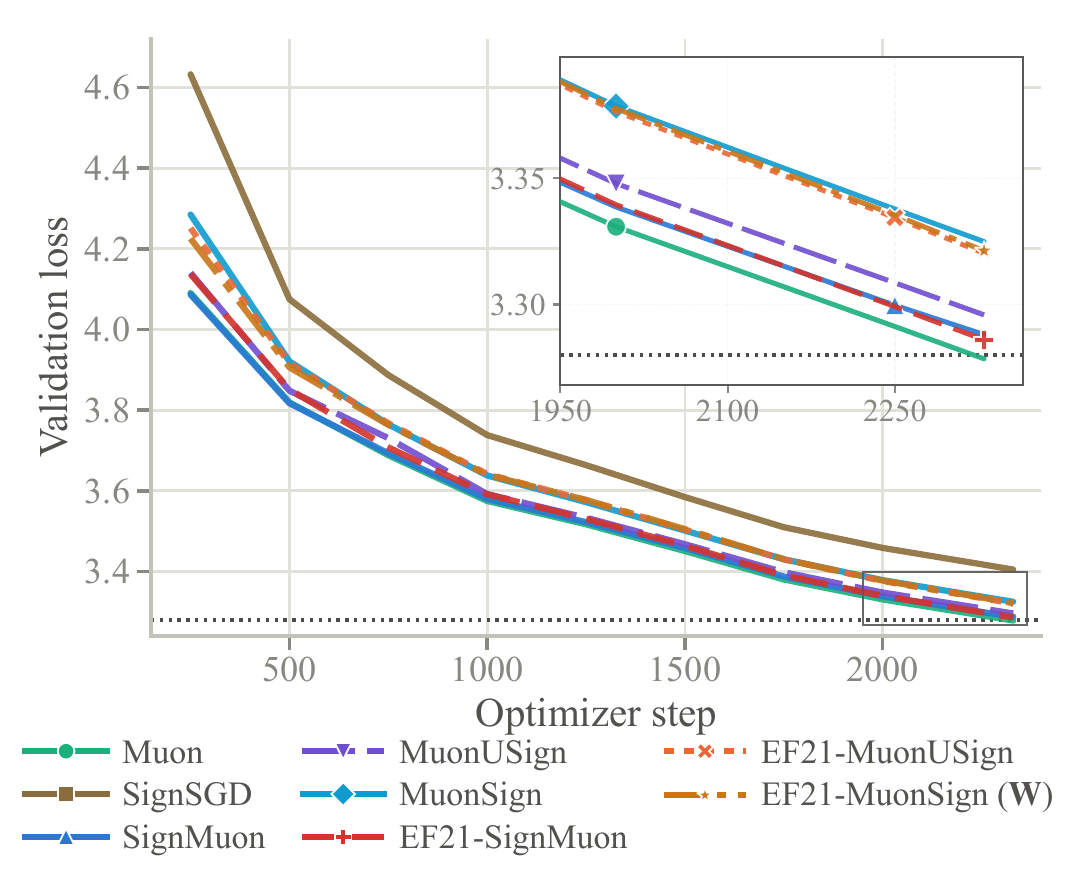}
\caption{NanoGPT speedrun, $8\times$H100: validation loss against optimizer step, one run per method; dotted is the target $3.28$, and EF21-MuonSign is drawn at its broadcast model $\mathbf{W}$. The inset magnifies the boxed tail, where SignSGD lies above the range shown.}
\label{fig:nanogpt}
\end{figure}

\begin{table}[!tb]
\centering\small
\setlength{\tabcolsep}{4pt}
\begin{tabular}{@{}llccc@{}}
\toprule
\textbf{Method} & \textbf{Step} & $\boldsymbol{\eta_0}$ & \textbf{Val.\ loss} & \textbf{Steps to $3.35$} \\
\midrule
Muon (record \#40)   & \textsc{lmo}  & 0.06 & \textbf{3.2785} & $1.00\times$ \\
EF21-SignMuon        & \textsc{lmo}  & 0.06 & 3.2860 & $1.02\times$ \\
SignMuon             & \textsc{sign} & 0.03 & 3.2881 & $1.02\times$ \\
MuonUSign            & \textsc{lmo}  & 0.06 & 3.2959 & $1.05\times$ \\
EF21-MuonUSign       & \textsc{lmo}  & 0.06 & 3.3203 & $1.13\times$ \\
EF21-MuonSign ($\mathbf{W}$) & \textsc{lmo} & 0.06 & 3.3213 & $1.14\times$ \\
MuonSign             & \textsc{sign} & 0.03 & 3.3249 & $1.14\times$ \\
SignSGD              & \textsc{sign} & 0.03 & 3.4049 & -- \\
\midrule
EF21-MuonSign ($\mathbf{X}$) & \textsc{lmo} & 0.06 & 5.5198 & -- \\
\bottomrule
\end{tabular}
\caption{NanoGPT after $2330$ steps ($611$M tokens), one run per method. Record \#40 reports $3.2780\pm0.0009$ over five seeds, so differences below ${\sim}0.003$ are noise; ``Steps to $3.35$'' is relative to Muon. The last row is the server model $\mathbf{X}$ of the same run as the $\mathbf{W}$ row.}
\label{tab:nanogpt}
\end{table}

We test where the matrices are large enough for the layer-rank dependence of \hyperref[cor:smooth]{Corollary~\ref*{cor:smooth}} to take effect: the modded-nanoGPT speedrun (record \#40), a $12$-layer transformer with $768\times3072$ hidden matrices, trained on $611$M FineWeb tokens on $8\times$H100 \citep{jordan2024moddednanogpt}. Each method replaces the record's Muon on the hidden matrices and gates, leaving all else untouched; our Muon run is record~\#40's own update rule, re-implemented without its Triton kernels. Rates are fixed a priori by the unit-gain rule at one $\eta_0$ per \emph{family}, so every contrast is matched-hyperparameter (\hyperref[app:nanogpt]{Appendix~\ref*{app:nanogpt}}).

Three conclusions follow (\hyperref[tab:nanogpt]{Table~\ref*{tab:nanogpt}}, \hyperref[fig:nanogpt]{Figure~\ref*{fig:nanogpt}}). First, composing the sign with the \textsc{lmo} transfers to language modelling: all six such methods improve on SignSGD by at least $0.08$ in validation loss, and the best two, EF21-SignMuon and SignMuon, lie within $0.01$ of full-precision Muon at equal wall-clock. Second, on the sign-after placement error feedback is free: EF21-SignMuon and SignMuon differ by $0.002$, below the ${\sim}0.003$ noise level, so the two are indistinguishable and jointly the strongest compressed methods here. The target both compress, $\operatorname{polar}(\tilde{\mathbf{M}}_t)$, is an orthogonal factor whose entries stay evenly spread in these runs, the regime in which the scaled sign loses least (\hyperref[app:nanogpt]{Appendix~\ref*{app:nanogpt}}). Third, EF21-MuonSign's two models~\eqref{eq:two_models} separate at this width. The run has a single client, so its broadcast model $\mathbf{W}$ is the only point at which a gradient is ever evaluated, and $\mathbf{W}$ is indistinguishable from EF21-MuonUSign: the compressed downlink costs nothing where training occurs. The server model $\mathbf{X}$, the iterate the guarantee bounds, settles $2.2$ nats above it, a persistent offset that \hyperref[app:nanogpt]{Appendix~\ref*{app:nanogpt}} localizes to one layer type, the zero-initialized output projection of each MLP block ($\mathtt{c\_proj}$). This is the mechanism of \hyperref[rem:dimension]{Remark~\ref*{rem:dimension}} (\hyperref[app:hyp]{Appendix~\ref*{app:hyp}}), not a tuning failure: only a spectrally contractive downlink compressor removes it.

The three settings therefore agree on the placement, and the placement they select is the one the theory excludes. Sign-after is what \hyperref[th:1]{Theorem~\ref*{th:1}} makes diverge and error feedback repairs for no $(L,\eta,\mu)$ (\hyperref[th:ef_div]{Theorem~\ref*{th:ef_div}}), yet in all three the strongest compressed method is a sign-after one, ahead of every variant carrying an unconditional guarantee.

\section{Conclusion}\label{sec:conclusion}

A sign step and a spectral step are each an LMO, sound alone; composed, they are an LMO for no norm, and the descent property each factor guarantees is lost in the composition. On small explicit instances, at every step size and momentum, every placement of the sign around the oracle can turn the update into an ascent direction on a linear objective: SignMuon after the LMO, MuonUSign before it, and MuonSign on both sides. Where error feedback is applied then decides whether it repairs this. On the oracle's \emph{output} it does not: the polar factor can move by a constant however small the step size, so one shared magnitude cannot track it, and for every $(L,\eta,\mu)$ some $L$-smooth objective makes EF21-SignMuon diverge. On the \emph{gradient} it does: EF21-MuonUSign attains the standard nonconvex rate at a one-bit uplink and EF21-MuonSign at one bit in each direction. Experiment ranks the methods in the opposite order, the strongest compressed method on all three architectures being a sign-after one, which carries no theoretical guarantee. We document that tension rather than resolve it: what the counterexamples preclude is a guarantee covering every problem in the class, not one holding under further conditions that ordinary problems satisfy. Identifying such conditions is the main open question, together with the convergence rate for the $(L^0,L^1)$ setting under a compressed downlink and whether \emph{any} one-bit compressor contracts in a layer norm.

%%%%%%%%%%%%%%%%%%%%%%%%%%%%%%%%%%%%%%%%%%%%%%%%%%%%%%%%%%%%
% References (placed after the Appendix; see note above).
% \bibliographystyle{unsrtnat}  % aaai2027.sty sets the style to aaai2027.bst automatically

% \bibliographystyle{unsrtnat}  % aaai2027.sty sets the style automatically
%%% <<< end signmuon_body.tex

%%% >>> begin arxiv/arxiv_backmatter.tex
%%%%%%%%%%%%%%%%%%%%%%%%%%%%%%%%%%%%%%%%%%%%%%%%%%%%%%%%%%%%%%%%%%%%%%%%%%%%
%  BACK MATTER FOR THE arXiv EDITION -- acknowledgments, author
%  contributions, code availability.
%
%  make_arxiv.py splices this in after the Conclusion and before the single
%  reference list. None of it exists in the anonymous AAAI sources, so all of
%  it is authored here.
%
%  The sections are unnumbered (\section*) so that the numbering of the body
%  is the same as in the AAAI PDF. \arxivrepo is defined in
%  arxiv_authors.tex.
%
%  The acknowledgment and the contribution split below are the real ones. If
%  funding or a compute provider has to be credited, the acknowledgment is
%  where it goes. make_arxiv.py prints a warning for every TODO left here.
%%%%%%%%%%%%%%%%%%%%%%%%%%%%%%%%%%%%%%%%%%%%%%%%%%%%%%%%%%%%%%%%%%%%%%%%%%%%

\section*{Acknowledgments}

We thank Alexander Tyurin for drawing our attention in June 2026 to the concurrent
SignMuon work of \citet{mishra2026signmuon}.

\section*{Author Contributions}

A.\,Kravatskiy proposed compressing Muon with a sign compressor in February
2026 and led the project. M.\,Smirnova ran the preliminary experiments and
wrote the first draft of the paper. The divergence counterexamples are due to
A.\,Kravatskiy; the reduction of EF21-MuonUSign and EF21-MuonSign to the
EF21-Muon framework of \citet{gruntkowska2025error} was carried out jointly.
Both authors ran the experiments reported here and expanded and edited the
manuscript into its present form.

\section*{Code Availability}

Code for this paper is available at \arxivrepo. The repository records the exact command and the
hyperparameters behind each reported number.
%%% <<< end arxiv/arxiv_backmatter.tex

%%% The single reference list.  It sits after the Conclusion, where
%%% AAAI puts it, and serves the appendix below as well: with both
%%% halves in one document there is nothing left to duplicate.
%%% aaai2027.sty selects aaai2027.bst on its own.
\bibliography{references}

\onecolumn
%%% >>> begin signmuon_appendix.tex
\appendix
\section{Appendix}\label{app}
% \subsection{Additional experiments / Proofs of Theorems}\label{app:exp}
\subsection{Preliminaries: LMO, Muon, and generalized smoothness}\label{app:prelim}

This subsection states in full the definitions that \hyperref[sec:proState]{Section~\ref*{sec:proState}} uses in abbreviated form.

Equip $\mathbb{R}^{m\times n}$ with the inner product $\langle \mathbf{M},\mathbf{D}\rangle=\sum_{i,j}{M}_{ij}{D}_{ij}$, let $\|\cdot\|$ be a norm on it, let $\mathcal{B}=\{\mathbf{D}\in\mathbb{R}^{m\times n}\,:\,\|\mathbf{D}\|\le1\}$ be its unit ball, and let $\|\mathbf{M}\|_{\mathrm{dual}}=\max_{\mathbf{D}\in\mathcal{B}}\langle\mathbf{M},\mathbf{D}\rangle$ be the dual norm. The \emph{Linear Minimization Oracle} (LMO) of $\|\cdot\|$ minimizes a linear form over the unit ball: it outputs
\begin{equation}\label{eq:lmo_def}
A(\mathbf{M})\in\arg\min_{\mathbf{D}\in\mathcal{B}}\ \langle \mathbf{M}, \mathbf{D}\rangle=-\bigl\{\mathbf{D}\in\mathcal{B}\,:\,\langle\mathbf{M},\mathbf{D}\rangle=\|\mathbf{M}\|_{\mathrm{dual}}\bigr\},
\end{equation}
Hence $\langle\mathbf{M},A(\mathbf{M})\rangle=-\|\mathbf{M}\|_{\mathrm{dual}}\le0$: the oracle returns a steepest-descent direction with respect to $\|\cdot\|$, normalized to the unit ball. The minimizer need not be unique, whence the inclusion; at rank-deficient $\mathbf{M}$ we use the rank-$r$ selection fixed in the next paragraph. An LMO method minimizes a differentiable objective $F$ by stepping along this direction: at the iterate $\mathbf{X}$ it forms an \emph{effective update direction} $\mathbf{M}$, a stochastic gradient or a momentum estimate in every method below, and moves along $A(\mathbf{M})$, the minimizer over $\mathcal{B}$ of the first-order model $\mathbf{D}\mapsto F(\mathbf{X})+\langle\mathbf{M},\mathbf{D}\rangle$. Each layer of \hyperref[sec:proState]{Section~\ref*{sec:proState}} carries its own norm and so its own oracle \citep{riabinin2025gluon}; for the role of the construction in optimizer design see \citep{bernstein2024old,bernstein2025deriving,pethick2025training,kovalev2025understanding,riabinin2025gluon,cesista2025schattenp,kravatskiy2025kyfannorms}.

Muon is the instance in which each matrix layer carries the spectral norm $\|\cdot\|_{2\to2}$, so that the oracle uses the matrix structure of the gradient and the update direction is obtained by orthogonalizing the gradient matrix \cite{bernstein2024old}. Let $\mathbf{M}=\mathbf{U}\boldsymbol{\Sigma}\mathbf{V}^\top$ be the singular value decomposition (SVD) of the matrix $\mathbf{M}\in\mathbb{R}^{m\times n}$, where $\mathbf{U}\in\mathbb{R}^{m\times r}$ and $\mathbf{V}\in\mathbb{R}^{n\times r}$ are the orthonormal matrices of singular vectors, $\boldsymbol{\Sigma}\in\mathbb{R}^{r\times r}$ is the diagonal matrix of singular values, and $r=\mathrm{rank}(\mathbf{M})$. Then the LMO direction is $A(\mathbf{M})=-\mathbf{U}\mathbf{V}^\top$: Muon selects an orthonormal update direction corresponding to the solution of the linear minimization problem induced by the spectral norm geometry. Orthogonalization equalizes the singular values of the update, so that the directions in which the gradient is weak are not crowded out by the dominant ones \citep{jordan2024muon}; large-scale pretraining studies report a corresponding efficiency gain over AdamW \citep{liu2025muon,shah2025practical}. The reference implementation obtains $\mathbf{U}_t\mathbf{V}_t^\top$ from a fixed number of Newton--Schulz iterations rather than from an SVD \citep{jordan2024muon}; \hyperref[alg:muon_lmo]{Algorithm~\ref*{alg:muon_lmo}} gives the procedure and the coefficients we use.

\paragraph{Choice of layer norm and the shape factor.} Muon is also presented in the RMS$\to$RMS operator norm, $\sqrt{n/m}\,\|\cdot\|_{2\to2}$ on $\mathbb{R}^{m\times n}$ \citep{bernstein2024old,pethick2025training}. The two readings give different oracles: the unit ball of $\sqrt{n/m}\,\|\cdot\|_{2\to2}$ is the spectral ball scaled by $\sqrt{m/n}$, so its LMO is $-\sqrt{m/n}\,\mathbf{U}\mathbf{V}^\top$. The unscaled $-\mathbf{U}\mathbf{V}^\top$ therefore belongs to the spectral norm and to no other, since a norm whose oracle it is has $\|\mathbf{M}\|_{\mathrm{dual}}=\langle\mathbf{M},\mathbf{U}\mathbf{V}^\top\rangle=\|\mathbf{M}\|_{*}$ for every $\mathbf{M}$ and hence equals $\|\cdot\|_{2\to2}$ by biduality. We take the spectral norm, in which \hyperref[as:2]{Assumption~\ref*{as:2}} and every statement below are stated, and treat the shape factor as a step-size question, for three reasons. It is a positive per-layer constant, so it changes neither $\operatorname{sign}(\cdot)$ nor the sign of any descent inner product, and no result of this paper is sensitive to it. The RMS$\to$RMS value $\sqrt{m/n}$ is not in fact the factor Muon uses: the reference implementation applies $\sqrt{\max(1,m/n)}$ \citep{jordan2024muon}, which agrees with it only for $m\ge n$ and is what our unit-gain rule returns (\hyperref[app:lrscale]{Appendix~\ref*{app:lrscale}}). And that rule must also scale the two sign-terminated placements, for which no norm supplies a scale at all: by \hyperref[th:1]{Theorems~\ref*{th:1}} and~\ref{th:3} they are oracles for none.

Finally, \hyperref[cor:l0l1]{Corollary~\ref*{cor:l0l1}} requires \hyperref[as:2]{Assumption~\ref*{as:2}} only in its weaker layer-wise $(L^0,L^1)$ form, which replaces the constant $L_i^{g}$ by $L^{0,g}_i + L^{1,g}_i\|\nabla_i g(\mathbf{X})\|_{*}$ \citep{riabinin2025gluon,pethick2025generalized}: for every layer $i$ and all parameter tuples $\mathbf{X},\mathbf{Y}$,
\begin{equation}
\|\nabla_i g(\mathbf{X}) - \nabla_i g(\mathbf{Y})\|_{*} \le \bigl(L^{0,g}_i + L^{1,g}_i\|\nabla_i g(\mathbf{X})\|_{*}\bigr)\,\|\mathbf{X}_i - \mathbf{Y}_i\|_{2\to 2},
\end{equation}
again with one pair of constants per layer for $g=f$ and per layer and client for $g=f_j$; at $L^{1,g}_i=0$ the display is \hyperref[as:2]{Assumption~\ref*{as:2}}. It is Assumptions~8--9 of \citet{gruntkowska2025error}, the hypotheses of the theorem \hyperref[cor:l0l1]{Corollary~\ref*{cor:l0l1}} quotes, stated there as here for \emph{arbitrary} pairs: the tuples may differ in every block although only block $i$ enters the bound. That quantifier is strong. Fixing $\mathbf{X}_i=\mathbf{Y}_i$ while varying the remaining blocks forces $\nabla_i g(\mathbf{X})=\nabla_i g(\mathbf{Y})$, so each $\nabla_i g$ depends on its own block alone, and a function satisfying the bound with finite constants is additively separable across layers. At $p=1$, the setting of every counterexample and synthetic measurement in this paper, the restriction is empty and the display is ordinary $(L^0,L^1)$-smoothness. For a multilayer network it is the framework's idealization, and we inherit it unweakened: the descent lemma behind \hyperref[thm:conv]{Theorem~\ref*{thm:conv}} applies the bound along a step in which every layer moves, which the restriction to pairs differing in one block would not license.

\subsection{Width-one blocks: vector parameters}\label{app:vector}

The setup of \hyperref[sec:proState]{Section~\ref*{sec:proState}} asks each block to be a matrix, and $\min(m_i,n_i)=1$ is permitted: a bias, a normalization gain, or any other one-dimensional parameter is the block $\mathbb{R}^{1\times n_i}$ or $\mathbb{R}^{m_i\times 1}$. Gluon states the same product space and leaves the norm on each block arbitrary, so it too admits them without comment \citep{riabinin2025gluon}; Scion is explicit, giving biases the RMS norm with oracle $\mathbf{b}/\|\mathbf{b}\|_{\mathrm{RMS}}$ \citep{pethick2025training}. Fixing the spectral norm on every block recovers nearly that: on a width-one block the rank is one, the spectral, nuclear and Euclidean norms coincide, and $\operatorname{polar}(\mathbf{g})=\mathbf{g}/\|\mathbf{g}\|_2$ for $\mathbf{g}\neq\mathbf{0}$. \hyperref[as:1]{Assumptions~\ref*{as:1}} and~\ref{as:2} then read unchanged, and the norm-equivalence constant of \hyperref[cor:smooth]{Corollary~\ref*{cor:smooth}} improves to $\bar\rho_i=\sqrt{r_i}=1$.

Substituting $\operatorname{polar}(\mathbf{g})=\mathbf{g}/\|\mathbf{g}\|_2$ into \eqref{eq:three_placements} gives, on $\mathbb{R}^{1\times n}$,
\begin{equation*}
\operatorname{sign}(\operatorname{polar}(\mathbf{g}))=\operatorname{sign}(\mathbf{g}),\qquad
\operatorname{polar}(\operatorname{sign}(\mathbf{g}))=\tfrac{1}{\sqrt{n}}\operatorname{sign}(\mathbf{g}),
\end{equation*}
and the two-sided placement returns $\operatorname{sign}(\mathbf{g})$ as well: all three are SignSGD up to a positive constant the step size absorbs. The ascent instances of \hyperref[th:1]{Theorems~\ref*{th:1}}--\ref{th:3} accordingly have $\min(m,n)\ge2$.

The unit-gain multipliers of \hyperref[app:lrscale]{Appendix~\ref*{app:lrscale}} need no special case either: $\lambda=\sqrt{m}/\|\mathbf{P}\|_F$ returns $(1,\,1/\sqrt{n})$ for the \textsc{lmo} and \textsc{sign} families on $\mathbb{R}^{1\times n}$, and $(\sqrt{m},\,1)$ on $\mathbb{R}^{m\times1}$.

\paragraph{Which parameters reach the methods.} Little of this affects our experiments, since the sign methods are applied where Muon is. Centralized and federated runs give the rule to parameters of two or more dimensions other than the classifier head, and route biases, BatchNorm scales and that head to AdamW as the auxiliary group, whose bandwidth cost \hyperref[app:commacct]{Appendix~\ref*{app:commacct}} counts. On nanoGPT we keep record~\#40's grouping unchanged: one-dimensional scalars, the embeddings and the head go to its distributed Adam, the hidden matrices and the two kinds of gate weight to the method under test. One gate is width-one, the $1\times12$ \texttt{smear\_gate}, and there $\operatorname{polar}$ is $\ell_2$ normalization, so the three sign-terminated methods reduce exactly to SignSGD on it and the \textsc{lmo} five to normalized momentum. In the CIFAR runs the corresponding parameters are one-dimensional, and our implementation returns those from the oracle unchanged rather than $\ell_2$-normalized, which differs from $\operatorname{polar}$ by a positive scale and so alters neither the sign nor the descent inner product.

\subsection{Divergence on linear objectives: the ascent criterion and momentum}\label{app:proof_reduction}

The criterion quoted in \hyperref[sec:theory]{Section~\ref*{sec:theory}} collapses each of the three methods to a single scalar inequality and removes momentum from the discussion entirely.

\begin{proposition}[Ascent criterion on linear objectives]\label{prop:reduction}
Run any of the three methods \eqref{eq:three_placements} on the linear objective \eqref{eq:lin_obj} with $\mathbf{G}\neq\mathbf{0}$, from an arbitrary $\mathbf{X}_0$, with any momentum coefficient $\mu\in[0,1)$ under either the Standard or the Nesterov rule. Then $\tilde{\mathbf{M}}_t=\gamma_t\,\mathbf{G}$ with $\gamma_t>0$, the update direction is the constant matrix $\mathbf{s}(\mathbf{G})$ obtained by substituting $\mathbf{G}$ for $\tilde{\mathbf{M}}_t$ in \eqref{eq:three_placements}, and
\begin{equation*}
f(\mathbf{X}_{t})-f(\mathbf{X}_{t-1})=-\eta_t\,\bigl\langle \mathbf{G},\ \mathbf{s}(\mathbf{G})\bigr\rangle .
\end{equation*}
In particular, if $\langle \mathbf{G},\mathbf{s}(\mathbf{G})\rangle<0$ then $f$ strictly increases at every iteration for any $\eta_t>0$, and $f(\mathbf{X}_t)\to+\infty$ whenever $\sum_t\eta_t=\infty$ (for instance, for any constant step size).
\end{proposition}

\paragraph{Proof of \hyperref[prop:reduction]{Proposition~\ref*{prop:reduction}}.} On the linear objective \eqref{eq:lin_obj} the gradient is globally constant, $\mathbf{G}_t=\mathbf{G}$, so the momentum buffer of \eqref{eq:signa_general} is $\mathbf{M}_t=(1-\mu)\sum_{i=0}^{t-1}\mu^i\,\mathbf{G}=(1-\mu^t)\,\mathbf{G}$. Both momentum rules then return a positive multiple of $\mathbf{G}$,
\begin{equation}
\tilde{\mathbf{M}}_t=\gamma_t\,\mathbf{G},\qquad
\gamma_t=
\begin{cases}
1-\mu^{t}, & \text{(Standard)},\\[2pt]
1-\mu^{t+1}, & \text{(Nesterov)},
\end{cases}
\end{equation}
both positive for $t\ge1$ because $\mu\in[0,1)$; the Nesterov case is $(1-\mu)\mathbf{G}+\mu(1-\mu^{t})\mathbf{G}$. (Under the heavy-ball convention every $\gamma_t$ is multiplied by $1/(1-\mu)$, which changes nothing below.) The elementwise $\operatorname{sign}(\cdot)$ and the Muon LMO $\operatorname{polar}(\cdot)$ are each invariant under multiplication by a positive scalar, so evaluating \eqref{eq:three_placements} at $\tilde{\mathbf{M}}_t=\gamma_t\mathbf{G}$ returns the same matrix as evaluating it at $\mathbf{G}$; that is, $\mathbf{s}_t=\mathbf{s}(\mathbf{G})$ for every $t$, independently of $\mu$ and of the momentum variant. Hence $f(\mathbf{X}_{t})-f(\mathbf{X}_{t-1})=\langle\mathbf{G},\mathbf{X}_{t}-\mathbf{X}_{t-1}\rangle=-\eta_t\langle\mathbf{G},\mathbf{s}(\mathbf{G})\rangle$, which is \eqref{eq:lin_increment}. If $\langle\mathbf{G},\mathbf{s}(\mathbf{G})\rangle<0$, then $f(\mathbf{X}_t)=f(\mathbf{X}_0)-\langle\mathbf{G},\mathbf{s}(\mathbf{G})\rangle\sum_{i=1}^{t}\eta_i$ increases strictly at every step and diverges to $+\infty$ whenever $\sum_t\eta_t=\infty$. $\blacksquare$

By \hyperref[prop:reduction]{Proposition~\ref*{prop:reduction}}, each divergence theorem reduces to exhibiting a single gradient $\mathbf{G}$ with $\langle\mathbf{G},\mathbf{s}(\mathbf{G})\rangle<0$; the three proofs below accomplish exactly this, and momentum requires no further comment.

\subsection{Proof of \hyperref[th:1]{Theorem~\ref*{th:1}} (Divergence of SignMuon)}\label{app:proof_th1}
% Back to 0: these three keep the numbers 1-3 they have always had, even though
% they are now stated after Theorem 4 (see the note at its statement).
\setcounter{theorem}{0}
\begin{theorem}[Divergence of SignMuon]\label{th:1}
There is a matrix $\mathbf{G}\in\mathbb{R}^{4\times4}$ with $\bigl\langle \mathbf{G},\operatorname{sign}(\operatorname{polar}(\mathbf{G}))\bigr\rangle=-\tfrac{42468}{103}<0$. Hence SignMuon ascends on \eqref{eq:lin_obj}: $f$ strictly increases at every iteration for all $\eta_t>0$, all $\mu\in[0,1)$, and both momentum variants.
\end{theorem}

For SignMuon $\mathbf{s}(\mathbf{G})=\operatorname{sign}(\operatorname{polar}(\mathbf{G}))$, so by \eqref{eq:lin_increment} it suffices to construct $\mathbf{G}\in\mathbb{R}^{4\times4}$ with $\langle \mathbf{G}, \operatorname{sign}(\operatorname{polar}(\mathbf{G})) \rangle < 0$. Group the SVD $\mathbf{G}=\mathbf{U}\boldsymbol{\Sigma}\mathbf{V}^\top$ of an invertible $\mathbf{G}$ as the polar decomposition $\mathbf{G}=\mathbf{Q}\mathbf{H}$, with $\mathbf{Q}=\mathbf{U}\mathbf{V}^\top=\operatorname{polar}(\mathbf{G})$ orthogonal and $\mathbf{H}=\mathbf{V}\boldsymbol{\Sigma}\mathbf{V}^\top\succ0$ symmetric. Only the symmetric part of $\mathbf{Q}^\top\operatorname{sign}(\mathbf{Q})$ then contributes to the trace against $\mathbf{H}$:
\begin{equation}\label{eq:qsignq}
\langle\mathbf{G},\operatorname{sign}(\operatorname{polar}(\mathbf{G}))\rangle=\operatorname{tr}\bigl(\mathbf{H}\mathbf{Q}^\top\operatorname{sign}(\mathbf{Q})\bigr)=\bigl\langle\mathbf{H},\operatorname{sym}(\mathbf{Q}^\top\operatorname{sign}(\mathbf{Q}))\bigr\rangle,
\end{equation}
so a counterexample with polar factor $\mathbf{Q}$ exists precisely when $\lambda_{\min}\bigl(\operatorname{sym}(\mathbf{Q}^\top\operatorname{sign}(\mathbf{Q}))\bigr)<0$: necessity because $\mathbf{H}\succ0$, sufficiency by taking $\mathbf{H}=\mathbf{w}\mathbf{w}^\top+\delta\mathbf{I}$ at a minimizing eigenvector $\mathbf{w}$ and $\delta>0$ small, which keeps $\mathbf{G}$ invertible and its polar factor unique. The question thus concerns orthogonal matrices alone.

We construct $\mathbf{G} \in \mathbb{R}^{4 \times 4}$ by defining a specific orthogonal matrix $\mathbf{O}$ and a specific rank-1 principal component $\mathbf{u}_1 \mathbf{v}_1^\top$. Let the orthogonal matrix $\mathbf{O}$ be given by the following exact rational numbers:
\begin{equation}
\mathbf{O} = \frac{1}{103} \begin{pmatrix} 101 & 20 & 2 & -2 \\ -20 & 97 & 20 & -20 \\ -2 & 20 & 2 & 101 \\ -2 & 20 & -101 & -2 \end{pmatrix}.
\end{equation}
Because no entry is zero, its element-wise sign matrix $\mathbf{S} = \operatorname{sign}(\mathbf{O})$ is uniquely defined. Now, let $\mathbf{u}_1$ and $\mathbf{v}_1$ be the following exact unit vectors:
\begin{equation}
\mathbf{u}_1 = \frac{1}{\sqrt{309}} \begin{pmatrix} 10 \\ -3 \\ 10 \\ 10 \end{pmatrix}, \quad \mathbf{v}_1 = \frac{1}{\sqrt{309}} \begin{pmatrix} 10 \\ 3 \\ -10 \\ 10 \end{pmatrix}.
\end{equation}
One can easily verify that $\mathbf{O} \mathbf{v}_1 = \mathbf{u}_1$, meaning $\mathbf{u}_1$ and $\mathbf{v}_1$ act perfectly as left and right singular vectors for this orthogonal space. The crucial feature of this geometry is that the Frobenius inner product between this rank-1 component and the sign matrix $\mathbf{S}$ yields an exact, strictly negative fraction:
\begin{equation}
\langle \mathbf{u}_1 \mathbf{v}_1^\top, \mathbf{S} \rangle = \mathbf{u}_1^\top \mathbf{S} \mathbf{v}_1 = -\frac{43}{103}.
\end{equation}

We construct the gradient matrix $\mathbf{G}$ by assigning a large singular value ($\sigma_1 = 1001$) to this pathological component and a singular value of $1$ to the rest of the orthogonal space. Exactly,
\begin{equation}\label{eq:th1_G}
\mathbf{G} := 1000\,\mathbf{u}_1\mathbf{v}_1^\top+\mathbf{O}
= \frac{1}{309}\begin{pmatrix}
100303 & 30060 & -99994 & 99994\\
-30060 & -8709 & 30060 & -30060\\
99994 & 30060 & -99994 & 100303\\
99994 & 30060 & -100303 & 99994
\end{pmatrix},
\end{equation}
or, to one decimal,
\begin{equation}
    \mathbf{G} \approx \begin{pmatrix}
        324.6 & 97.3 & -323.6 & 323.6 \\
        -97.3 & -28.2 & 97.3 & -97.3 \\
        323.6 & 97.3 & -323.6 & 324.6 \\
        323.6 & 97.3 & -324.6 & 323.6 \\
      \end{pmatrix}.
\end{equation}
Since $\mathbf{O}\mathbf{v}_1=\mathbf{u}_1$ and $\mathbf{v}_1^\top\mathbf{v}_1=1$,
\[
\mathbf{G}\mathbf{v}_1 = 1000\,\mathbf{u}_1(\mathbf{v}_1^\top\mathbf{v}_1)+\mathbf{O}\mathbf{v}_1 = 1001\,\mathbf{u}_1,
\]
so $(\mathbf{u}_1,\mathbf{v}_1)$ is a singular pair of $\mathbf{G}$ with $\sigma_1=1001$. For any $\mathbf{w}\perp\mathbf{v}_1$ we have
$\mathbf{G}\mathbf{w}=\mathbf{O}\mathbf{w}$; as $\mathbf{O}$ is orthogonal and $\mathbf{O}\mathbf{v}_1=\mathbf{u}_1$, the restriction
$\mathbf{O}|_{\mathbf{v}_1^\perp}$ is an isometry onto $\mathbf{u}_1^\perp$, hence $\sigma_2=\sigma_3=\sigma_4=1$. Moreover, with
$\mathbf{u}_k=\mathbf{O}\mathbf{v}_k$ for $k\ge2$,
\[
\textstyle\sum_{k\ge2}\mathbf{u}_k\mathbf{v}_k^\top=\mathbf{O}\bigl(\mathbf{I}-\mathbf{v}_1\mathbf{v}_1^\top\bigr)=\mathbf{O}-\mathbf{u}_1\mathbf{v}_1^\top,  
\]
so $\mathbf{U}_G\mathbf{V}_G^\top=\mathbf{u}_1\mathbf{v}_1^\top+(\mathbf{O}-\mathbf{u}_1\mathbf{v}_1^\top)=\mathbf{O}$. All four singular values are positive, so $\mathbf{G}$ is invertible and $\operatorname{polar}(\mathbf{G})=\mathbf{O}$ is its unique polar factor. The sign matrix \eqref{eq:muonsign_S} of \hyperref[th:2]{Theorems~\ref*{th:2}}--\ref{th:3} is full rank as well, so in all three proofs the argument of $\operatorname{polar}$ is invertible and its polar factor is unique: no statement depends on the selection rule fixed above for rank-deficient arguments.

Substituting the resulting matrix $\mathbf{G}$ into the descent condition yields:
\begin{equation}
\langle \mathbf{G}, \mathbf{S} \rangle = \langle 1000 \, \mathbf{u}_1 \mathbf{v}_1^\top + \mathbf{O}, \mathbf{S} \rangle = 1000 \langle \mathbf{u}_1 \mathbf{v}_1^\top, \mathbf{S} \rangle + \langle \mathbf{O}, \mathbf{S} \rangle.
\end{equation}
The inner product $\langle \mathbf{O}, \mathbf{S} \rangle$ is equivalent to the $L_1$ norm (sum of absolute values) of $\mathbf{O}$, which equals exactly $\frac{532}{103}$. Therefore:
\begin{equation}
\langle \mathbf{G}, \mathbf{S} \rangle = 1000 \left(-\frac{43}{103}\right) + \frac{532}{103} = -\frac{42468}{103} \approx -412.31.
\end{equation}
Here $\mathbf{S}=\operatorname{sign}(\mathbf{O})=\operatorname{sign}(\operatorname{polar}(\mathbf{G}))=\mathbf{s}(\mathbf{G})$, so $\langle \mathbf{G}, \mathbf{s}(\mathbf{G}) \rangle=-\tfrac{42468}{103}<0$. By \hyperref[prop:reduction]{Proposition~\ref*{prop:reduction}}, SignMuon strictly ascends, $f(\mathbf{X}_{t})-f(\mathbf{X}_{t-1})=\tfrac{42468}{103}\,\eta_t>0$ at every iteration, for every $\eta_t>0$, every $\mu\in[0,1)$, and both momentum variants; $f(\mathbf{X}_t)\to+\infty$ under any non-summable step size. $\blacksquare$

\subsection{Proof of \hyperref[th:2]{Theorem~\ref*{th:2}} (Divergence of MuonUSign)}\label{app:proof_th2}
\begin{theorem}[Divergence of MuonUSign]\label{th:2}
There is a matrix $\mathbf{G}\in\mathbb{R}^{5\times5}$ with $\bigl\langle \mathbf{G},\operatorname{polar}(\operatorname{sign}(\mathbf{G}))\bigr\rangle\approx-13.89<0$. Hence MuonUSign ascends on \eqref{eq:lin_obj} for all $\eta_t>0$, all $\mu\in[0,1)$, and both momentum variants.
\end{theorem}

MuonUSign (\cref{alg:muon_usign}) applies the sign \emph{before} the LMO, so $\mathbf{s}(\mathbf{G})=\operatorname{polar}(\operatorname{sign}(\mathbf{G}))$. By \eqref{eq:lin_increment} it suffices to construct $\mathbf{G}\in\mathbb{R}^{5\times5}$ with $\langle\mathbf{G},\operatorname{polar}(\operatorname{sign}(\mathbf{G}))\rangle<0$.

Here the step depends on $\mathbf{G}$ only through the magnitudes $|G_{ij}|$ and the sign pattern $\mathbf{S}=\operatorname{sign}(\mathbf{G})$; write $\mathbf{D}=\operatorname{polar}(\mathbf{S})$ for the direction it produces. Under the randomized convention $\mathbf{S}\in\{\pm1\}^{m\times n}$ throughout, and the entrywise identity $G_{ij}=|G_{ij}|\,S_{ij}$ holds without exception, both sides vanishing wherever $G_{ij}=0$. Consequently
\begin{equation}\label{eq:musign_pattern}
\bigl\langle\mathbf{G},\operatorname{polar}(\operatorname{sign}(\mathbf{G}))\bigr\rangle=\textstyle\sum_{i,j}|G_{ij}|\,S_{ij}D_{ij},
\end{equation}
and likewise $\langle\mathbf{G},\operatorname{sign}(\mathbf{D})\rangle=\sum_{i,j}|G_{ij}|\,S_{ij}\operatorname{sign}(D_{ij})$ for the MuonSign step of \hyperref[th:3]{Theorem~\ref*{th:3}}. Every summand in \eqref{eq:musign_pattern} is nonnegative unless some entry is \emph{mismatched}, $S_{ij}D_{ij}<0$. It therefore suffices to exhibit one sign matrix carrying a mismatched entry: inflating $|G_{ij}|$ there, with the other magnitudes held fixed, drives the sum below zero.

Fix the full-rank sign matrix $\mathbf{S}\in\{-1,1\}^{5\times 5}$,
\begin{equation}\label{eq:muonsign_S}
\mathbf{S} = \begin{pmatrix}
-1 & -1 &  1 &  1 &  1 \\
-1 & -1 &  1 & -1 & -1 \\
 1 & -1 &  1 &  1 & -1 \\
 1 &  1 & -1 & -1 &  1 \\
 1 &  1 &  1 & -1 &  1
\end{pmatrix},
\end{equation}
and let $\mathbf{D}=\operatorname{polar}(\mathbf{S})$ be its (unique, since $\mathbf{S}$ is full rank) polar factor. A direct computation gives $D_{4,2}=-1/\sqrt{17}\approx-0.2425<0$ while $\operatorname{sign}(D_{i,j})=S_{i,j}$ at all $24$ other entries; equivalently, $\operatorname{sign}(\mathbf{D})$ and $\mathbf{S}$ disagree at exactly the single entry $(4,2)$, where $S_{4,2}=+1$. We exploit this lone mismatch. Define
\begin{equation}\label{eq:muonsign_G}
\mathbf{G} = \epsilon\,\mathbf{S} + (M-\epsilon)\,\mathbf{e}_4\mathbf{e}_2^\top, \qquad \epsilon>0,\; M>\epsilon,
\end{equation}
so that $G_{4,2}=M>0$ and $G_{i,j}=\epsilon S_{i,j}$ otherwise; hence $\operatorname{sign}(\mathbf{G})=\mathbf{S}$ and $\operatorname{polar}(\operatorname{sign}(\mathbf{G}))=\mathbf{D}$ for every $M>0$. The descent inner product splits as
\begin{equation}
\langle \mathbf{G}, \mathbf{D}\rangle
= \epsilon\!\!\sum_{(i,j)\neq(4,2)}\!\! S_{i,j}D_{i,j} \;+\; M\,D_{4,2}
= \epsilon\,C + M\,D_{4,2},
\end{equation}
where $C:=\sum_{(i,j)\neq(4,2)}|D_{i,j}|=10.366>0$ is fixed (every such entry agrees in sign with $\mathbf{S}$), while $M\,D_{4,2}=-M/\sqrt{17}\to-\infty$. Thus $\langle\mathbf{G},\mathbf{D}\rangle<0$ for any $M>\sqrt{17}\,C\epsilon\approx42.7\,\epsilon$; with $\epsilon=1,\ M=100$ one obtains $\langle\mathbf{G},\mathbf{D}\rangle=-13.89$, for the exact polar factor that the theorem is stated over. By \hyperref[prop:reduction]{Proposition~\ref*{prop:reduction}}, MuonUSign strictly ascends on $f(\mathbf{X})=\langle\mathbf{G},\mathbf{X}\rangle$ for every $\eta_t>0$, every $\mu\in[0,1)$, and both momentum variants; $f(\mathbf{X}_t)\to+\infty$ under any non-summable step size. $\blacksquare$

\subsection{Proof of \hyperref[th:3]{Theorem~\ref*{th:3}} (Divergence of MuonSign)}\label{app:proof_th3}
\begin{theorem}[Divergence of MuonSign]\label{th:3}
For the \emph{same} matrix $\mathbf{G}\in\mathbb{R}^{5\times5}$ as in \hyperref[th:2]{Theorem~\ref*{th:2}}, $\bigl\langle \mathbf{G},\operatorname{sign}(\operatorname{polar}(\operatorname{sign}(\mathbf{G})))\bigr\rangle=-76<0$. Hence MuonSign ascends on \eqref{eq:lin_obj} for all $\eta_t>0$, all $\mu\in[0,1)$, and both momentum variants.
\end{theorem}

MuonSign (\cref{alg:muon_sign}) signs the polar factor as well, so $\mathbf{s}(\mathbf{G})=\operatorname{sign}(\operatorname{polar}(\operatorname{sign}(\mathbf{G})))=\operatorname{sign}(\mathbf{D})$. We reuse the \emph{same} $\mathbf{S}$ and $\mathbf{G}$ of \eqref{eq:muonsign_S}--\eqref{eq:muonsign_G}: since $\operatorname{sign}(\mathbf{G})=\mathbf{S}$, the bidirectional step is the constant matrix $\operatorname{sign}(\mathbf{D})$, which agrees with $\mathbf{S}$ at all $24$ entries except $(4,2)$, where $\operatorname{sign}(D_{4,2})=-1=-S_{4,2}$. Using $S_{i,j}^2=1$ everywhere,
\begin{equation}
\begin{aligned}
\langle \mathbf{G}, \operatorname{sign}(\mathbf{D})\rangle
&= \epsilon\!\!\sum_{(i,j)\neq(4,2)}\!\! S_{i,j}\operatorname{sign}(D_{i,j}) + M\,\operatorname{sign}(D_{4,2})\\
&= 24\,\epsilon - M .
\end{aligned}
\end{equation}
This is negative for every $M>24\epsilon$; with $\epsilon=1,\ M=100$ it equals exactly $-76$. By \hyperref[prop:reduction]{Proposition~\ref*{prop:reduction}}, MuonSign strictly ascends on $f(\mathbf{X})=\langle\mathbf{G},\mathbf{X}\rangle$ for every $\eta_t>0$, every $\mu\in[0,1)$, and both momentum variants; $f(\mathbf{X}_t)\to+\infty$ under any non-summable step size. In particular, the \emph{same} $5\times5$ linear instance is an ascent instance for the uplink-only placement (MuonUSign) and for the bidirectional one (MuonSign) alike. $\blacksquare$

%%% >>> begin signmuon_mishra_comparison.tex (was \input on signmuon_appendix.tex:183)
% =====================================================================
%  Comparison with the concurrent convergence claim of Mishra, Trivedi
%  and Kumar (2026) for SignMuon.  Drop-in appendix subsection for
%  the one-column appendix builds.  Self-contained: refers only to this
%  paper's Theorem~\ref{th:1}, Proposition~\ref{prop:reduction}, the
%  linear objective~\eqref{eq:lin_obj}, Remark~\ref{rem:efsm_bdd}, and
%  \citet{mishra2026signmuon}.
% =====================================================================
\subsection{Comparison with the convergence claim for SignMuon}\label{app:mishra}

Concurrently, \citet{mishra2026signmuon} introduced the sign-after-LMO method
under the same name \emph{SignMuon}: their algorithm forms the momentum
$\mathbf{M}_t$, computes its polar factor, and steps along
$\mathbf{S}_t=\operatorname{sign}(\operatorname{polar}(\mathbf{M}_t))$,
normalized to $\mathbf{D}_t=\mathbf{S}_t/\sqrt{mn}$; after absorbing
$1/\sqrt{mn}$ into $\eta$ this is the SignMuon step of
\hyperref[sec:theory]{Section~\ref*{sec:theory}}. Their abstract and contributions attribute to this
method an $\mathcal{O}(1/\sqrt{T})$ stationarity guarantee; the theorem that
establishes the rate is stated, accurately, for a \emph{gradient-sign
instantiation}. The distance between the two statements is the subject of this
subsection. Two facts resolve it, and neither contradicts
\hyperref[th:1]{Theorem~\ref*{th:1}}: the finalized rate is proved for an update that computes
no polar factor and carries no momentum, and the one bound of theirs that does
apply to the SignMuon update is an inequality whose right-hand side exceeds
its left-hand side on the instance of \hyperref[th:1]{Theorem~\ref*{th:1}}, at every iteration
and for every step size, so that it is satisfied there while the method
ascends.

\paragraph{The two components of their analysis.}
Their stationarity measure is
$\mathcal{G}_T=\tfrac1T\sum_{t=0}^{T-1}\mathbb{E}\bigl[\|\nabla
f(\mathbf{X}_t)\|_1/\sqrt{mn}\bigr]$, controlled through the per-entry
sign-error probabilities
$q_{ij,t}=\Pr\bigl(\mathbf{S}_{t,ij}\neq\operatorname{sign}([\nabla
f(\mathbf{X}_t)]_{ij})\mid\mathbf{X}_t\bigr)$ of the transmitted sign matrix
$\mathbf{S}_t$. The first component is generic. For an \emph{arbitrary} sign
oracle, the conditional identity
\begin{equation}\label{eq:mishra_identity}
\mathbb{E}\bigl[\langle \nabla f(\mathbf{X}_t),\,\mathbf{D}_t\rangle \,\big|\, \mathbf{X}_t\bigr]
=\frac{1}{\sqrt{mn}}\sum_{i,j}\bigl|[\nabla f(\mathbf{X}_t)]_{ij}\bigr|\,\bigl(1-2q_{ij,t}\bigr)
\end{equation}
and the descent lemma of spectral smoothness telescope into
\begin{equation}\label{eq:mishra_generic}
\mathcal{G}_T
\;\le\;
\frac{f(\mathbf{X}_0)-f^\ast}{\eta T}
+\frac{L_\ast}{2}\,\eta
+\underbrace{\frac{2}{T}\sum_{t=0}^{T-1}\mathbb{E}\!\left[\frac{1}{\sqrt{mn}}\sum_{i,j}\bigl|[\nabla f(\mathbf{X}_t)]_{ij}\bigr|\,q_{ij,t}\right]}_{=:\,R_T}.
\end{equation}
The second component estimates the residual $R_T$, and it is here that the
oracle is fixed: for $\mathbf{S}_t=\operatorname{sign}(\widetilde{\mathbf{G}}_t)$,
where $\widetilde{\mathbf{G}}_t$ is an unbiased stochastic gradient with
coordinatewise variance at most $\sigma_{ij}^2/n_b$ at mini-batch size $n_b$,
a Markov--Jensen argument yields $\bigl|[\nabla
f(\mathbf{X}_t)]_{ij}\bigr|\,q_{ij,t}\le\sigma_{ij}/\sqrt{n_b}$, hence
$R_T\le2\|\sigma\|_1/\sqrt{mn\,n_b}$ with
$\|\sigma\|_1=\sum_{i,j}\sigma_{ij}$. The residual vanishes as
$n_b\to\infty$, and the choice $n_b=T$ produces the
$\mathcal{O}(1/\sqrt{T})$ rate.

\paragraph{The finalized rate is not a result about SignMuon.}
The oracle of the second component transmits the sign of the stochastic
gradient itself. That update invokes neither the momentum buffer nor the polar
factor; as an algorithm it is SignSGD at batch size $n_b$, in single-worker
and majority-vote form, normalized by $1/\sqrt{mn}$ and analysed under
spectral rather than coordinatewise smoothness, which their own comparison
identifies as the sole improvement over \citet{bernstein2018signsgd}. Nor is
the restriction incidental. The Markov--Jensen estimate bounds the probability
of a sign error by the ratio of noise to signal,
$\sigma_{ij}/(\sqrt{n_b}\,|[\nabla f(\mathbf{X}_t)]_{ij}|)$, and is therefore
available exactly when sampling noise is the only mechanism by which a
transmitted sign can disagree with the gradient's. For the SignMuon oracle the
disagreement is structural rather than stochastic: even with exact gradients
($\sigma\equiv0$), $q_{ij,t}$ is the indicator that
$\operatorname{sign}(\operatorname{polar}(\mathbf{M}_t))$ and
$\operatorname{sign}(\nabla f(\mathbf{X}_t))$ differ at $(i,j)$, a quantity
that no batch size reduces. The $\mathcal{O}(1/\sqrt{T})$ rate accordingly
attaches to the gradient-sign update, and to SignMuon their analysis offers
only \eqref{eq:mishra_generic} with $R_T$ unestimated.

\paragraph{Scope of the generic bound.}
Inequality \eqref{eq:mishra_generic} is valid for every sign oracle, and for
that reason asserts nothing until $R_T$ is estimated. By
\eqref{eq:mishra_identity}, the summand of $R_T$ at time $t$ exceeds the
corresponding summand of $\mathcal{G}_T$ precisely when
$\mathbb{E}[\langle\nabla f(\mathbf{X}_t),\mathbf{D}_t\rangle\mid\mathbf{X}_t]\le0$,
that is, precisely when the expected step fails to be a descent direction.
Whenever this occurs at every $t$, the right-hand side of
\eqref{eq:mishra_generic} exceeds the left-hand side termwise and the
inequality holds irrespective of how the iterates behave. The bound therefore
has content only where the transmitted sign is already positively aligned with
the gradient in expectation; that alignment is the property a convergence
proof for SignMuon would have to establish, and it is the property
\hyperref[th:1]{Theorem~\ref*{th:1}} refutes.

\paragraph{On the instance of \hyperref[th:1]{Theorem~\ref*{th:1}}.}
Consider the linear objective~\eqref{eq:lin_obj} with the $4\times4$ gradient
$\mathbf{G}$ of \hyperref[th:1]{Theorem~\ref*{th:1}}. By \hyperref[prop:reduction]{Proposition~\ref*{prop:reduction}} the
gradient equals $\mathbf{G}$ and the update direction equals the constant
matrix $\operatorname{sign}(\operatorname{polar}(\mathbf{G}))$ at every
iteration, whatever the momentum, so the oracle is deterministic and $q_{ij}$
is the indicator that
$\operatorname{sign}(\operatorname{polar}(\mathbf{G}))_{ij}\neq
\operatorname{sign}(\mathbf{G}_{ij})$. Identity \eqref{eq:mishra_identity}
then evaluates exactly:
\begin{equation}\label{eq:mishra_instance}
\sum_{i,j}|\mathbf{G}_{ij}|\,(1-2q_{ij})
=\bigl\langle \mathbf{G},\,\operatorname{sign}(\operatorname{polar}(\mathbf{G}))\bigr\rangle
=-\tfrac{42468}{103}\;<\;0 .
\end{equation}
Dividing by $\sqrt{mn}=4$, every term of $R_T$ exceeds the corresponding term
of $\mathcal{G}_T$ by the same amount, so that
\begin{equation}\label{eq:mishra_dominates}
R_T=\mathcal{G}_T+\tfrac{10617}{103}
\qquad\text{for every } T .
\end{equation}
The right-hand side of \eqref{eq:mishra_generic} therefore exceeds the
left-hand side by at least $10617/103$ for every $\eta>0$, every $L_\ast$ and
every $T$: the inequality is satisfied and constrains nothing. What the
trajectory actually does is read off the same identity,
$f(\mathbf{X}_t)-f(\mathbf{X}_{t-1})=-\eta\langle\mathbf{G},\mathbf{D}_t\rangle
=\eta\cdot\tfrac{10617}{103}>0$, the divergence of \hyperref[th:1]{Theorem~\ref*{th:1}}: the
excess $R_T-\mathcal{G}_T$ of \eqref{eq:mishra_dominates} and the
per-iteration ascent rate
$\eta^{-1}\bigl(f(\mathbf{X}_t)-f(\mathbf{X}_{t-1})\bigr)$ are the same
number.
The one assumption of theirs the linear instance lacks is lower boundedness,
and the modification of \hyperref[rem:efsm_bdd]{Remark~\ref*{rem:efsm_bdd}} applies unchanged: $f$ is
unbounded below only on a half-space the iterates never enter, where a smooth
bounded replacement restores the assumption without moving the trajectory or
any quantity above.

In their terms, \hyperref[th:1]{Theorem~\ref*{th:1}} exhibits a smooth instance on
which the sign-error probabilities $q_{ij,t}$ of the SignMuon oracle, averaged
over the entries with weights $\bigl|[\nabla f(\mathbf{X}_t)]_{ij}\bigr|$,
exceed $\tfrac12$ at every iteration. Any rate extracted from
\eqref{eq:mishra_generic} requires that weighted average to stay below
$\tfrac12$ by a uniform margin, and no assumption of theirs implies such a bound for
$\operatorname{sign}(\operatorname{polar}(\mathbf{M}_t))$. Their theorems
stand as guarantees for majority-vote SignSGD under spectral smoothness. A
convergence guarantee for SignMuon they are not, and \hyperref[th:1]{Theorem~\ref*{th:1}} shows
that none is available at this level of generality.
%%% <<< end signmuon_mishra_comparison.tex

\subsection{Extended comparison with S-Muon}\label{app:smuon}
% Moved here from Related Work to save main-part space.
Taking the norms dual to the Ky Fan $k$-norms, \citet{kravatskiy2025kyfannorms} obtain the \emph{Fanion} family, whose updates $\sum_{i\le k}\mathbf{u}_i\mathbf{v}_i^\top$ interpolate between the rank-one step of the nuclear norm and Muon's full-rank $\mathbf{U}\mathbf{V}^\top$ at $k=\min(m,n)$; a conic combination of LMO algorithms is again an LMO algorithm, for the norm dual to the corresponding combination of dual norms. Their S-Muon is one such combination, $\tau\,\mathbf{U}\mathbf{V}^\top+(1-\tau)\,c\,\operatorname{sign}(\mathbf{M}_t)$ with fixed $\tau\in[0,1]$, $c>0$ (their notation differs; we reserve $\alpha$ for compressor contraction and $\eta$ for the learning rate). There the sign enters \emph{inside} the oracle, so the step is still an LMO for an explicit norm and inherits the convergence theory of one; our three placements act \emph{around} the oracle. A caution from the same work applies to us directly: they exhibit an LMO method (rank-one Neon) markedly worse than Muon in practice despite sharing its convergence asymptotics in the bounds of \citet{kovalev2025understanding} and \citet{riabinin2025gluon}, from which our own guarantee descends. A rate of the form $\mathcal{O}(T^{-1/2})$ is not a prediction of the performance a method will attain.

%%% >>> begin ef21_signmuon_divergence.tex (was \input on signmuon_appendix.tex:189)
% =====================================================================
%  Divergence of EF21-SignMuon (Theorem th:ef_div).
%  Drop-in appendix subsection for v2_SignMuon_AAAI.tex (double-column).
%  Structure: Theorem -> Proof idea (narrative) -> Proof in three labelled
%  parts (rescaling; the limit cycle; the realization).  Exact rational
%  arithmetic; the construction and its all-method comparison are implemented
%  in code/counterexamples/problems.py (figure via run_counterexamples.py).
% =====================================================================
\providecommand{\smat}[4]{\bigl(\begin{smallmatrix}#1&#2\\#3&#4\end{smallmatrix}\bigr)}

\subsection{Proof of \hyperref[th:ef_div]{Theorem~\ref*{th:ef_div}} (Divergence of EF21-SignMuon)}\label{app:proof_th4}

% The statement is repeated here for the reader. The counter is set by hand so
% that the restatement carries the main text's number 4 and does not consume a
% new one; ef21_musign_reduction.tex sets it again for Theorem 5. No \label is
% attached: th:ef_div stays the main-text statement.
\setcounter{theorem}{3}
\begin{theorem}[Divergence of EF21-SignMuon]
For every $L>0$, step size $\eta>0$, momentum coefficient $\mu\in[0,1)$, and
either momentum variant, there is an $L$-smooth (\hyperref[as:2]{Assumption~\ref*{as:2}}),
bounded-below (\hyperref[as:1]{Assumption~\ref*{as:1}}) function
$f:\mathbb{R}^{2\times2}\!\to\mathbb{R}$ on which EF21-SignMuon started at
$\mathbf{X}_0=\mathbf{0}$ diverges: for an explicit constant $c=c(f)>0$,
\begin{equation}\label{eq:exact_rate_app}
f(\mathbf{X}_{t+2})-f(\mathbf{X}_t)=c\,L\eta^2>0\qquad(t\ge3),
\end{equation}
so $f(\mathbf{X}_t)\to+\infty$. In particular, no step-size rule
$\eta=\eta(L,\mu)$ using only the smoothness and momentum constants can make
the method convergent.
\end{theorem}

Recall from the main text that EF21-SignMuon (\hyperref[ef21_signmuon]{Algorithm~\ref*{ef21_signmuon}})
does not step along the polar factor
$\mathbf D_t=\operatorname{polar}(\tilde{\mathbf M}_t)$ itself, but along the
error-feedback estimate $\mathbf d_t^{\mathrm{est}}$ of \eqref{eq:efsm_est},
followed by $\mathbf X_t=\mathbf X_{t-1}-\eta\,\mathbf d_t^{\mathrm{est}}$. A
\emph{single} magnitude $\alpha_t$ rescales the signs of \emph{all} entries at
once, and it is this coupling that the counterexample exploits. We now prove
\hyperref[th:ef_div]{Theorem~\ref*{th:ef_div}}.

\subsubsection*{Proof idea}

\emph{The mechanism.} The error-feedback update \eqref{eq:efsm_est} moves
\emph{every} entry of $\mathbf d_t^{\mathrm{est}}$ by the same magnitude
$\alpha_t$; only the signs are individual. Suppose then that one entry must
track a target alternating
between $+1$ and $-1$ while another must track a small constant
$-\varepsilon$. The alternating entry keeps its residual, and with it
$\alpha_t$, at $\Theta(1)$; the constant entry is therefore displaced by
$\pm\Theta(1)$ at every step and can only oscillate about its target, never
settle on it. Which side of the target the oscillation favours is decided by
its \emph{phase}, and a one-bit sign carries no information by which a phase
could be corrected. In the unfavourable phase the estimate of the constant
entry has a time average of the sign \emph{opposite} to $-\varepsilon$, and
the iterate driven by that estimate moves the wrong way forever.

\emph{From the sketch to an instance.} The sketch is not yet a counterexample:
in EF21-SignMuon the quantity tracked is the polar factor
$\mathbf D_t=\operatorname{polar}(\tilde{\mathbf M}_t)$, of unit spectral
norm, and the gradients behind it must all come from one smooth
function. Both constraints are met at size $2\times2$ by the reflections
$\bar{\mathbf D}^{\pm}=\smat{a}{\pm b}{\pm b}{-a}$ with $a^2+b^2=1$: one
matrix the oracle can emit carries both roles of the sketch at once, the
large alternating entries on the off-diagonal ($b=\tfrac{24}{25}$) and the
small constant ones on the diagonal ($a=\tfrac7{25}$), coupled by the shared
$\alpha_t$. On these targets the estimate enters a period-two cycle in which
the time average of the $(2,2)$-entry is positive although every target value
is $-\tfrac7{25}$, so $(\mathbf X_t)_{22}$ travels to $-\infty$, the direction
in which the objective increases (\hyperref[fig:divergence_plot]{Figure~\ref*{fig:divergence_plot}}, right). Nor does the mechanism rest on
degeneracy: the momentum matrices of the divergent tail have condition number
$\tfrac{16}{9}$ throughout.

\emph{The role of the preamble.} The unfavourable phase must be arranged.
From $\mathbf d_0^{\mathrm{est}}=\mathbf 0$, the alternating targets alone
lock the estimate into a period-two cycle whose diagonal average has the
\emph{correct} sign. The recursion \eqref{eq:efsm_est} is piecewise affine,
the pieces indexed by the sign pattern of the residual, and the harmless
cycle and the wrong-sign one lie in different pieces, which the dynamics
cannot join. The target sequence of Part~2 therefore opens with two preamble
steps, whose sole purpose is to place the estimate in the piece containing
the wrong-sign cycle.

\emph{Eliminating the parameters.} The step size and the smoothness constant
only rescale the trajectory, so $\eta=1$ and one value of $L$ suffice
(Part~1). Momentum determines only which gradients produce a given target
sequence, not how the recursion \eqref{eq:efsm_est} responds to it; solving
the momentum recursion for the gradients therefore settles every
$(\mu,\text{variant})$ at once (Part~3).

Accordingly the proof has three parts, and they are independent. Part~1
removes $L$ and $\eta$ by rescaling. Part~2 carries the dynamics, and is a
finite computation in exact rational arithmetic: on one fixed sequence of LMO
targets, the estimate enters a limit cycle whose diagonal has the wrong sign.
Part~3 is an existence argument only: it exhibits a smooth function whose
gradients generate that target sequence for every $\mu$ and either momentum
variant. A reader willing to grant that some smooth objective produces the
targets can stop after Part~2.

\subsubsection*{Proof}

\paragraph{Part 1: rescaling.}
\begin{lemma}[Scale reduction]\label{lem:efsm_reduction}
If $\tilde f$ is $\tilde L$-smooth with EF21-SignMuon iterates
$\tilde{\mathbf X}_t$ at step size $1$, then
$f(\mathbf X):=\tfrac{L\eta^2}{\tilde L}\,\tilde f(\mathbf X/\eta)$ is $L$-smooth
and its run at step size $\eta$ (same $\mu$, same variant, from $\mathbf 0$)
satisfies $\mathbf X_t=\eta\,\tilde{\mathbf X}_t$ and
$f(\mathbf X_t)=\tfrac{L\eta^2}{\tilde L}\,\tilde f(\tilde{\mathbf X}_t)$.
\end{lemma}
\paragraph{Proof.}
$\nabla f(\mathbf X)=\tfrac{L\eta}{\tilde L}\nabla\tilde f(\mathbf X/\eta)$, so
$f$ is $L$-smooth. Suppose $\mathbf X_s=\eta\tilde{\mathbf X}_s$ for $s<t$: the
two gradients then differ by the positive factor $\tfrac{L\eta}{\tilde L}$,
which $\mathbf M_t$ and $\tilde{\mathbf M}_t$ (positive combinations of past
gradients) inherit and $\operatorname{polar}$ discards. Hence $\mathbf D_t$ and
$\mathbf d_t^{\mathrm{est}}$ match the normalized run, and
$\mathbf X_t=\eta\tilde{\mathbf X}_t$. $\blacksquare$

\medskip\noindent
It therefore suffices to exhibit, for each $(\mu,\text{variant})$, a $C^\infty$
$\tilde L$-smooth $\tilde f$ on which the \emph{normalized} run ($\eta=1$) obeys
\eqref{eq:exact_rate_app}; that $\tilde f$ may be taken bounded below
(\hyperref[as:1]{Assumption~\ref*{as:1}}) is shown afterwards (\hyperref[rem:efsm_bdd]{Remark~\ref*{rem:efsm_bdd}}). We fix
$\eta=1$ from now on.

\paragraph{Part 2: the limit cycle of the estimate.}
The dynamical core of the proof is the behavior of the recursion
\eqref{eq:efsm_est} on the fixed target sequence
\begin{equation}\label{eq:target_seq}
\mathbf D_1=\mathbf S_1,\quad \mathbf D_2=\mathbf S_2,\quad
\mathbf D_t=\bar{\mathbf D}^{(-1)^t}\ (t\ge3),
\end{equation}
\begin{equation}\label{eq:targets}
\begin{gathered}
\bar{\mathbf D}^{\pm}=\smat{7/25}{\pm24/25}{\pm24/25}{-7/25},\\
\mathbf S_1=\smat{-4/5}{3/5}{-3/5}{-4/5},\quad
\mathbf S_2=\smat{3/5}{-4/5}{0}{0}.
\end{gathered}
\end{equation}
The divergence originates in the tail ($t\ge3$): the reflections
$\bar{\mathbf D}^{\pm}$ share the diagonal $(\tfrac7{25},-\tfrac7{25})$, while
their off-diagonal $\pm\tfrac{24}{25}$ reverses sign at every step. The rotation
$\mathbf S_1$ and the rank-one $\mathbf S_2$ form a two-step preamble. Only
$\mathbf S_2$ requires comment: at a rank-deficient argument the spectral-ball
LMO is not unique, and \hyperref[sec:proState]{Section~\ref*{sec:proState}} resolves it through the thin
SVD that retains only the nonzero singular directions,
$\operatorname{polar}(\mathbf M)=\mathbf U\mathbf V^\top$ with one column of
$\mathbf U,\mathbf V$ per nonzero singular value of $\mathbf M$; under that
convention a rank-one matrix of unit spectral norm, such as $\mathbf S_2$, is
its own polar factor, and this is also what the implementation computes. The next lemma says
what the preamble is for, and that something like it is unavoidable.

\begin{lemma}[The cycle reached without the preamble]\label{lem:efsm_twocycles}
Write $\bar{\mathbf D}^{\pm}=\smat{a}{\pm b}{\pm b}{-a}$ with $a^2+b^2=1$ and
$0<a<b$. On the purely alternating targets
$\mathbf D_t=\bar{\mathbf D}^{(-1)^t}$ ($t\ge1$), the recursion
\eqref{eq:efsm_est} started from $\mathbf d_0^{\mathrm{est}}=\mathbf 0$ enters a
period-two cycle immediately, namely
$\mathbf d_t^{\mathrm{est}}=\tfrac{a+b}{2}\smat{1}{-1}{-1}{-1}$ for odd $t$ and
$\tfrac{b-a}{2}\smat{-1}{1}{1}{1}$ for even $t$. Over a period its
$(2,2)$-entry averages $-\tfrac a2$: the \emph{same} sign as every target
value $-a$, at half the magnitude.
\end{lemma}
\paragraph{Proof.}
Put $m=\tfrac{a+b}2$. The residual $\bar{\mathbf D}^{-}-\mathbf 0$ has signs
$\smat{+}{-}{-}{-}$ and mean modulus $m$, so
$\mathbf d_1^{\mathrm{est}}=m\smat{1}{-1}{-1}{-1}$. Next,
$\bar{\mathbf D}^{+}-\mathbf d_1^{\mathrm{est}}
=\smat{a-m}{b+m}{b+m}{m-a}$ has signs $\smat{-}{+}{+}{+}$ (as $a<m$) and mean
modulus $b$, giving $\mathbf d_2^{\mathrm{est}}=\tfrac{b-a}2\smat{-1}{1}{1}{1}$.
Repeating once returns $\mathbf d_1^{\mathrm{est}}$. The average of the two
$(2,2)$-entries is $\tfrac12\bigl(-m+\tfrac{b-a}2\bigr)=-\tfrac a2$.
$\blacksquare$

That cycle is harmless: the alternation by itself does not diverge, and no
choice of $(a,b)$ makes it do so. The divergence comes from a \emph{second}
period-two cycle of the same recursion, whose residuals carry a \emph{uniform}
sign pattern
($\smat{+}{+}{+}{+}$ and $\smat{-}{-}{-}{-}$) where those of the harmless
cycle are mixed. Since \eqref{eq:efsm_est} is affine on each sign-pattern
cell, the dynamics cannot pass from one cycle to the other; the preamble
exists solely to place $\mathbf d_2^{\mathrm{est}}$ in the cell of the
wrong-sign cycle, which is what \hyperref[lem:efsm_cycle]{Lemma~\ref*{lem:efsm_cycle}} verifies.

\begin{lemma}[Wrong-sign limit cycle]\label{lem:efsm_cycle}
On the targets \eqref{eq:target_seq}, the recursion \eqref{eq:efsm_est} from
$\mathbf d_0^{\mathrm{est}}=\mathbf 0$ enters at $t=3$ the exact period-two cycle
$\mathbf d_t^{\mathrm{est}}=\mathbf d_B$ (odd $t$), $\mathbf d_A$ (even $t$),
where
\begin{equation}\label{eq:dcycle}
\mathbf d_B=\tfrac1{200}\smat{-61}{-201}{-61}{-61},\qquad
\mathbf d_A=\tfrac1{200}\smat{131}{-9}{131}{131}.
\end{equation}
Over a period its $(2,2)$-entry averages
$\tfrac12\bigl((\mathbf d_A)_{22}+(\mathbf d_B)_{22}\bigr)
=+\tfrac7{40}$, opposite in sign to every target value
$(\mathbf D_t)_{22}=-\tfrac7{25}$.
\end{lemma}
\paragraph{Proof.}
Substituting \eqref{eq:target_seq} into \eqref{eq:efsm_est} gives the values
\begin{equation}\label{eq:efsm_table}
\renewcommand{\arraystretch}{1.25}
\begin{array}{c|cccc}
t & 1 & 2 & 3 & 4\ (\text{then }2\text{-periodic})\\\hline
\alpha_t & \tfrac7{10} & \tfrac{21}{20} & \tfrac{131}{200} & \tfrac{24}{25}\\
\mathbf d_t^{\mathrm{est}} & \tfrac1{10}\smat{-7}{7}{-7}{-7}
 & \tfrac1{20}\smat{7}{-7}{7}{7} & \mathbf d_B & \mathbf d_A
\end{array}
\end{equation}
Each residual $\mathbf D_t-\mathbf d_{t-1}^{\mathrm{est}}$ has strictly nonzero
entries, so the signs are unambiguous; e.g.\ at $t=3$ it is
$\tfrac1{100}\smat{-7}{-61}{-131}{-63}$, all negative. From $t=3$ the targets
are $2$-periodic and the pair $(\mathbf d_B,\mathbf d_A)$ reproduces itself:
$\bar{\mathbf D}^{+}-\mathbf d_B$ and $-(\bar{\mathbf D}^{-}-\mathbf d_A)$ are
both entrywise positive with mean $\tfrac{24}{25}$, so
$\mathbf d_B+\tfrac{24}{25}\mathbf J=\mathbf d_A$ and
$\mathbf d_A-\tfrac{24}{25}\mathbf J=\mathbf d_B$ ($\mathbf J$ all-ones). $\blacksquare$

\medskip\noindent
Since $\mathbf X_t=\mathbf X_{t-1}-\mathbf d_t^{\mathrm{est}}$ (recall
$\eta=1$), over one period the $(2,2)$-coordinate changes by
$-(\mathbf d_A+\mathbf d_B)_{22}=-\tfrac7{20}$. Consequently a term $-\gamma
W_{22}$ in $\tilde f$, with $\gamma:=\tfrac7{12}$, increases by
$\gamma\cdot\tfrac7{20}=\tfrac{49}{240}$ per period. This is the divergence,
provided a genuine smooth function produces the targets \eqref{eq:target_seq};
Part~3 constructs one.

\paragraph{Part 3: realization by a smooth function.}
Fix $(\mu,\text{variant})$ and set $\gamma=\tfrac7{12}$, $\nu=\tfrac1{1+2\mu}$,
and
\begin{equation}\label{eq:Amu}
\begin{gathered}
A=\frac{1+\mu}{1-\mu}\ \text{(standard)},\\
A=\frac{1+\mu}{(1-\mu)(1+2\mu)}\ \text{(Nesterov)}.
\end{gathered}
\end{equation}
We build $\tilde f=g+\sum_{k=1}^{3}b_k$ from a periodic-plus-linear
\emph{field} $g$ and three \emph{localized corrections} $b_k$, all explicit.

\begin{figure*}[t]
\centering
\includegraphics[width=\textwidth]{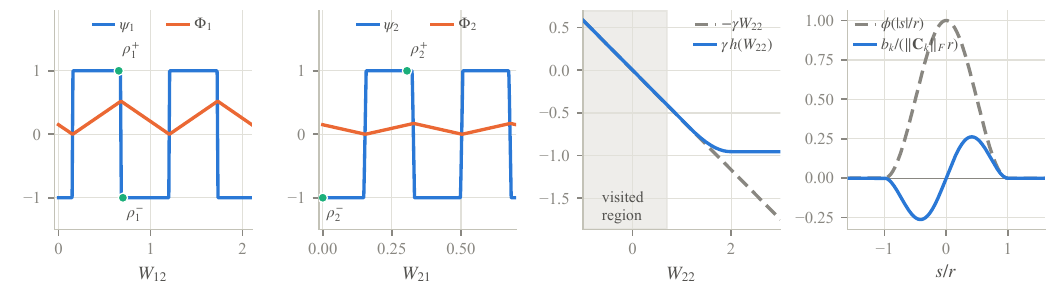}
\caption{The components of the objective $\tilde f=g+\sum_k b_k$ of Part~3, as
implemented. \emph{First two panels:} the periodic ramps $\psi_1,\psi_2$ and
their bounded antiderivatives $\Phi_1,\Phi_2$ over two periods; the marked
residues $\rho_i^{+}$ (kept at odd iterates) and $\rho_i^{-}$ (even iterates)
lie on the plateaus where $\psi_i=\pm1$ exactly, so from $t\ge4$ the
off-diagonal gradient entries alternate between $+A$ and $-A$. The drawn
$\psi_i$ is the implementation's ramp; the trajectory samples only the
plateaus, on which any $C^\infty$ choice agrees with it. \emph{Third panel:}
the divergence term $-\gamma W_{22}$ and its bounded replacement
$\gamma\,h(W_{22})$ of \hyperref[rem:efsm_bdd]{Remark~\ref*{rem:efsm_bdd}}; the two agree on the visited
region $\{W_{22}\le\tfrac7{10}\}$ (shaded), so the floor changes no iterate
while restoring \hyperref[as:1]{Assumption~\ref*{as:1}}. \emph{Fourth panel:} the correction
$b_k$ of \eqref{eq:bump} along the ray
$\mathbf Z_k+s\,\mathbf C_k/\|\mathbf C_k\|_F$, normalized by
$\|\mathbf C_k\|_F\,r$, with the cutoff $\phi$: $b_k$ vanishes at $\mathbf Z_k$
while $\nabla b_k(\mathbf Z_k)=\mathbf C_k$, and its support, the ball
$\{\|\mathbf W-\mathbf Z_k\|_F\le r\}$, contains no iterate other than its own
center.}
\label{fig:ef21_construction}
\end{figure*}

\emph{The field} is
\begin{equation}\label{eq:field}
g(\mathbf W)=-\gamma\,W_{22}+A\,\Phi_1(W_{12})+A\,\Phi_2(W_{21}),
\end{equation}
where $\Phi_i(w)=\int_0^{w}\psi_i$ and $\psi_i\colon\mathbb R\to[-1,1]$ is a
fixed $C^\infty$, $p_i$-periodic function with $\int_0^{p_i}\psi_i=0$,
\[
p_1=\tfrac{21}{20},\quad p_2=\tfrac7{20},\qquad \delta=\tfrac1{100},
\]
equal to $+1$ on $[\rho_i^{+}\!-\delta,\rho_i^{+}\!+\delta]$ and $-1$ on
$[\rho_i^{-}\!-\delta,\rho_i^{-}\!+\delta]$ (mod $p_i$), where
\[
\rho_1^{+}=\tfrac{131}{200},\ \rho_1^{-}=\tfrac{140}{200};\qquad
\rho_2^{+}=\tfrac{61}{200},\ \rho_2^{-}=0.
\]
The two intervals are disjoint mod $p_i$ (their centers are
$\tfrac9{200}>2\delta$ apart), so such a $\psi_i$ exists; the zero-mean
condition, met by balancing the rest of the period, makes $\Phi_i$ periodic
(hence bounded). On each of the two intervals $\Phi_i'=\psi_i=\pm1$ exactly.
\hyperref[fig:ef21_construction]{Figure~\ref*{fig:ef21_construction}} draws both ramps, together with the
remaining components of $\tilde f$.

\emph{The corrections} pin the first three gradients. Fix once a $C^\infty$
cutoff $\phi\colon\mathbb R\to[0,1]$ with $\phi\equiv1$ on $[0,\tfrac12]$ and
$\phi\equiv0$ on $[1,\infty)$, put $r=\tfrac1{50}$, and set
\begin{equation}\label{eq:bump}
b_k(\mathbf W)=\bigl\langle\mathbf C_k,\,\mathbf W-\mathbf Z_k\bigr\rangle\,
\phi\!\bigl(\|\mathbf W-\mathbf Z_k\|_F/r\bigr),
\end{equation}
centered at the first three iterates (from \hyperref[lem:efsm_cycle]{Lemma~\ref*{lem:efsm_cycle}})
\[
\mathbf Z_1=\mathbf 0,\quad
\mathbf Z_2=\tfrac1{10}\smat{7}{-7}{7}{7},\quad
\mathbf Z_3=\tfrac1{20}\smat{7}{-7}{7}{7}.
\]
Each $b_k$ is $C^\infty$, supported in $\{\|\mathbf W-\mathbf Z_k\|_F\le r\}$;
since its linear factor vanishes at $\mathbf Z_k$ while $\phi(0)=1$,
$\nabla b_k(\mathbf Z_k)=\mathbf C_k$. We choose
$\mathbf C_k:=\hat{\mathbf G}_k-\nabla g(\mathbf Z_k)$ with the explicit
$\hat{\mathbf G}_k$ of \eqref{eq:transient} below, so that
$\nabla\tilde f(\mathbf Z_k)=\hat{\mathbf G}_k$. Every term of $\tilde f$ has a
bounded Hessian, so $\tilde f$ is $C^\infty$ and $\tilde L$-smooth for a finite
$\tilde L(\mu,\text{variant})$.

\begin{lemma}[Realization]\label{lem:efsm_realize}
For every $\mu\in[0,1)$ and either variant, EF21-SignMuon on this $\tilde f$
($\eta=1$, from $\mathbf 0$) generates exactly the targets
\eqref{eq:target_seq}.
\end{lemma}
\paragraph{Proof.}
\emph{The required gradients.} The buffer recursion
$\mathbf M_t=\mu\mathbf M_{t-1}+(1-\mu)\mathbf G_t$ of \eqref{eq:signa_general}
can be solved for the gradients: prescribing the buffers $(\mathbf M_t)_{t\ge1}$
forces $\mathbf G_t=(\mathbf M_t-\mu\mathbf M_{t-1})/(1-\mu)$, and these are
the gradients the function must deliver. Define accordingly the
\emph{transient gradients}
\begin{equation}\label{eq:transient}
\hat{\mathbf G}_t=\frac{\mathbf M_t-\mu\mathbf M_{t-1}}{1-\mu}\quad(t=1,2,3;\ \mathbf M_0=\mathbf 0)
\end{equation}
and the \emph{field gradient}
$\mathbf G^{\pm}=\smat{0}{\pm A}{\pm A}{-\gamma}$, where the prescribed buffer
values $\mathbf M_{1,2,3}$ are given below. With
$\bar{\mathbf M}^{\pm}=\smat{0}{\pm1}{\pm1}{-\gamma}$, the factorization
\begin{equation}\label{eq:polarM}
\begin{gathered}
\bar{\mathbf M}^{\pm}
=\bar{\mathbf D}^{\pm}\,\smat{24/25}{\mp7/25}{\mp7/25}{337/300}\\
(\text{second factor }\succ0,\ \det=1)
\end{gathered}
\end{equation}
shows $\operatorname{polar}(\bar{\mathbf M}^{\pm})=\bar{\mathbf D}^{\pm}$, with
singular values $\tfrac34$ and $\tfrac43$; this is the condition number
$\tfrac{16}9$ cited in the proof idea.

\emph{Standard momentum} ($\tilde{\mathbf M}_t=\mathbf M_t$). Take
$\mathbf M_1=\mathbf S_1$, $\mathbf M_2=\mathbf S_2$,
$\mathbf M_3=\bar{\mathbf M}^{-}$; then $\hat{\mathbf G}_{1,2,3}$ are the
explicit matrices \eqref{eq:transient}. Feeding $\mathbf G^{(-1)^t}$ for $t\ge4$
keeps $\mathbf M_t=\bar{\mathbf M}^{(-1)^t}$, because
$(1-\mu)\mathbf G^{\pm}=\bar{\mathbf M}^{\pm}-\mu\bar{\mathbf M}^{\mp}$ with $A$
as in \eqref{eq:Amu}. As $\operatorname{polar}(\mathbf S_1)=\mathbf S_1$
(orthogonal) and $\operatorname{polar}(\mathbf S_2)=\mathbf S_2$ (rank one),
\eqref{eq:polarM} yields the targets \eqref{eq:target_seq}.

\emph{Nesterov momentum} ($\tilde{\mathbf M}_t=(1+\mu)\mathbf M_t-\mu\mathbf M_{t-1}$).
It is here that the orthogonality of $\mathbf S_1$ is needed. The Nesterov
direction involves two consecutive buffers, so prescribing the tail $t\ge4$
already fixes $\mathbf M_3$, and only the preamble remains to absorb the
mismatch at $t=3$; and the set of matrices with a given polar factor
$\mathbf D$ is $\{\mathbf D\mathbf H:\mathbf H\succ0\}$, which is
three-dimensional when $\mathbf D$ is orthogonal but only one-dimensional (a
positive scalar) when $\mathbf D$ has rank one. A preamble of two rank-one
targets leaves too little freedom: a symbolic check shows that it admits no
realization once $\mu\gtrsim0.3$. One orthogonal target supplies enough
freedom, of which the construction uses a single scalar, the $\tau$ below.
Take $\mathbf M_1=\tfrac1{1+\mu}\mathbf S_1\mathbf H_1$,
$\mathbf M_2=\tfrac1{1+\mu}(\mathbf S_2+\mu\mathbf M_1)$,
$\mathbf M_3=\bar{\mathbf M}'^{-}$ and $\mathbf M_t=\bar{\mathbf M}'^{(-1)^t}$
($t\ge4$), where $\bar{\mathbf M}'^{\pm}=\smat{0}{\pm\nu}{\pm\nu}{-\gamma}$,
$\mathbf H_1=\operatorname{diag}(1,1+\tau)$, and, for $\mu>0$,
\[
\tau=\Bigl(\tfrac{140(1+\mu)^2}{1+2\mu}-44\Bigr)\tfrac{1+\mu}{117\mu}>0 .
\]
(At $\mu=0$ the Nesterov rule reads $\tilde{\mathbf M}_t=\mathbf M_t$ and is the standard case already treated, so nothing is left to prove there.)
Then $\tilde{\mathbf M}_1=\mathbf S_1\mathbf H_1$ and
$\tilde{\mathbf M}_2=\mathbf S_2$ have polar factors $\mathbf S_1,\mathbf S_2$,
and $\tilde{\mathbf M}_t=\bar{\mathbf M}^{(-1)^t}$ for $t\ge4$ (since
$\nu(1+2\mu)=1$). The one nontrivial step is $t=3$. Set
$\mathbf H_3:=\bar{\mathbf D}^{-}\tilde{\mathbf M}_3$; since
$(\bar{\mathbf D}^{-})^2=\mathbf I$, this is the same as
$\tilde{\mathbf M}_3=\bar{\mathbf D}^{-}\mathbf H_3$. The stated $\tau$ is
exactly the value making $\mathbf H_3$ symmetric,
and $\mathbf H_3\succ0$ for all $\mu\in(0,1)$: under $\mu=\tfrac{s}{1+s}$ ($s>0$)
the numerators of its two leading minors are polynomials in $s$ with
nonnegative coefficients. Hence
$\operatorname{polar}(\tilde{\mathbf M}_3)=\bar{\mathbf D}^{-}$, completing
\eqref{eq:target_seq}.

\emph{The function delivers these gradients.} It remains to verify
$\nabla\tilde f(\tilde{\mathbf X}_{t-1})=\hat{\mathbf G}_t$ for $t\le3$ and
$=\mathbf G^{(-1)^t}$ for $t\ge4$, by induction along the run: as long as the
gradients match this prescription, the iterates are those computed in
Part~2, and the prescription need only be checked at those points. By
\hyperref[lem:efsm_cycle]{Lemma~\ref*{lem:efsm_cycle}} the iterates
$\tilde{\mathbf X}_0,\tilde{\mathbf X}_1,\tilde{\mathbf X}_2$ are exactly the
centers $\mathbf Z_1,\mathbf Z_2,\mathbf Z_3$, where
$\nabla\tilde f=\hat{\mathbf G}_{1,2,3}$ by construction. The three balls are
disjoint ($\|\mathbf Z_j-\mathbf Z_k\|_F\ge\tfrac7{10}>2r$), and every later
iterate has $(1,2)$-entry $\ge\tfrac{131}{200}$ while the centers have it
$\le0$, so no ball is re-entered. For $t\ge4$ the query lies in the field,
where
\[
\nabla g(\mathbf W)=\begin{pmatrix}0 & A\psi_1(W_{12})\\ A\psi_2(W_{21}) & -\gamma\end{pmatrix}.
\]
The period shift $\tilde{\mathbf X}_{t+2}-\tilde{\mathbf X}_t
=-(\mathbf d_A+\mathbf d_B)$ advances $W_{12}$ by exactly $+p_1$ and $W_{21}$
by exactly $-p_2$ per period, so the two coordinates travel to $+\infty$ and
$-\infty$ respectively while, mod $p_i$, they hold the residues $\rho_i^{+}$
at odd indices and $\rho_i^{-}$ at even ones; there $\psi_i=\pm1$, giving
$\nabla g=\mathbf G^{(-1)^t}$. $\blacksquare$

\paragraph{Proof of \hyperref[th:ef_div]{Theorem~\ref*{th:ef_div}}.}
By \hyperref[lem:efsm_realize]{Lemma~\ref*{lem:efsm_realize}} the normalized run produces the targets
\eqref{eq:target_seq}, so by \hyperref[lem:efsm_cycle]{Lemma~\ref*{lem:efsm_cycle}} its estimate locks onto
the wrong-sign cycle and $\tilde{\mathbf X}_{t+2}-\tilde{\mathbf X}_t$ is a
constant shift. Along it $\Phi_1,\Phi_2$ return to their values and the
corrections $b_k$ vanish, so only the linear term acts:
$\tilde f(\tilde{\mathbf X}_{t+2})-\tilde f(\tilde{\mathbf X}_t)
=-\gamma(-\tfrac7{20})=\tfrac{49}{240}>0$. By \hyperref[lem:efsm_reduction]{Lemma~\ref*{lem:efsm_reduction}},
$f$ then obeys \eqref{eq:exact_rate_app} with $c=\tfrac{49}{240\tilde L}$. As
$(L,\eta,\mu,\text{variant})$ were arbitrary, no rule $\eta=\eta(L,\mu)$ can
prevent divergence. $\blacksquare$

\begin{remark}[Boundedness below]\label{rem:efsm_bdd}
Only the linear term $-\gamma W_{22}$ makes $\tilde f$ unbounded below, and only
as $W_{22}\to+\infty$, a region the iterates never enter, since
$(\tilde{\mathbf X}_t)_{22}\le\tfrac7{10}$ throughout (it decreases after the
transient). Replacing $-\gamma W_{22}$ by any $C^\infty$ function that agrees
with it on $\{W_{22}\le1\}$ and is constant on $\{W_{22}\ge2\}$ therefore leaves
the whole trajectory, and \eqref{eq:exact_rate_app} with it, unchanged while rendering
$\tilde f$ bounded below (\hyperref[as:1]{Assumption~\ref*{as:1}}); the theorem is stated with
this modification in force.
\end{remark}

\begin{figure*}[t]
\centering
\includegraphics[width=\textwidth]{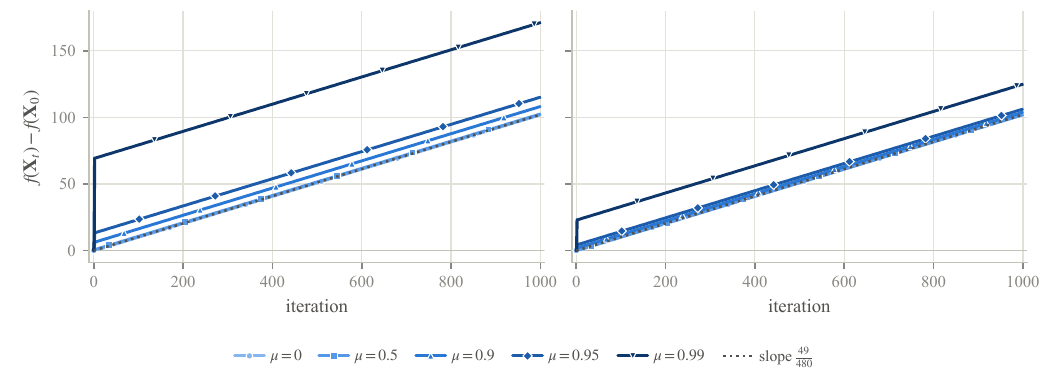}
\caption{Momentum does not prevent the divergence of EF21-SignMuon. For each momentum coefficient
$\mu$ and variant, EF21-SignMuon is run on the corresponding instance of
\hyperref[th:ef_div]{Theorem~\ref*{th:ef_div}} and $f(\mathbf X_t)-f(\mathbf X_0)$ is plotted;
\emph{left:} standard momentum, \emph{right:} Nesterov. Every setting diverges
at the common rate $\tfrac{49}{480}$ (dotted). Subtracting $f(\mathbf X_0)$
removes the only genuinely $\mu$-dependent offset (the field constant $A(\mu)$
scales a bounded periodic term); what remains is a bounded transient, largest as
$\mu\to1$, on top of the shared linear divergence.}
\label{fig:ef21_momentum}
\end{figure*}

\begin{remark}[The construction is not convex]\label{rem:efsm_nonconvex}
\hyperref[th:1]{Theorems~\ref*{th:1}}--\ref{th:3} run on a linear, hence convex, objective;
$\tilde f$ is nonconvex, through the periodic terms $A\Phi_i$ and the
corrections $b_k$. The nonconvexity is forced by the run rather than chosen by
the realization: along the divergent trajectory,
\[
\bigl\langle\nabla\tilde f(\tilde{\mathbf X}_6)-\nabla\tilde f(\tilde{\mathbf X}_3),\
\tilde{\mathbf X}_6-\tilde{\mathbf X}_3\bigr\rangle=-\tfrac{9A}{50}<0,
\]
violating the gradient monotonicity that every convex function obeys, so
\emph{no} convex function generates these iterates and gradients. The theorem
is stated under \hyperref[as:1]{Assumptions~\ref*{as:1}}--\ref{as:2} because that is where it is
needed: EF21-MuonUSign and EF21-MuonSign converge under exactly these
hypotheses (\hyperref[thm:conv]{Theorem~\ref*{thm:conv}}), so the divergence and the guarantees
concern one problem class. Whether some convex instance, necessarily through a
different trajectory, also defeats EF21-SignMuon we leave open.
\end{remark}

\begin{remark}[Verification]\label{rem:efsm_scope}
The construction is checked in two independent ways: symbolically, in exact
rational arithmetic, and numerically, by running the float64 reference
implementation of \hyperref[ef21_signmuon]{Algorithm~\ref*{ef21_signmuon}} on the assembled $\tilde f$.
The right panel of \hyperref[fig:divergence_plot]{Figure~\ref*{fig:divergence_plot}} confirms that at $\mu=0$
EF21-SignMuon is the only one of the eight methods that diverges, the others
(SignMuon, MuonUSign, MuonSign, EF21-MuonUSign, EF21-MuonSign, SignSGD, Muon)
staying bounded; \hyperref[fig:ef21_momentum]{Figure~\ref*{fig:ef21_momentum}} confirms that EF21-SignMuon
diverges at the exact rate $\tfrac{49}{480}$ for every
$\mu\in\{0,\tfrac12,0.9,0.95,0.99\}$ under both standard and Nesterov momentum, as the
reduction predicts.
\end{remark}
%%% <<< end ef21_signmuon_divergence.tex

%%% >>> begin ef21_musign_reduction.tex (was \input on signmuon_appendix.tex:191)
% =====================================================================
%  Convergence of EF21-MuonUSign and EF21-MuonSign by exact reduction to
%  the EF21-Muon framework of Gruntkowska, Burlachenko & Richtarik (2025),
%  arXiv:2510.00643 (source: EF_Muon_Friends/template.tex).
%
%  We do NOT restate their per-layer bounds. The reduction rests on one
%  non-trivial fact -- the scaled sign is a contractive compressor (Lemma) --
%  plus bookkeeping (LMO step, momentum, loop order). Rates and stepsize
%  rules are in our notation; the constants are their Theorems 19 and 24
%  under the substitution of Table 3.
%  Assumptions as:1/as:2 are stated in the main text (Problem Statement).
% =====================================================================
\subsection{Convergence of EF21-MuonUSign and EF21-MuonSign}\label{app:hyp}

We do not analyse the two error-feedback methods from scratch. EF21-MuonUSign and EF21-MuonSign are \emph{exact instances} of EF21-Muon \citep[Algorithm~3]{gruntkowska2025error}, already analysed in the layer-wise, stochastic, federated setting. Three conditions must be verified before its guarantees transfer: our step is their LMO step, our loop is their loop (\hyperref[prop:instance]{Proposition~\ref*{prop:instance}}), and our messages come from contractive compressors (\hyperref[lem:signcontr]{Lemma~\ref*{lem:signcontr}}). Only the last is non-trivial. \hyperref[tab:dict]{Table~\ref*{tab:dict}} is the change of variables.

\paragraph{Notation and constants.} For the layer tuple $\mathbf{X} = [\mathbf{X}_1,\dots,\mathbf{X}_p]$, $\mathbf{X}_i \in \mathbb{R}^{m_i\times n_i}$, of the Problem Statement write $d_i := m_i n_i$, $r_i := \min(m_i,n_i)$, $d_{\max} := \max_i d_i$; a second subscript selects a layer ($\mathbf{X}_{t,i}$, $\mathbf{g}_{t,i}$). We use $\|\mathbf{Y}\|_{2\to2} \le \|\mathbf{Y}\|_F \le \|\mathbf{Y}\|_{*} \le \sqrt{\operatorname{rank}\mathbf{Y}}\,\|\mathbf{Y}\|_F$. Smoothness constants: $L_i$ for $f$ and $L_{i,j}$ for $f_j$ in \hyperref[as:2]{Assumption~\ref*{as:2}}, with $\tilde{L}_i^2 := \frac{1}{N}\sum_j L_{i,j}^2$ and $L := \max_i L_i$; in the $(L^0,L^1)$ form the pairs are again per layer for $f$ and per layer and client for $f_j$, with $L^1_{i,\max} := \max_j L^1_{i,j}$.

% Assumptions 1-2 are stated in the main text; in the standalone supplement
% the assumption counter starts fresh, so it is set by hand to keep this one
% numbered 3 in both builds.
\setcounter{assumption}{2}
\begin{assumption}[\textbf{Stochastic gradient}]\label{as:3}
Each client's stochastic gradient is unbiased, $\mathbb{E}_{\xi}[\nabla f_j(\mathbf{X};\xi)] = \nabla f_j(\mathbf{X})$, with bounded variance $\mathbb{E}_{\xi}\|\nabla f_j(\mathbf{X};\xi) - \nabla f_j(\mathbf{X})\|_*^2 \le \sigma^2$.
\end{assumption}

\hyperref[as:1]{Assumptions~\ref*{as:1}}--\ref{as:3} are Assumptions~1--2 and~6--10 of \citet{gruntkowska2025error} with the layer norms taken spectral; our variance bound is stated in the nuclear norm and implies theirs via $\|\cdot\|_F \le \|\cdot\|_{*}$. The clause $f_j \ge f_j^*$ of \hyperref[as:1]{Assumption~\ref*{as:1}} is needed only for \hyperref[cor:l0l1]{Corollary~\ref*{cor:l0l1}}.

\subsubsection{Main result}

% Theorems 1-4 are the divergence results (three in Appendix A.2-A.4, one in the
% body); this is Theorem 5, as everything citing it assumes.
\setcounter{theorem}{4}
\begin{theorem}[\textbf{Convergence of the EF21 methods}]\label{thm:conv}
Run the federated \hyperref[alg:fed_serverlmo]{Algorithm~\ref*{alg:fed_serverlmo}} with the EF21 uplink and EMA momentum $\mu \in [0,1)$. Then:
\begin{itemize}
\item[\emph{(i)}] \emph{(smooth; EF21-MuonUSign and EF21-MuonSign)} under \hyperref[as:1]{Assumptions~\ref*{as:1}}--\ref{as:3}, with the ``sharp'' learning rate $\eta_{t,i} = \gamma_i\|\mathbf{g}_{t,i}\|_{*}$ and tuned $(\gamma_i,\mu)$, both methods reach $\tfrac{1}{T}\sum_{t<T}\mathbb{E}\|\nabla f(\mathbf{X}_t)\|_*^2 = \mathcal{O}(T^{-1/2})$ (\hyperref[cor:smooth]{Corollary~\ref*{cor:smooth}});
\item[\emph{(ii)}] \emph{(generalized smooth; EF21-MuonUSign only)} under $(L^0,L^1)$-smoothness, EF21-MuonUSign with a plain \emph{constant} learning rate reaches $\min_{t\le T}\sum_i\mathbb{E}\|\nabla_i f(\mathbf{X}_t)\|_* = \mathcal{O}(T^{-1/4})$ (\hyperref[cor:l0l1]{Corollary~\ref*{cor:l0l1}}).
\end{itemize}
Both statements are written for a per-layer constant common to all layers ($\gamma_i \equiv \gamma$, $\eta_i \equiv \eta$); for unequal constants the left-hand sides carry the step-size weights of \hyperref[cor:smooth]{Corollaries~\ref*{cor:smooth}}--\ref{cor:l0l1}.
\end{theorem}

Part~(i) implies the rate announced in the main text, since the minimum over $t\le T$ is at most the average. The centralized \hyperref[central_alg_ef]{Algorithms~\ref*{central_alg_ef}}--\ref{central_alg_ud} are the federated method at $N=1$, so both parts cover the centralized runs as well.

\subsubsection{The reduction}

\paragraph{The step is their LMO step.}
For $\mathbf{G} = \mathbf{U}\boldsymbol{\Sigma}\mathbf{V}^\top$ the framework's oracle over the spectral ball of radius $\tau$ is $\mathrm{lmo}(\mathbf{G}) = \mathbf{X} - \tau\mathbf{U}\mathbf{V}^\top$, which is our server step with $\tau = \eta_{t,i}$, since $\mathbf{D}_{t,i} = -A(\mathbf{g}_{t,i}) = \mathbf{U}_{t,i}\mathbf{V}_{t,i}^\top$:
\begin{equation}\label{eq:lmo_step_ours}
\mathbf{X}_{t,i} = \mathbf{X}_{t-1,i} - \eta_{t,i}\,\mathbf{D}_{t,i} = \mathrm{lmo}_{\mathcal{B}(\mathbf{X}_{t-1,i},\eta_{t,i})}(\mathbf{g}_{t,i}).
\end{equation}
Their ``sharp'' step $\mathbf{X} - \gamma\,\mathbf{G}^{\sharp}$ uses $\mathbf{G}^{\sharp} = \|\mathbf{G}\|_{*}\mathbf{U}\mathbf{V}^\top$, so a constant $\gamma_i$ amounts to the schedule $\eta_{t,i} = \gamma_i\|\mathbf{g}_{t,i}\|_{*}$ and a constant $\eta_{t,i}$ to the plain Muon rate: the two differ as learning-rate choices, not as algorithms, and either run is an instance of the framework. The guarantees, however, are attached to specific choices. Part~(i) of \hyperref[thm:conv]{Theorem~\ref*{thm:conv}} assumes the nuclear-norm schedule, which our implementation does not use; what the corollaries then cover of the constant-rate runs we actually perform is recorded in the Scope paragraph below. The spectral norm-equivalence constants are $\underline{\rho}_i = 1$, $\bar{\rho}_i = \sqrt{r_i}$, defined by $\underline{\rho}_i\|\mathbf{Y}\|_{2\to2}\le\|\mathbf{Y}\|_F\le\bar{\rho}_i\|\mathbf{Y}\|_{2\to2}$. At rank-deficient $\mathbf{G}$ the LMO is non-unique, but any selection serves: writing $\boldsymbol{\Delta} := \mathrm{lmo}_{\mathcal{B}(\mathbf{X},\tau)}(\mathbf{G}) - \mathbf{X}$ for the displacement it produces, the analysis uses only $\langle\mathbf{G},\boldsymbol{\Delta}\rangle = -\tau\|\mathbf{G}\|_{*}$ and $\|\boldsymbol{\Delta}\|_{2\to2}\le\tau$, with $\mathbf{G}=\mathbf{0}$ read as the zero step ($\mathbf{0}^{\sharp}=\mathbf{0}$).

\paragraph{The momentum is their momentum.}
\hyperref[alg:fed_serverlmo]{Algorithm~\ref*{alg:fed_serverlmo}} uses $\mathbf{M}_t = \mu\mathbf{M}_{t-1} + (1-\mu)\mathbf{G}_t$, which is \eqref{eq:signa_general} and is the framework's momentum with $\beta = 1-\mu$. Rescaling the momentum stream by a constant, as the heavy-ball convention does by the factor $1/(1-\mu)$, alters nothing: the EF21 recursion is positively homogeneous and the Muon LMO scale-invariant, so the factor leaves the iterates unchanged under a constant $\eta_{t,i}$ and is absorbed into $\gamma_i$ under the ``sharp'' schedule. (The Nesterov branch is a different filter; see the end of this appendix.)

\begin{proposition}[\textbf{Exact instance}]\label{prop:instance}
Fix $\mu$, a learning-rate schedule, and the LMO selection above. Started from $\mathbf{M}_0^{(j)} = \mathbf{g}_0^{(j)} = \mathbf{g}_0 = \mathbf{0}$, $\mathbf{W}_0 = \mathbf{X}_0$, \hyperref[alg:fed_serverlmo]{Algorithm~\ref*{alg:fed_serverlmo}} with the EF21 uplink and either downlink mode ($\mathcal{C}^{\downarrow} \in \{\text{exact},\text{EF21-P}\}$) produces the same trajectory as Algorithm~3 of \citet{gruntkowska2025error} with spectral norms, scaled-sign worker compressors, identity/scaled-sign server compressor, $\beta = 1-\mu$, and radii $t_i^k = \eta_{k,i}$, up to a one-step index shift $X^{t+1} = \mathbf{X}_t$.
\end{proposition}
\paragraph{Proof} The loops differ only in where the round is cut: they order it \emph{step $\to$ downlink $\to$ gradient $\to$ uplink}, we order it \emph{downlink $\to$ gradient $\to$ uplink $\to$ step}. Their iteration $k=0$ is vacuous under our initialization: $\mathbf{g}_0 = \mathbf{0}$ gives $X^1 = X^0 = \mathbf{X}_0$ and a zero downlink residual, so $W^1 = \mathbf{X}_0$. Thereafter their iteration $k = t$ performs our round $t$ verbatim, with the same momentum, the same compressed residual and the same LMO step \eqref{eq:lmo_step_ours}, giving $X^{t+1} = \mathbf{X}_t$, $W^{t+1} = \mathbf{W}_t$ by induction. Running their method for $K = T+1$ iterations therefore yields $\{\mathbf{X}_0,\mathbf{X}_0,\mathbf{X}_1,\dots,\mathbf{X}_{T-1}\}$, and any average or minimum over these equals ours up to one duplicated nonnegative term. $\blacksquare$

Zero initialization also makes their initial-error constant $\Psi^0$ explicit: the $\|\mathbf{M}_0 - \mathbf{g}_0\|$ term vanishes and the gradient-deviation terms become $\|\nabla_i f(\mathbf{X}_0)\|$. The $\beta^{-1}$ surviving in $\Psi^0$ enters only through the $\Psi^0/T$ term, which $T^{-1/2}$ dominates.

\subsubsection{The scaled sign is a contractive compressor}

The framework requires every transmitted message to originate in a \emph{contractive compressor}: a map with $\mathbb{E}\|\mathcal{C}(\mathbf{Y})-\mathbf{Y}\|^2 \le (1-\alpha)\|\mathbf{Y}\|^2$ for some $\alpha \in (0,1]$ \citep[Def.~1]{gruntkowska2025error}. This is precisely the property that fails for a bare sign and holds once it is scaled.

\begin{lemma}[\textbf{Contractivity of the scaled sign}]\label{lem:signcontr}
For every $\mathbf{Y} \in \mathbb{R}^{m\times n}$ ($d = mn$), the scaled sign $\mathcal{C}(\mathbf{Y}) = \operatorname{mean}(|\mathbf{Y}|)\operatorname{sign}(\mathbf{Y})$, with exact zeros resolved to $\pm1$ as in \hyperref[sec:theory]{Section~\ref*{sec:theory}}, satisfies
\begin{equation}\label{eq:sign_identity}
\|\mathcal{C}(\mathbf{Y})-\mathbf{Y}\|_F^2 = \|\mathbf{Y}\|_F^2 - \tfrac{1}{d}\|\mathbf{Y}\|_1^2 \le \bigl(1-\tfrac{1}{d}\bigr)\|\mathbf{Y}\|_F^2 ,
\end{equation}
so $\mathcal{C}$ is Euclidean-contractive with $\alpha = 1/d$. The identity holds for every draw of the random signs, not merely in expectation, and gives the exact contraction $\alpha(\mathbf{Y}) = \|\mathbf{Y}\|_1^2/(d\|\mathbf{Y}\|_F^2)$, which is $\Theta(1)$ for dense $\mathbf{Y}$ (e.g.\ $\to 2/\pi$ for i.i.d.\ Gaussian entries) and equals $1/d$ exactly at any $1$-sparse $\mathbf{Y}$.
\end{lemma}
\paragraph{Proof} Write $c := \|\mathbf{Y}\|_1/d$ and $s_{kl}=\pm1$ for the transmitted signs. A nonzero entry contributes $(|Y_{kl}|-c)^2$ and a zero entry contributes $(c\,s_{kl})^2 = c^2$ whichever sign was drawn, so
$\|\mathcal{C}(\mathbf{Y})-\mathbf{Y}\|_F^2 = \|\mathbf{Y}\|_F^2 - 2c\|\mathbf{Y}\|_1 + c^2\|\mathbf{Y}\|_0 + c^2(d-\|\mathbf{Y}\|_0) = \|\mathbf{Y}\|_F^2 - 2c\|\mathbf{Y}\|_1 + dc^2$, which is \eqref{eq:sign_identity} on substituting $c$. The bound then follows from $\|\mathbf{Y}\|_1 \ge \|\mathbf{Y}\|_F$, with equality exactly at the $1$-sparse $\mathbf{Y}$; the Gaussian limit uses $\mathbb{E}|Y_{kl}| = \sqrt{2/\pi}\,\sigma$. $\blacksquare$

The randomized $\operatorname{sign}(0)$ is what makes \eqref{eq:sign_identity} hold with equality for every $\mathbf{Y}$: under the ternary convention the $\|\mathbf{Y}\|_0$-terms in the proof do not cancel, the error depends on the sparsity of $\mathbf{Y}$, and $\alpha = 1/d$ becomes an infimum rather than an attained value.

The scaling is essential to \hyperref[lem:signcontr]{Lemma~\ref*{lem:signcontr}}: a bare sign contracts for no $\alpha$ at all, since as $\|\mathbf{Y}\|_F \to 0$ on a fixed support, $\|\operatorname{sign}(\mathbf{Y})-\mathbf{Y}\|_F \to \sqrt{\|\mathbf{Y}\|_0}$. That is the dividing line between our divergent and convergent methods: the majority-vote methods transmit unscaled signs of full quantities, the EF21 variants the scaled sign of a \emph{residual}, at the cost of one extra scalar per layer per round. Per layer the uplink scaled sign lies in $\mathbb{B}_2(1/d_i) \subseteq \mathbb{B}_2(\alpha)$ with $\alpha := 1/d_{\max}$; the downlink is exact for EF21-MuonUSign ($\alpha_P = 1$) and the same scaled sign for EF21-MuonSign.

\subsubsection{Transferred guarantees}

All requirements hold, so the framework's theorems apply through \hyperref[tab:dict]{Table~\ref*{tab:dict}}. We state the rates and stepsize rules; the explicit non-asymptotic bounds are those of Theorems~19 (smooth) and~24 ($(L^0,L^1)$-smooth) of \citet{gruntkowska2025error}, evaluated at the constants of \hyperref[tab:dict]{Table~\ref*{tab:dict}}.

\begin{table}[!t]
\centering\footnotesize
\setlength{\tabcolsep}{4pt}
\begin{tabular}{@{}lll@{}}
\toprule
 & \textbf{This paper} & \textbf{EF21-Muon} \\
\midrule
clients / rounds & $N$, $T$ & $n$, $K = T+1$ \\
iterate & $\mathbf{X}_t$ & $X^{t+1}$ \\
momentum & $\mu$ & $1 - \beta$ \\
learning rate & $\eta_{t,i}$ & radius $t_i^{t}$ \\
norm equivalence & $1,\ \sqrt{r_i}$ & $\underline{\rho}_i,\ \bar{\rho}_i$ \\
uplink compr. & scaled sign & $\mathbb{B}_2(\alpha_D)$ \\
downlink compr. & exact / sign & $\mathcal{I}$ / $\mathbb{B}_2(\alpha_P)^{\dagger}$ \\
compression $\alpha$ & $1/d_{\max}$ & $\alpha_D$ \\
\bottomrule
\end{tabular}
\caption{Change of variables from our notation to that of \citet{gruntkowska2025error}. Layer norms $\|\cdot\|_{(i)},\|\cdot\|_{(i)\star}$ are read as $\|\cdot\|_{2\to2},\|\cdot\|_{*}$ and smoothness constants $L_i^0,\tilde{L}_i^0$ as $L_i,\tilde{L}_i$; the index shift is \hyperref[prop:instance]{Proposition~\ref*{prop:instance}}. $^{\dagger}$Their Theorem~19 as printed asks for a server compressor in $\mathbb{B}(\alpha_P)$, contractive in the \emph{layer} norm; the Euclidean class $\mathbb{B}_2(\alpha_P)$, to which the scaled sign does belong, is admitted by their Remark~23, which introduces the factor $\bar\rho_i^2$ appearing in \hyperref[cor:smooth]{Corollary~\ref*{cor:smooth}} and \hyperref[rem:normmismatch]{Remark~\ref*{rem:normmismatch}}.}
\label{tab:dict}
\end{table}

\begin{corollary}[\textbf{Smooth case; EF21-MuonUSign and EF21-MuonSign}]\label{cor:smooth}
Let \hyperref[as:1]{Assumptions~\ref*{as:1}}--\ref{as:3} hold. Run EF21-MuonUSign or EF21-MuonSign (\hyperref[alg:fed_serverlmo]{Algorithm~\ref*{alg:fed_serverlmo}} with the EF21 uplink and the exact or the scaled-sign downlink, respectively) with EMA momentum and $\eta_{t,i} = \gamma_i\|\mathbf{g}_{t,i}\|_{*}$, for any $\gamma_i$ below the per-layer threshold of \citet[Thm.~19]{gruntkowska2025error} under \hyperref[tab:dict]{Table~\ref*{tab:dict}}. Then, applying the momentum tuning of \citet[Cor.~2]{gruntkowska2025error} per layer (their Corollary~1 supplies the layer-wise initialization, their Corollary~2 the tuning at $p=1$),
\begin{equation*}
\tfrac{1}{T}\textstyle\sum_{t<T}\sum_{i=1}^p w_i\,\mathbb{E}\|\nabla_i f(\mathbf{X}_t)\|_*^2 = \mathcal{O}(T^{-1/2}),
\end{equation*}
with $w_i := \gamma_i/(\tfrac{1}{p}\sum_l\gamma_l)$; for a common $\gamma_i$ all $w_i = 1$ and the left side is $\tfrac{1}{T}\sum_{t<T}\mathbb{E}\|\nabla f(\mathbf{X}_t)\|_*^2$.
\end{corollary}
\paragraph{Proof} By \hyperref[prop:instance]{Proposition~\ref*{prop:instance}} the run is an instance of Algorithm~3 with $K = T+1$, and by \hyperref[lem:signcontr]{Lemma~\ref*{lem:signcontr}} its compressors satisfy $\alpha_D = \alpha$ and $\alpha_P \in \{1,\alpha\}$ (EF21-MuonSign's Euclidean downlink is admissible by their Remark~23, which multiplies the $\alpha_P$-terms of the threshold constant $\zeta_i$ by $\bar{\rho}_i^2 = r_i$). Substituting $\underline{\rho}_i = 1$, $\bar{\rho}_i = \sqrt{r_i}$ into \citet[Thm.~19]{gruntkowska2025error} gives the threshold and the rate; the duplicated $\mathbf{X}_0$-term changes the constant by at most the factor $\frac{T+1}{T}$. $\blacksquare$

The threshold in question is $\gamma_i \le (2L_i + 2\sqrt{\zeta_i})^{-1}$, with $\zeta_i$ the constant of \citet[Thm.~19]{gruntkowska2025error}. Under \hyperref[tab:dict]{Table~\ref*{tab:dict}}, $\zeta_i$ grows polynomially in the layer dimension, through $\bar{\rho}_i^2 = r_i$ and the uplink $\alpha = 1/d_{\max}$; for EF21-MuonSign its $\alpha_P$-terms carry the additional factor $r_i$ of Remark~23, shrinking the admissible $\gamma_i$ by a further $\sqrt{r_i}$. That additional factor is structural: the scaled sign is spectrally contractive for no parameter at all (\hyperref[rem:normmismatch]{Remark~\ref*{rem:normmismatch}}).

\begin{corollary}[\textbf{$(L^0,L^1)$ case; EF21-MuonUSign}]\label{cor:l0l1}
Let \hyperref[as:2]{Assumption~\ref*{as:2}} hold in its $(L^0,L^1)$ form. Run EF21-MuonUSign (exact downlink) with $S := T+2$, momentum $\mu = 1 - S^{-1/2}$, and the \emph{constant} per-layer rate $\eta_{t,i} \equiv \eta_i/S^{3/4}$ for any $\eta_i \le 1$ with $\eta_i^2\sqrt{r_i}\,(L^1_{i,\max})^2 = \mathcal{O}(1)$.\footnote{\citet[Thm.~24]{gruntkowska2025error} imposes four upper bounds on $\eta_i^2$; at $\beta = (K+1)^{-1/2}$ the two quoted here are the binding ones, the first bound relaxing as $T$ grows. Their second bound, as printed, carries a factor $K+1$ in the denominator that the proof does not use: the display it is chosen to ensure requires only $\eta_i^2 \le (1-\sqrt{1-\alpha_D})\,\underline{\rho}_i\,(K+1)^{1/2}/(24\sqrt{1-\alpha_D}\,\bar{\rho}_i(L^1_{i,\max})^2)$, which also relaxes as $T$ grows, and this is the requirement we work from.} Then
\begin{equation*}
\min_{0\le t\le T}\ \textstyle\sum_{i=1}^p v_i\,\mathbb{E}\|\nabla_i f(\mathbf{X}_t)\|_* = \mathcal{O}(T^{-1/4}),\quad v_i := \tfrac{\eta_i}{\frac{1}{p}\sum_l\eta_l}.
\end{equation*}
\end{corollary}
\paragraph{Proof} Their Theorem~24 requires the identity server compressor, i.e.\ the exact downlink of EF21-MuonUSign (so EF21-MuonSign is excluded). Applying it with $K = T+1$, $\beta = S^{-1/2}$ and \hyperref[tab:dict]{Table~\ref*{tab:dict}} gives the schedule and the rate. $\blacksquare$

Under plain smoothness ($L^1 = 0$) every constraint above collapses to $\eta_i \le 1$. The \emph{constant}-learning-rate EF21-MuonUSign we actually run is therefore covered as it stands, and at no loss of order: $\mathcal{O}(T^{-1/4})$ for the norm corresponds to $\mathcal{O}(T^{-1/2})$ for its square, the quantity of part~(i). What the ``sharp'' schedule changes is the bounded quantity, an average of the squared norm in place of a minimum of the norm, not the effective speed.

\subsubsection{Scope}

The reduction covers the two error-feedback methods and no others, which matches our negative results. The majority-vote methods send unscaled signs and fall outside the framework: their server aggregates by a vote, $\operatorname{sign}(\sum_j\mathbf{s}_t^{(j)})$, where Algorithm~3 averages, and it is the average that \hyperref[prop:instance]{Proposition~\ref*{prop:instance}} matches. EF21-SignMuon compresses the LMO \emph{output} rather than the gradient, so it tracks the non-Lipschitz polar factor and the momentum-tracking step breaks (\hyperref[th:ef_div]{Theorem~\ref*{th:ef_div}}). EF21-MuonSign gets the smooth guarantee but not the $(L^0,L^1)$ one, whose theory assumes an uncompressed downlink.

One further gap is ours and not the framework's. \hyperref[cor:smooth]{Corollary~\ref*{cor:smooth}} reaches EF21-MuonSign only under the ``sharp'' schedule $\eta_{t,i} = \gamma_i\|\mathbf{g}_{t,i}\|_{*}$, while \hyperref[cor:l0l1]{Corollary~\ref*{cor:l0l1}}, the constant-rate statement, excludes it for want of an identity server compressor. Our experiments run a cosine-annealed constant $\eta_0$, so the EF21-MuonSign runs are covered by neither, whereas constant-rate EF21-MuonUSign is covered by \hyperref[cor:l0l1]{Corollary~\ref*{cor:l0l1}} at $L^1 = 0$. The trajectories remain exact instances of Algorithm~3 either way; what the constant rate costs is the step-size hypothesis of the rate, not the reduction.

The remaining gap is the Nesterov branch, which our language-model runs use (\hyperref[app:nanogpt]{Appendix~\ref*{app:nanogpt}}): with $\beta := 1-\mu$ it steers by $\tilde{\mathbf{M}}_t = (1-\beta)\mathbf{M}_t + \beta\mathbf{G}_t$ rather than by the buffer $\mathbf{M}_t$, and the framework's analysis is written for a direction that itself satisfies the recursion $\mathbf{M}_t = (1-\beta)\mathbf{M}_{t-1} + \beta\mathbf{G}_t$, which $\tilde{\mathbf{M}}_t$ does not. The discrepancy is small and explicit: writing $\mathbf{n}_t := \mathbf{G}_{t,i} - \nabla_i f(\mathbf{X}_t)$, we have $\nabla_i f - \tilde{\mathbf{M}}_t = (1-\beta)(\nabla_i f - \mathbf{M}_t) - \beta\,\mathbf{n}_t$, and since $\mathbb{E}\langle \nabla_i f - \mathbf{M}_t,\mathbf{n}_t\rangle = -\beta\,\mathbb{E}\|\mathbf{n}_t\|_2^2$,
\begin{equation*}
\mathbb{E}\|\nabla_i f - \tilde{\mathbf{M}}_t\|_2^2 \le (1-\beta)^2\,\mathbb{E}\|\nabla_i f - \mathbf{M}_t\|_2^2 + 3\beta^2\sigma_i^2 .
\end{equation*}
The first term is exactly the deviation \citet[Thm.~19]{gruntkowska2025error} already tracks, contracted rather than enlarged; the second is dominated by the $\beta\sigma_i^2$ term already in its bound. Nesterov should therefore degrade the constants rather than the rate. We state this as an expectation, not as a corollary, since the tracking recursion for $\|\tilde{\mathbf{M}}_{t+1}-\tilde{\mathbf{M}}_t\|$ would have to be redone as well.

\begin{remark}[\textbf{Worst-case and realized contraction}]\label{rem:dimension}
The uplink parameter $\alpha = 1/d_{\max}$ enters the threshold of \hyperref[cor:smooth]{Corollary~\ref*{cor:smooth}} through the $1/\alpha^2$-terms of $\zeta_i$, so the admissible $\gamma_i$ shrinks linearly in $d_{\max}$, the dimension of the largest layer; the same dimension factor arises for Top-$1$ compressors in Euclidean EF21 \citep{richtarik2021ef21}. On the uplink this is a worst case only: by \hyperref[lem:signcontr]{Lemma~\ref*{lem:signcontr}}, $\alpha = 1/d$ requires a residual concentrated on a single coordinate, whereas the momentum residual is dense, with $\alpha(\boldsymbol{\Delta}) = \Theta(1)$, and the framework is stated to extend to iteration-dependent $\alpha$ \citep[Rem.~12]{gruntkowska2025error}, under which the dense value would enter in place of the worst case.

The downlink residual of EF21-MuonSign does not stay dense. Unlike the uplink residual, which every round's gradient refreshes, it is produced by the compressor's own recursion
\begin{equation*}
\boldsymbol{\Delta}_{t+1}^{\downarrow} = \boldsymbol{\Delta}_t^{\downarrow} - \eta_{t+1} \mathbf{D}_{t+1} - \operatorname{mean}|\boldsymbol{\Delta}_t^{\downarrow}|\operatorname{sign}(\boldsymbol{\Delta}_t^{\downarrow}),
\end{equation*}
which corrects every coordinate by the same scalar $\operatorname{mean}|\boldsymbol{\Delta}_t^{\downarrow}|$. A coordinate whose per-step drive $\eta_{t+1}(\mathbf{D}_{t+1})_{kl}$ exceeds that scalar receives a correction smaller than its drive at every step, while the remaining coordinates keep the scalar small, so the residual concentrates on few coordinates; by the exact expression $\alpha(\mathbf{Y}) = \|\mathbf{Y}\|_1^2/(d\|\mathbf{Y}\|_F^2)$ of \hyperref[lem:signcontr]{Lemma~\ref*{lem:signcontr}}, concentration is precisely what lowers $\alpha$. \hyperref[sec:exp_nanogpt]{Section~\ref*{sec:exp_nanogpt}} measures the effect: on the layers built from a zero initialization, $\alpha(\boldsymbol{\Delta}^{\downarrow})$ falls to $1.2\times10^{-4}$, about four orders of magnitude below the uplink value on the same layers, yet still far above the floor $1/d = 4.2\times10^{-7}$. The dimension dependence of \hyperref[cor:smooth]{Corollary~\ref*{cor:smooth}} is therefore not attained on the downlink either, but the downlink lacks the $\Theta(1)$ contraction that suppresses it on the uplink.
\end{remark}

\begin{remark}[\textbf{The scaled sign is not layer-norm contractive}]\label{rem:normmismatch}
\hyperref[cor:smooth]{Corollary~\ref*{cor:smooth}} reaches EF21-MuonSign only through \citet[Rem.~23]{gruntkowska2025error} and its factor $\bar{\rho}_i^2$, for two reasons. First, the alternative, a server compressor contractive in the \emph{layer} norm, is unavailable: the scaled sign is contractive there for no $\alpha_P>0$; for $\mathbf{Y}_\delta = (1-\delta)\mathbf{I}_n + \delta\mathbf{J}_n$ ($\mathbf{J}_n$ all ones, $0<\delta<1$, $n \ge 2$), every entry is positive, so $\operatorname{sign}(\mathbf{Y}_\delta) = \mathbf{J}_n$, $\operatorname{mean}|\mathbf{Y}_\delta| = (1+(n-1)\delta)/n$, and
\begin{equation*}
\mathcal{C}(\mathbf{Y}_\delta) - \mathbf{Y}_\delta = (1-\delta)\bigl(\tfrac{1}{n}\mathbf{J}_n - \mathbf{I}_n\bigr),
\end{equation*}
a matrix with eigenvalues $0$ and $-(1-\delta)$; hence $\|\mathcal{C}(\mathbf{Y}_\delta)-\mathbf{Y}_\delta\|_{2\to2} = 1-\delta$ while $\|\mathbf{Y}_\delta\|_{2\to2} = 1+(n-1)\delta$, and the ratio tends to $1$ as $\delta\downarrow 0$. Second, the sufficient condition of \citet[App.~D]{gruntkowska2025error}, which certifies a compressor in $\mathbb{B}_2(\alpha)$ as layer-norm contractive when $\alpha > 1-1/r_i$, is out of reach: for a $768\times3072$ layer it demands $\alpha > 0.9987$, where the scaled sign attains $2/\pi$ on dense inputs (\hyperref[lem:signcontr]{Lemma~\ref*{lem:signcontr}}). The factor $\bar{\rho}_i^2 = r_i$ in $\zeta_i$ can therefore be removed only by replacing the compressor, for instance by random dropout ($\alpha_P = p$) or a Top-$K$ SVD compressor ($\alpha_P = 1-\sigma_{K+1}^2/\sigma_1^2$), both layer-norm contractive \citep[App.~D]{gruntkowska2025error} and neither one-bit.
\end{remark}

\begin{remark}[\textbf{From Muon to Gluon}]\label{rem:gluon}
Only $(\underline{\rho}_i,\bar{\rho}_i) = (1,\sqrt{r_i})$ is spectral-norm-specific: \hyperref[lem:signcontr]{Lemma~\ref*{lem:signcontr}} is Euclidean and the framework's theorems hold for arbitrary layer norms. A new geometry has to supply only its LMO and its norm-equivalence pair; both corollaries then hold with those constants in place of $(1,\sqrt{r_i})$, giving \emph{EF21-GluonUSign} and \emph{EF21-GluonSign} for the Gluon setting \citep{riabinin2025gluon}. Admissible geometries abound: $\ell_1\to\ell_\infty$ for embeddings \citep{pethick2025training}, the Schatten-$p$ norms \citep{cesista2025schattenp}, the Ky Fan duals of the Fanion family \citep{kravatskiy2025kyfannorms}, and, by the closure property of the last work, the norms whose LMO is a conic combination of these LMOs; their unit ball is the corresponding Minkowski sum, whence $\bar{\rho}_i \le \sum_j\alpha_j\bar{\rho}_i^{(j)}$. Convergence is thus a property of error feedback together with scaled-sign compression, not of the spectral geometry.
\end{remark}

We assume the exact spectral LMO, as do all analyses of Muon-type methods \citep{li2025note,kovalev2025understanding,riabinin2025gluon,gruntkowska2025error}; in practice the polar factor is approximated by polynomial iterations \citep{amsel2025polar,grishina2025chebyshev}, in our case the five Newton--Schulz steps of \hyperref[alg:muon_lmo]{Algorithm~\ref*{alg:muon_lmo}}, with vector parameters and the last layer trained by AdamW as usual \citep{jordan2024muon}. That substitution is not covered by \hyperref[thm:conv]{Theorem~\ref*{thm:conv}} either; \citet{shulgin2025inexactmuon} analyse the inexact Muon update directly and find the method tolerant of oracle error.
%%% <<< end ef21_musign_reduction.tex

\subsection{Reproducibility Details}\label{app:repro}
\paragraph{Choice of benchmarks.}
Each benchmark answers a question the others cannot. On the convex quadratic of \hyperref[app:images_task]{Appendix~\ref*{app:images_task}} the smoothness constant and the minimizer are known in closed form, so the alignment the counterexamples attack can be measured there rather than inferred. CIFAR-10 with a ResNet-18 \citep{krizhevsky2009cifar} carries no such special property: we use it because it is the benchmark on which the sign-compression line of work is quoted \citep{bernstein2018signsgd,bernstein2019signsgdmajority,karimireddy2019efsign} and on which SignMuon's own concurrent proposal is evaluated \citep{mishra2026signmuon}; it is small enough to run every method at several step sizes and several seeds, which is what the claims about seed spread require. The federated split of the same data at $N=11$ tests the methods in the setting they are designed for, a bandwidth-limited link between clients and a server. NanoGPT supplies what CIFAR cannot: a transformer language model at practical scale, with matrices wide enough for the layer-rank term of \hyperref[cor:smooth]{Corollary~\ref*{cor:smooth}} to be visible, which is where the two models of EF21-MuonSign separate.

\paragraph{Computing infrastructure.}
The experiments were not all run on the same machine. The synthetic study ran on one NVIDIA RTX A4000 (128-core AMD EPYC 7543 host, $472$\,GB RAM; Linux 6.12, Python 3.12, PyTorch 2.7.0, CUDA 12.8, driver 575.51.03); the language-modelling runs used a rented $8\times$H100 SXM node, specified in full in \hyperref[app:nanogpt]{Appendix~\ref*{app:nanogpt}}. The centralized and federated CIFAR-10 runs were executed on single-GPU workstations; each run records its machine, commit and wall time in its \texttt{metrics.json}. All $126$ centralized runs and all $175$ federated runs behind the tables and figures below were executed on the same machine and at one commit: an NVIDIA RTX A4500 ($19.6$\,GB) in a $32$-core AMD Ryzen 9 5950X host with $62.7$\,GB of memory, under Linux 5.15, Python 3.12.11, PyTorch 2.5.1+cu124, CUDA 12.4, driver 560.35.05.

\paragraph{Randomness and seeds.}
The network experiments are seeded through a single routine that seeds Python's \texttt{random} module, NumPy and PyTorch on all CUDA devices; the federated runs additionally pin cuDNN to deterministic kernels, while the centralized sweep leaves cuDNN autotuning enabled, since it reports the spread across seeds rather than a bitwise-reproducible trajectory. The synthetic study forks and re-seeds its own generator per configuration. The federated experiment of \hyperref[tab:exp_3]{Table~\ref*{tab:exp_3}} uses five seeds ($0$--$4$) per method and the centralized experiment of \hyperref[tab:cifar_main]{Table~\ref*{tab:cifar_main}} three ($0$--$2$), each reported as mean $\pm$ one sample standard deviation across seeds; the weight-decay ablations use seed $0$. The nanoGPT runs of \hyperref[tab:nanogpt]{Table~\ref*{tab:nanogpt}} are single runs at the speedrun's own unpinned initialization, and the table quotes the five-seed spread published upstream; our released script accepts an explicit seed but pins the generator only, since deterministic kernels would forfeit the wall-clock time the same table reports. The synthetic study of \hyperref[app:images_task]{Appendix~\ref*{app:images_task}} averages over three draws of the problem (seeds $1337$--$1339$) at a fixed $\mathbf{X}_0$ (seed $42$), its claims concerning random instances. Differences smaller than the seed spread are not claimed as results, and we report spreads rather than significance tests because at these seed counts no test could reject: a paired Wilcoxon signed-rank test over $n$ seeds has smallest attainable two-sided exact $p$-value $2^{1-n}$, which is $0.0625$ at five seeds and $0.25$ at three, above the $5\%$ level in both cases. The comparisons we claim are separated by several standard deviations. Learning-rate selection is performed once, at seed $0$, on a validation split disjoint from the test set, and the selected rate is reused unchanged for every seed.

\paragraph{Step-size schedules.}
Each experiment inherits the schedule conventional to its domain, and none of them is the constant rate our rates are proved for. The ResNet runs anneal cosinally to zero, the standard for this architecture and the premise of every accuracy we can be compared against; the nanoGPT runs keep record \#40's stable-then-decay schedule, flat for the first $55\%$ of steps and decaying linearly to a tenth of the base rate $\eta_0$ of \eqref{eq:unit_gain_main} thereafter (\hyperref[app:nanogpt]{Appendix~\ref*{app:nanogpt}}), because altering it would forfeit the reproduction that validates our port. The discrepancy is deliberate: matching the analysed step size would sacrifice comparability on both benchmarks and close only one of several gaps between the theory and the runs, the others being Newton--Schulz in place of an exact oracle, momentum, and normalization layers. One measurement is exempt: the growth-exponent diagnostic of \hyperref[app:lrscale]{Appendix~\ref*{app:lrscale}}, which is run at a constant rate because under a decaying one the accumulated update saturates and the fit reports the schedule instead of the alignment.

\paragraph{Learning-rate selection.}
The only tuned hyperparameter is $\eta_0$. Every method with a norm-fixed step is tuned and reported under the unit-gain rule~\eqref{eq:unit_gain_main}, so that $\eta_0$ is the per-step RMS gain for each of them; SGD and Adam have no norm-fixed step and run at one global rate. The rule is a heuristic, and \hyperref[app:lrscale]{Appendix~\ref*{app:lrscale}} bounds and measures what depends on it: the three \textsc{sign} methods keep their order when re-tuned from scratch under one global rate and under $\mu$P (\hyperref[tab:rule_ablation]{Table~\ref*{tab:rule_ablation}}). Momentum is fixed at $0.9$ and weight decay at $0$ in the primary tables, the setting the sweep of \citet{mishra2026signmuon} itself selects; the regularized case is an ablation (\hyperref[app:images_cen]{Appendix~\ref*{app:images_cen}} centralized, and below for the federated study). The auxiliary group, biases, normalization parameters and the classifier head, is trained by AdamW at $10^{-3}$ for every method, a convention rather than a verified common optimum: a sweep at matched budget places SignMuon's optimum at $10^{-3}$ and Muon's at $2\times10^{-3}$, a difference of $0.16$ points that lies within the seed spread, so the auxiliary rate is method-dependent, to a degree the sweep does not quantify.

In the centralized study, selection uses a fixed $45$k/$5$k train/validation partition and validation accuracy averaged over the last five epochs; the test set is never consulted during tuning. Each method starts from the same five-point $1$--$2$--$5$ lattice (three points per order of magnitude); an optimum at a grid endpoint triggers a widening and a re-run, up to four times, extending Muon and EF21-MuonUSign to seven points and SignSGD to nine. Selection runs use the same $75$-epoch cosine schedule as the reported ones, so the tuning and reporting horizons coincide. The selected $\eta_0$ is then retrained on the full $50$k training set at three seeds, and we report the mean and standard deviation of the test accuracy over the last five epochs.

\paragraph{Tuning the federated study.}
So that no placement is handicapped by its step size, every method receives the same tuning budget: a five-point $1$--$2$--$5$ lattice in $\eta_0$, ranked on a $5$k validation split held out of the $50$k \emph{before} the client partition, at the full $2000$-round horizon \hyperref[tab:exp_3]{Table~\ref*{tab:exp_3}} reports, with the grid widened and the method re-tuned whenever an optimum occurred at an endpoint. SignMuon is the one method that required the widening, settling at $\eta_0=0.1$ on seven points; every selected rate is interior to its own grid. The reported runs then use the full $50$k at the selected rate, so no test image is ever scored during selection.

One selection margin requires comment. Eight of the nine methods carrying a per-layer multiplier (the eight of \hyperref[app:lrscale]{Appendix~\ref*{app:lrscale}} and the server-side-LMO control) separate their selected rate from the runner-up by $0.15$ to $1.07$ validation points; SignMuon separates $0.1$ from $0.05$ by $0.02$ at the single tuning seed, which is no separation at all, so its row in \hyperref[tab:exp_3]{Table~\ref*{tab:exp_3}} is to be read as either of two adjacent lattice points. The ambiguity does not extend further, the next points out lying $0.5$ and $1.3$ points behind.

The selected rates span a factor of ten, from $\eta_0 = 0.01$ for SignSGD and EF21-MuonUSign to $0.1$ for Muon and SignMuon, with both families covering that range. This spread is between \emph{methods} and carries no verdict on the per-layer rule, whose claim concerns layer \emph{shape}; the measurement that does bear on the rule, re-tuning under competing conventions, is \hyperref[tab:rule_ablation]{Table~\ref*{tab:rule_ablation}}. Weight decay is $0$ in the reported table. Switching on a decoupled $5\times10^{-4}$ at seed $0$ moves Muon by $-0.32$ points, SignMuon by $+0.27$ and the server-side-LMO control by $+0.05$: Muon leads in either setting, and the other two exchange places by margins below the $0.24$ standard deviation SignMuon carries over five seeds, so the ablation separates nothing that the primary table does not.

\paragraph{Accuracy and the threshold column.}
Alongside final accuracy we report the number of epochs (rounds, in the federated tables) to reach a fixed test accuracy. The two measure distinct quantities: with the methods spanning about a point and a half at $75$ epochs (\hyperref[tab:cifar_central]{Table~\ref*{tab:cifar_central}}), final accuracy is close to the noise floor, whereas a threshold crossing on a monotonically rising curve separates the methods by factors rather than by tenths of a point, and is the analogue of the ``steps to $3.35$'' column of \hyperref[tab:nanogpt]{Table~\ref*{tab:nanogpt}}.

\paragraph{Conventions with numerical consequences.}
Three implementation conventions can displace reported numbers and are recorded here. (i)~The Muon LMO is computed in \texttt{bfloat16} (five Newton--Schulz steps) unless stated otherwise; for the methods that sign the LMO \emph{output}, entries of $\operatorname{polar}(\cdot)$ near zero may flip at this precision, so that their trajectories carry a precision-dependent component (\texttt{--lmo-dtype float32} is available). (ii)~In the federated runs, BatchNorm running statistics are never updated: local models are discarded each round and BN runs in inference mode during gradient accumulation, so the statistics stay at their initialization for the entire run, in training and evaluation alike. The result is a fixed normalization with learnable affine parameters, self-consistent between train and test, applied identically to every method. It is also one reason a channel may remain inactive across an entire local batch and so contribute an exactly zero row to the momentum. (iii)~For EF21-MuonSign, training metrics are logged at the broadcast model $\mathbf{W}$ (where gradients must be evaluated) while validation and test metrics are evaluated at the server model $\mathbf{X}$ of \eqref{eq:two_models}, the iterate the guarantee bounds, except where a table states otherwise.

\subsection{The smooth convex problem}\label{app:images_task}

% Moved here from the main text (Experiments) to save main-part space. Every
% number below is read off results/synthetic/SUMMARY.md, run 2026-07-29.
On a deterministic $L$-smooth convex quadratic we measure the scalar the guarantees rest on: the alignment $\rho_t = \langle \nabla F(\mathbf{X}_t), \mathbf{D}_t\rangle / (\|\nabla F(\mathbf{X}_t)\|_F\|\mathbf{D}_t\|_F)$ between the gradient and the step. \hyperref[th:1]{Theorems~\ref*{th:1}}--\ref{th:3} construct instances driving it negative; on random instances the three methods they cover keep $\rho_t \ge 0.142$ throughout a tuned trajectory, so the construction is not one that random data reproduces. Only EF21-MuonSign becomes negative, on $0.95\%$ of steps and to $-0.026$; it is also the one method here possessing a convergence guarantee (\hyperref[app:hyp]{Appendix~\ref*{app:hyp}}), obtained without per-step descent.

The same experiment separates the two effects that a fixed-target iteration count confounds. Sign compression of an \textsc{lmo} step secures a lower accuracy floor, not a faster rate: SignMuon's floor lies a factor $1.74$--$1.80$ below SignSGD's at every step size, and the two share $\|\mathbf{S}\|_F$ exactly, so the separation resides in the floor rather than in the step length.

\paragraph{Construction.}
To isolate the effect of matrix structure from stochastic noise and from the complexity of DNN architectures, we use a deterministic $L$-smooth convex quadratic:
% One line, not a multline: this appendix is only ever set \onecolumn (both
% v2_SignMuon_AAAI_supp and _full switch before \input-ing it), so the break
% the two-column measure needed only opened a gap across the page.
\begin{equation}
F(\mathbf{X}) = \frac{1}{2} \|\mathbf{A}^{1/2} \mathbf{X} \mathbf{B}^{1/2}\|_F^2
= \frac{1}{2} \langle \mathbf{X}, \mathbf{A} \mathbf{X} \mathbf{B} \rangle
\to \min_{\mathbf{X} \in \mathbb{R}^{m \times n}},
\label{eq:synthetic_problem}
\end{equation}
with $\mathbf{A} \in \mathbb{S}_{++}^m$, $\mathbf{B} \in \mathbb{S}_{++}^n$ symmetric and $\mathbf{X}_0$ drawn entrywise from $\mathcal{N}(0, 0.01)$. Eigenvalues are sampled uniformly from $(0,1)$ in a Haar-random eigenbasis, so the matrices are almost surely positive \emph{definite} and the minimizer is $\mathbf{X}^\star = \mathbf{0}$ with $F^\star = 0$.

Two facts about this instance are exact rather than estimated, and both are used below. The Hessian of $F$ is the Kronecker product $\mathbf{B} \otimes \mathbf{A}$, so its eigenvalues are the products $\lambda_i(\mathbf{A})\lambda_j(\mathbf{B})$: the Frobenius smoothness constant is $L = \max_{ij}\lambda_i(\mathbf{A})\lambda_j(\mathbf{B}) \le 1$ and the strong-convexity constant is $\sigma = \min_{ij}\lambda_i(\mathbf{A})\lambda_j(\mathbf{B}) > 0$. The uniform draw leaves the resulting condition number $L/\sigma$ uncontrolled (it is near $3.7 \times 10^{4}$ at the $m = n = 100$ every measurement in this subsection uses), so where conditioning is the variable we instead use log-spaced spectra with $L = 1$ and $L/\sigma$ set exactly. And since $\nabla F(\mathbf{X}) = \mathbf{A}\mathbf{X}\mathbf{B}$ in closed form, the gradient can be evaluated at any point without an autograd graph, which permits the bidirectional method to be scored on its exact model $\mathbf{X}$ while its gradient is taken at the broadcast model $\mathbf{W}$, as its algorithm requires.

\paragraph{The fixed-target criterion.}
The natural criterion, fewest iterations to $F(\mathbf{X}) \le 10^{-3}$ within $T_{\max} = 5000$ with learning rate and momentum tuned per method, does not measure a convergence rate. Eight of the ten methods take a \emph{norm-fixed} step, $\|\operatorname{sign}(\cdot)\|_F = \sqrt{mn}$ and $\|\mathrm{polar}(\cdot)\|_F = \sqrt{r}$, so at a constant $\eta$ the iterate settles into a ball of radius $\eta\|\mathbf{S}\|_F$ and $F$ plateaus; Adam, bounded entrywise by ${\approx}\eta$, plateaus as well, and SGD, whose step vanishes with the gradient, is the only method that does not. Write $F_\infty$ and $g_\infty$ for the settled values of $F(\mathbf{X}_t)$ and $\|\nabla F(\mathbf{X}_t)\|_F$. Measured directly, $g_\infty \propto \eta$ for every method possessing a floor, $F_\infty \propto \eta^2$ for every such method but SignSGD, whose exponent is $1.33$, and the iteration count is $\text{const}/\eta$. The tuner accordingly returns the largest $\eta$ whose plateau falls under the target, and the resulting ranking is one of \emph{accuracy floors}. Over the seven step sizes at which both methods were run, SignMuon's $g_\infty$ lies a factor $1.74$ to $1.80$ below SignSGD's, flat in $\eta$ as two floors of equal exponent must be, while its $F_\infty$ lies $6.7$ to $71$ times below, the two $F$-exponents differing. Tuned per budget, SignMuon holds the smaller $\min_t\|\nabla F\|_F$ at every horizon but $T = 250$, where SignSGD prevails by $2\%$, a margin below what three draws resolve. We therefore report floor and rate separately, reading the descent lemma
\begin{equation}
F(\mathbf{X}_{t+1}) \le F(\mathbf{X}_t) - \eta \langle \nabla F(\mathbf{X}_t), \mathbf{D}_t\rangle + \tfrac{\eta^2 L}{2}\|\mathbf{D}_t\|_F^2
\label{eq:descent_lemma}
\end{equation}
as the statement that separates them: the second term is the floor, the first is the rate.

\begin{figure*}[!t]
    \centering
    \includegraphics[width=\textwidth]{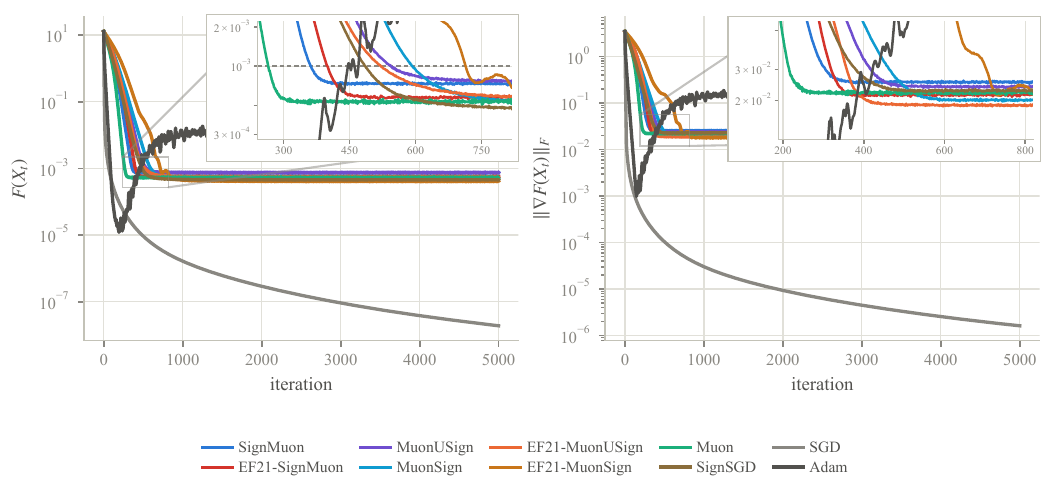}
    \caption{Trajectories at the optima of \hyperref[tab:synthetic_tuned]{Table~\ref*{tab:synthetic_tuned}}, geometric mean over three problem draws. Every normalized method plateaus within a few hundred iterations and stays there, which renders the crossing time a ranking of floors; SGD, whose step vanishes with the gradient, is the only curve still descending at $T_{\max}$, and Adam oscillates about its floor rather than settling onto it. The insets magnify the arrival window, over the full axis a fold of curves in the first sixth of the range, and are scaled to the eight norm-fixed methods: the dashed line is the target $F \le 10^{-3}$, and the order in which the curves cut it is the iteration count of \hyperref[tab:synthetic_tuned]{Table~\ref*{tab:synthetic_tuned}}, Muon first and EF21-MuonSign last. Those counts are means of the per-draw crossing times rather than crossings of the plotted mean, and for the slowest method the two differ by enough to see.}
    \label{fig:synthetic_main}
\end{figure*}

\hyperref[tab:synthetic_tuned]{Table~\ref*{tab:synthetic_tuned}} nonetheless reports the criterion. Read as a ranking of floors, it orders the six placements identically within both families: sign after the \textsc{lmo} is the least expensive, sign on both channels the most. Error feedback adds $10\%$ on SignMuon and $17\%$ on MuonSign, and removes $6\%$ on MuonUSign, the one placement whose \textsc{lmo} already receives a compressed argument. Muon surpasses all six placements, and SGD ($66$) and Adam ($85$) surpass every normalized step by factors of three to ten, a quadratic with an exactly known Hessian being precisely the case for which a scaled gradient step is designed. \hyperref[fig:synthetic_main]{Figure~\ref*{fig:synthetic_main}} plots the trajectories behind those counts, and shows the plateau that makes the criterion a ranking of floors.

\begin{table}[!t]
\centering\small
\setlength{\tabcolsep}{4pt}
\begin{tabular}{@{}lrrrl@{}}
\toprule
\textbf{Algorithm} & \textbf{iters} & \textbf{best $F$} & \textbf{$\min\|\nabla F\|_F$} & \textbf{tuned $(\eta,\mu)$} \\
\midrule
SignMuon       & $364$ & $6.9\cdot10^{-4}$ & $2.4\cdot10^{-2}$ & $(1.0\cdot10^{-3},\, 0)$ \\
MuonUSign      & $556$ & $7.1\cdot10^{-4}$ & $2.3\cdot10^{-2}$ & $(1.0\cdot10^{-2},\, 0.5)$ \\
MuonSign       & $595$ & $4.8\cdot10^{-4}$ & $1.9\cdot10^{-2}$ & $(6.8\cdot10^{-4},\, 0.5)$ \\
EF21-SignMuon  & $401$ & $5.4\cdot10^{-4}$ & $2.1\cdot10^{-2}$ & $(6.8\cdot10^{-3},\, 0)$ \\
EF21-MuonUSign & $521$ & $5.3\cdot10^{-4}$ & $1.8\cdot10^{-2}$ & $(6.8\cdot10^{-3},\, 0)$ \\
EF21-MuonSign  & $695$ & $3.8\cdot10^{-4}$ & $1.7\cdot10^{-2}$ & $(4.6\cdot10^{-3},\, 0)$ \\
Muon           & $267$ & $4.7\cdot10^{-4}$ & $2.0\cdot10^{-2}$ & $(1.0\cdot10^{-2},\, 0.5)$ \\
SignSGD        & $486$ & $4.3\cdot10^{-4}$ & $2.2\cdot10^{-2}$ & $(6.8\cdot10^{-4},\, 0.5)$ \\
SGD            & $66$  & $1.9\cdot10^{-8}$ & $1.6\cdot10^{-6}$ & $(2.15,\, 0.5)$ \\
Adam           & $85$  & $9.2\cdot10^{-6}$ & $6.5\cdot10^{-4}$ & $(6.8\cdot10^{-2},\, \text{---})$ \\
\bottomrule
\end{tabular}
\caption{Fixed-target criterion under the current protocol: fewest iterations to $F \le 10^{-3}$ within $T_{\max} = 5000$ at $m = n = 100$, over three problem draws, $(\eta,\mu,\text{schedule})$ tuned on logarithmic grids spanning five orders of magnitude. A constant schedule is selected for every method. Best $F$ and $\min\|\nabla F\|_F$ are minima over all $5000$ iterations rather than values at the crossing, so they report the floor the method settles into.}
\label{tab:synthetic_tuned}
\end{table}

\paragraph{Alignment.}
Equation~\eqref{eq:descent_lemma} makes progress contingent on a single scalar, the normalized alignment between the gradient and the step actually taken,
\begin{equation}
\rho_t := \frac{\langle \nabla F(\mathbf{X}_t), \mathbf{D}_t\rangle}{\|\nabla F(\mathbf{X}_t)\|_F \, \|\mathbf{D}_t\|_F} \in [-1, 1].
\label{eq:alignment}
\end{equation}
\hyperref[th:1]{Theorems~\ref*{th:1}}--\ref{th:3} are constructions that drive $\rho_t$ negative. \hyperref[tab:synthetic_alignment]{Table~\ref*{tab:synthetic_alignment}} reports its distribution along the tuned trajectory on random instances, which is the empirical counterpart of those theorems and the one measurement here that is about the methods rather than about the tuning protocol. Three closed forms anchor it: $\rho = 1$ for SGD, $\rho = \|\mathbf{G}\|_1 / (\|\mathbf{G}\|_F\sqrt{mn})$ for SignSGD, and $\rho = \|\mathbf{G}\|_* / (\|\mathbf{G}\|_F\sqrt{r})$ for Muon. The six sign-around-the-\textsc{lmo} methods admit none, which is the subject of this paper.

\begin{table*}[!t]
\centering
\small
\begin{tabular}{@{}lrrrrrl@{}}
\toprule
\textbf{Algorithm} & \textbf{$\min_t \rho_t$} & \textbf{median $\rho_t$} & \textbf{mean $\rho_t$} & \textbf{\% of steps with $\rho_t<0$} & \textbf{closed form} & \textbf{tuned $(\eta,\mu)$} \\
\midrule
SignMuon        & $0.208$  & $0.381$ & $0.399$ & $0.00\%$ & ---      & $(2.2\cdot10^{-4},\, 0)$ \\
MuonUSign       & $0.161$  & $0.357$ & $0.378$ & $0.00\%$ & ---      & $(3.2\cdot10^{-3},\, 0)$ \\
MuonSign        & $0.142$  & $0.360$ & $0.364$ & $0.00\%$ & ---      & $(2.2\cdot10^{-4},\, 0)$ \\
EF21-SignMuon   & $0.333$  & $0.659$ & $0.573$ & $0.00\%$ & ---      & $(1.5\cdot10^{-3},\, 0)$ \\
EF21-MuonUSign  & $0.346$  & $0.558$ & $0.553$ & $0.00\%$ & ---      & $(1.5\cdot10^{-3},\, 0)$ \\
EF21-MuonSign   & $-0.026$ & $0.489$ & $0.419$ & $0.95\%$ & ---      & $(1.5\cdot10^{-3},\, 0)$ \\
Muon            & $0.454$  & $0.666$ & $0.633$ & $0.00\%$ & $0.695$  & $(1.5\cdot10^{-3},\, 0)$ \\
SignSGD         & $0.174$  & $0.503$ & $0.499$ & $0.00\%$ & $0.794$  & $(2.2\cdot10^{-4},\, 0.5)$ \\
SGD             & ---      & ---     & ---     & ---      & $1$      & $(3.2,\, 0.95)$ \\
\bottomrule
\end{tabular}
\caption{Alignment $\rho_t$ of Equation~\eqref{eq:alignment} along the tuned trajectory, over three problem draws at $m = n = 100$. The closed form is evaluated at $\mathbf{X}_0$ without momentum and is comparable only to rows whose tuned $\mu$ is $0$; the six sign-around-the-\textsc{lmo} methods admit none. SGD's closed form is $1$ but its step is not instrumented, so its distribution column is empty; Adam is omitted, its step lying outside the descent lemma. Every method of \hyperref[th:1]{Theorems~\ref*{th:1}}--\ref{th:3} stays bounded away from zero; the single negative excursion is EF21-MuonSign's, which \hyperref[app:hyp]{Appendix~\ref*{app:hyp}} proves convergent through the EF21 estimator rather than per-step descent.}
\label{tab:synthetic_alignment}
\end{table*}

\paragraph{Closed-form checks.}
Three measurements test quantities the theory predicts in closed form. (i)~The floor: balancing the two terms of \eqref{eq:descent_lemma} at $\langle\nabla F,\mathbf{D}\rangle = \rho\|\nabla F\|_F\|\mathbf{D}\|_F$ gives $g_\infty = \eta\,L\|\mathbf{S}\|_F/(2\rho)$, a slope of $1$ in $\log\eta$ with that coefficient. SignMuon and SignSGD share $\|\mathbf{S}\|_F = \sqrt{mn}$ exactly, so any separation between their floors is attributable to $\rho$ alone. (ii)~The budget exponent: tuning $(\eta, \mu, \text{schedule})$ separately at each horizon $T$ and fitting $\text{err} \propto T^{-p}$, $\eta^\star \propto T^{-q}$, with $\text{err} := \min_{t\le T}\|\nabla F(\mathbf{X}_t)\|_*^2$ the squared dual norm our theorems bound. The nonconvex bound gives $p = q = \tfrac12$; a strongly convex problem gives $p = q = 1$, and this instance, having $\sigma > 0$, need not occupy the nonconvex regime. The fit reports both at once: the eight norm-fixed methods tune $\eta^\star$ as the nonconvex bound prescribes, $q$ scattering about $\tfrac12$, while the attained error falls at $p$ between $1.76$ and $2.14$, a step size chosen for the worst case of the smoothness class applied to an instance far easier than that worst case. (iii)~The stability edge: the largest stable $\eta$, with SGD as a control required to reproduce the textbook $2/L$. Reported as the step length $\eta_{\max}\|\mathbf{S}\|_F$, this would be family-independent were the operative bound the Frobenius ball; the spread measures how far that bound stands from the geometry in which each step actually resides.

\begin{table*}[!t]
\centering
\small
\begin{tabular}{@{}lrrrrrrr@{}}
\toprule
& \multicolumn{2}{c}{\textbf{floor}} & \multicolumn{4}{c}{\textbf{budget}} & \textbf{stability} \\
\cmidrule(lr){2-3}\cmidrule(lr){4-7}\cmidrule(lr){8-8}
\textbf{Algorithm} & $\mathrm{d}\log g_\infty / \mathrm{d}\log\eta$ & $R^2$ & $p$ & $R^2$ & $q$ & $R^2$ & $\eta_{\max}\|\mathbf{S}\|_F$ \\
\midrule
SignMuon        & $1.000$ & $1.000$ & $1.88$ & $0.996$ & $0.39$  & $0.94$ & $13.4$ \\
MuonUSign       & $1.000$ & $1.000$ & $1.86$ & $0.995$ & $0.39$  & $0.94$ & $14.9$ \\
MuonSign        & $1.000$ & $1.000$ & $1.93$ & $0.993$ & $0.44$  & $0.94$ & $11.8$ \\
EF21-SignMuon   & $1.000$ & $1.000$ & $2.07$ & $1.000$ & $0.55$  & $1.00$ & $18.9$ \\
EF21-MuonUSign  & $1.000$ & $1.000$ & $2.04$ & $1.000$ & $0.55$  & $1.00$ & $14.9$ \\
EF21-MuonSign   & $1.000$ & $1.000$ & $2.14$ & $1.000$ & $0.55$  & $1.00$ & $14.9$ \\
Muon            & $1.000$ & $1.000$ & $2.03$ & $0.999$ & $0.55$  & $1.00$ & $17.6$ \\
SignSGD         & $0.989$ & $1.000$ & $1.76$ & $0.990$ & $0.39$  & $0.94$ & $11.1$ \\
SGD             & \multicolumn{2}{c}{\emph{no floor}} & $8.97$ & $0.985$ & $-0.11$ & $0.50$ & $2.06^{\dagger}$ \\
Adam            & $1.000$ & $1.000$ & $2.65$ & $0.915$ & $-0.66$ & $0.62$ & $> 100^{\ddagger}$ \\
\midrule
\emph{predicted} & $1$ & & $1/2$ or $1$ & & $1/2$ or $1$ & & $2/L$ \\
\bottomrule
\end{tabular}
\caption{Floor exponent, budget exponents and stability edge against the predicted values, at $m = n = 100$ over three problem draws. The floor exponent attains its prediction exactly for eight of the nine methods possessing a floor and to within $1.1\%$ for SignSGD, fitted on the seven step sizes at which the plateau is reached (four for Adam). SGD possesses no floor, its step vanishing with the gradient, so its floor columns are empty and its $\eta_{\max}$ is the control reproducing $2/L$; neither baseline admits a power law in $\eta^\star$. $^\dagger$~SGD has no $\|\mathbf{S}\|_F$; the entry is $\eta_{\max}$, $1.02$ times $2/L$. $^\ddagger$~The search reached its ceiling without encountering an edge: Adam's step is bounded by ${\approx}\eta$ irrespective of the gradient, so it oscillates rather than diverges.}
\label{tab:synthetic_dynamics}
\end{table*}

\begin{figure*}[!t]
    \centering
    \includegraphics[width=\textwidth]{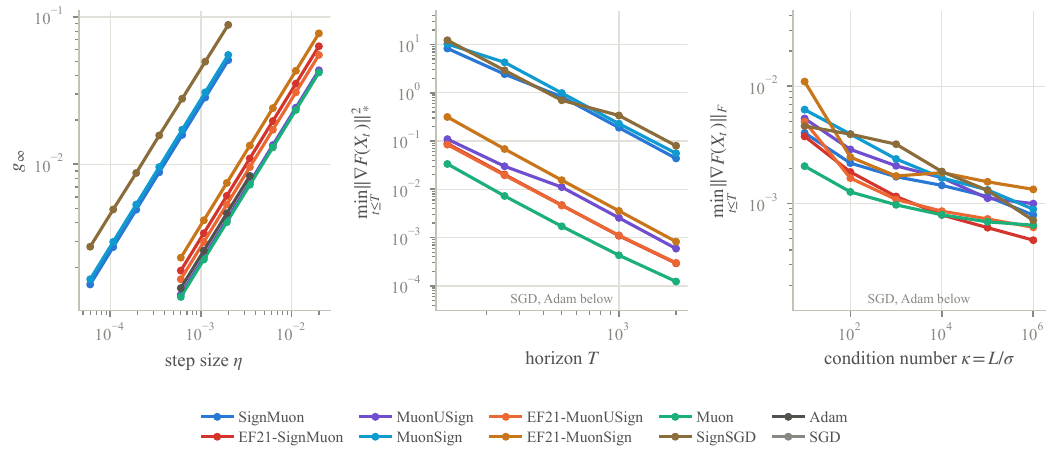}
    \caption{Log--log fits behind \hyperref[tab:synthetic_dynamics]{Table~\ref*{tab:synthetic_dynamics}}, at $m = n = 100$ over three problem draws. \emph{Left}, the accuracy floor $g_\infty$ of a constant step, the settled $\|\nabla F\|_F$: the lines are parallel because every exponent is $1$, and the two clusters are the two step lengths, the three methods with $\|\mathbf{S}\|_F = \sqrt{mn}$ on the left and the five with $\sqrt{r}$ on the right, each swept over the window its length calls for. \emph{Centre}, error at each retuned horizon, in the norm dual to each method's own \textsc{lmo} ball, whence the two bands: $\ell_1$ above and nuclear below, comparable within a family but not across families. \emph{Right}, tuned at each $\kappa$; over the five orders of magnitude swept, the eight norm-fixed methods vary by a factor of $3.2$ to $8.3$, against the sixteen orders of magnitude spanned by SGD and Adam. Both of the last two panels are scaled to those eight, leaving SGD and Adam below the frame: on a quadratic with an exactly known Hessian a scaled gradient step converges to machine precision. The gap reaches nineteen orders of magnitude at $\kappa = 10$, where SGD attains $2\cdot10^{-22}$ against the $2\cdot10^{-3}$ of the closest norm-fixed method.}
    \label{fig:synthetic_dynamics}
\end{figure*}

\paragraph{Conditioning.}
Conditioning governs the dynamics of a quadratic, and the construction above leaves it to chance, so the right panel of \hyperref[fig:synthetic_dynamics]{Figure~\ref*{fig:synthetic_dynamics}} sweeps $\kappa$ over five orders of magnitude, from $10$ to $10^6$, at fixed $L = 1$. The eight norm-fixed methods are insensitive to it: over the entire sweep the attained $\|\nabla F\|$ moves by a factor of $3.2$ (Muon) to $8.3$ (EF21-MuonSign), and the fitted $\mathrm{d}\log\|\nabla F\|/\mathrm{d}\log\kappa$ lies between $-0.096$ (Muon) and $-0.17$ (EF21-SignMuon). The factor is the measurement and the slope a summary of it, which is the order to read them in: six of the eight fits have $R^2 \ge 0.90$, but EF21-MuonUSign's is $0.83$ and EF21-MuonSign's $0.67$, so for those two rows the slope is not itself a quantity we would quote. The sign is mildly negative because at fixed $L$ a larger $\kappa$ entails a smaller $\sigma$, hence a flatter landscape to occupy rather than a harder one to descend; their floors are fixed by $\eta\|\mathbf{S}\|_F$, which is independent of the spectrum. SGD and Adam are sensitive to it, jointly spanning sixteen orders of magnitude across the sweep, which is why the panel is scaled to the eight and leaves those two below it. SGD's fitted slope of $3.6$ is not an exponent: its first three points lie at or below $10^{-20}$, which on this problem constitutes exact convergence rather than a measurement.

\paragraph{Protocol and reproduction.}
Every measurement above is at $m = n = 100$ over three problem draws (seeds $1337$--$1339$, $\mathbf{X}_0$ from seed $42$), with the \textsc{lmo} taken in \texttt{bfloat16} as in the network experiments, and is reported as the geometric mean over the draws for the error metrics and the arithmetic mean for iteration counts.

Learning rates are searched on logarithmic grids spanning five orders of magnitude, one per step-norm family, each ending past the largest stability edge measured for that family so that no stable step size falls outside the search. A linear grid spanning a single order of magnitude, which an earlier version of this experiment used, is narrow enough to censor an optimum at a grid boundary, giving an upper bound rather than a tuned value; the search flags any optimum at an edge. No learning rate reported above is flagged. The windows are derived from $\|\mathbf{S}\|_F$ rather than from the method name, which matters for MuonUSign and EF21-SignMuon: both take a step of length $\sqrt{r}$ despite the \textsc{sign} in their names, and under a window assigned by name both were searched an order of magnitude below the range their step length calls for.

The remaining flags are on momentum, at the top of its grid ($0.99$), and belong to SGD and SignSGD at the largest condition numbers: an ill-conditioned quadratic asks for heavy momentum, so those two rows of the $\kappa$ sweep bound the dependence from one side rather than measuring it. The other eight methods select $0.5$ or less at every $\kappa$, and the three carrying an error-feedback estimator select $0$ throughout.

Every number in this subsection comes from a single scripted run of the released benchmark, which records each stage together with the commit, GPU and wall time it ran under.

\subsection{Centralized training results}\label{app:images_cen}
\begin{figure*}[!t]
    \centering
    \includegraphics[width=\textwidth]{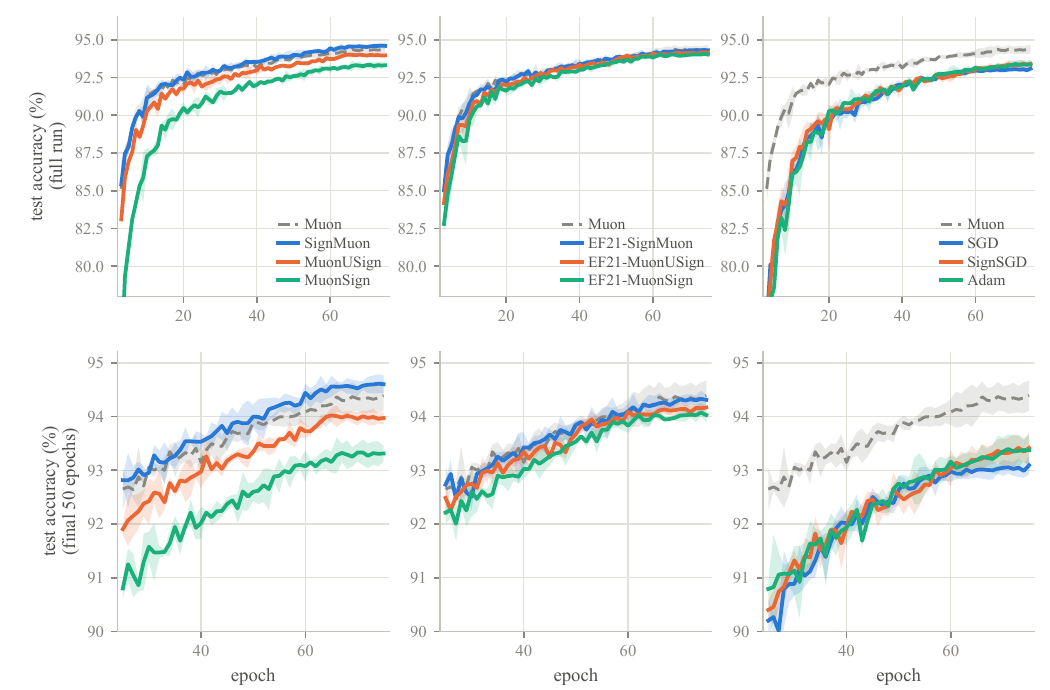}
    \caption{Test accuracy over the whole run (top row) and from epoch~$25$ (bottom row), the same data at two scales, grouped as in \hyperref[fig:cifar_results]{Figure~\ref*{fig:cifar_results}}. Bands are $\pm1$ standard deviation over three seeds; Muon is the gray dashed reference in every panel. Series are named in each panel's legend.}
    \label{fig:cifar_curves_appendix}
\end{figure*}

\begin{figure*}[!t]
    \centering
    \includegraphics[width=\textwidth]{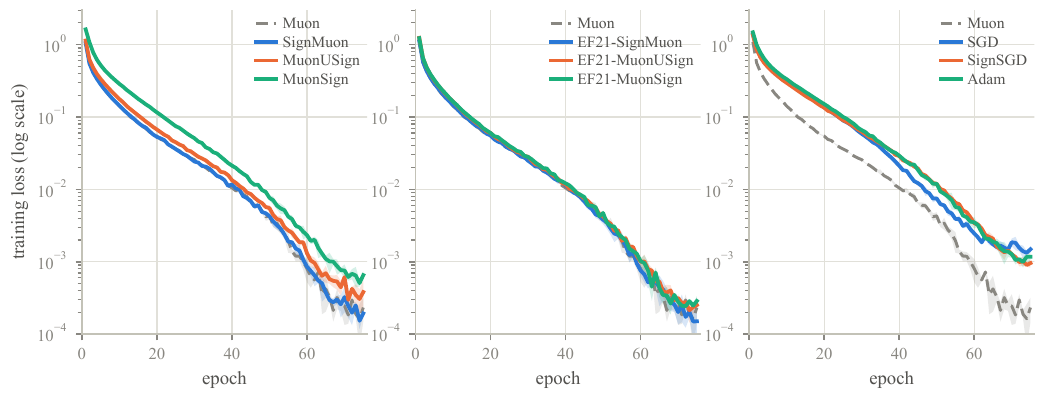}
    \caption{Training loss (log scale), grouped as in \hyperref[fig:cifar_results]{Figure~\ref*{fig:cifar_results}} and over the same runs. Within a panel the methods are hard to separate (eight of the ten fall below $10^{-3}$, differing mainly in how early), which is why the body figure carries the accuracy alone. What the loss does show is the separation between families: in the right panel Muon reaches $10^{-2}$ six to eleven epochs before SGD, SignSGD and Adam and terminates four to six times lower, while in the left panel SignMuon tracks Muon and the sign-before and both-sides orderings lie above it, in the same order as their accuracies. The three EF21 variants (centre) are indistinguishable in loss, whereas their accuracies span $0.27$ points.}
    \label{fig:cifar_train_loss}
\end{figure*}

\begin{figure*}[!t]
    \centering
    \includegraphics[width=\textwidth]{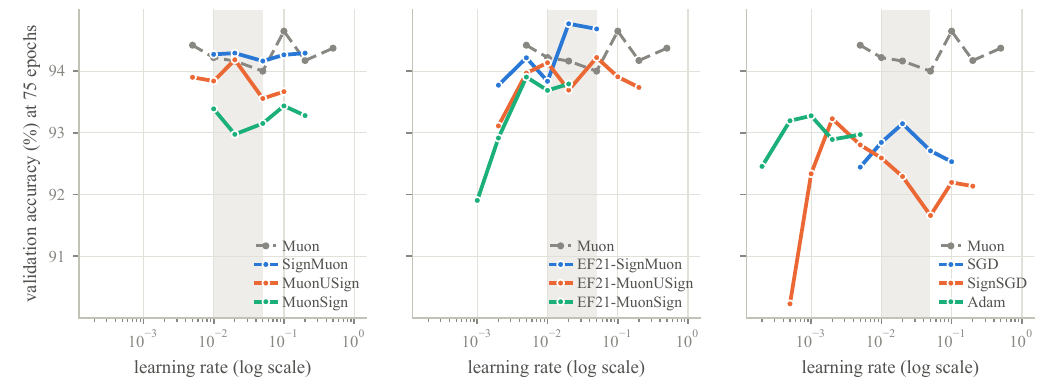}
    \caption{Learning-rate sweep at the $75$-epoch selection horizon, on the tuning split and on the validation accuracy that selection ranked, grouped as in \hyperref[fig:cifar_results]{Figure~\ref*{fig:cifar_results}}, with Muon repeated as the gray dashed reference. The shaded band marks $\eta_0 \in [0.01, 0.05]$, the only interval every method but Adam was swept over; the nine vary by $0.10$ to $0.93$ points within it, SGD ($0.44$) among them, so an \emph{absolute} window no longer discriminates between the families, their optima lying two orders of magnitude apart. The matched comparison is the relative one. Within a factor of five of its own optimum SignMuon varies by $0.13$ points and SignSGD by $3.00$; the two emit steps of identical Frobenius length, so the difference is attributable to the \textsc{lmo}. The remaining eight lie between these, from $0.46$ (MuonSign) to $2.00$ (EF21-MuonSign). Taking the sign of the oracle's output does not reduce this insensitivity: SignMuon's sweep is the flattest of the ten.}
    \label{fig:cifar_lr}
\end{figure*}

\hyperref[tab:cifar_central]{Table~\ref*{tab:cifar_central}} reports the centralized CIFAR-10 results on ResNet-18 in full: the selected $\eta_0$, the final train and test accuracies, the epoch count to the $90\%$ threshold, and the per-epoch cost. \hyperref[fig:cifar_results]{Figure~\ref*{fig:cifar_results}} starts at epoch~$25$, where the methods are a few points apart rather than twenty. Every run uses batch $128$, momentum $0.9$, auxiliary rate $10^{-3}$, and zero weight decay, so that the matrix rule and $\eta_0$ are the only quantities that vary. \hyperref[fig:cifar_curves_appendix]{Figure~\ref*{fig:cifar_curves_appendix}} gives the accuracy curves from epoch~$3$, \hyperref[fig:cifar_train_loss]{Figure~\ref*{fig:cifar_train_loss}} the training loss, and \hyperref[fig:cifar_lr]{Figure~\ref*{fig:cifar_lr}} the learning-rate sweep behind the selection.

\begin{table*}[!t]
\centering
\small
\renewcommand{\arraystretch}{1.2}
\begin{tabular}{|c|l|c|c|c|c|c|c|}
\hline
\textbf{Dataset} & \textbf{Optimizer} & \textbf{Epochs} & $\boldsymbol{\eta_0}$ & \textbf{Train Acc} & \textbf{Test Acc} & \textbf{Ep.\ to $\mathbf{90\%}$} & \textbf{s/epoch} \\
\hline
\multirow{10}{*}{CIFAR-10}
& SignMuon & \multirow{10}{*}{75} & 0.02 & $99.99$ & $94.60 \pm 0.15$ & 7.7 & 16.6 \\
\cline{2-2} \cline{4-8}
& Muon & & 0.1 & $99.99$ & $94.35 \pm 0.27$ & 7.7 & 14.4 \\
\cline{2-2} \cline{4-8}
& EF21-SignMuon & & 0.02 & $100.00$ & $94.31 \pm 0.11$ & 9.0 & 16.6 \\
\cline{2-2} \cline{4-8}
& EF21-MuonUSign & & 0.05 & $99.99$ & $94.14 \pm 0.07$ & 10.3 & 16.7 \\
\cline{2-2} \cline{4-8}
& EF21-MuonSign & & 0.005 & $99.99$ & $94.04 \pm 0.10$ & 11.0 & 18.1 \\
\cline{2-2} \cline{4-8}
& MuonUSign & & 0.02 & $99.99$ & $93.98 \pm 0.12$ & 10.3 & 16.2 \\
\cline{2-2} \cline{4-8}
& Adam & & 0.001 & $99.97$ & $93.37 \pm 0.27$ & 20.0 & 12.7 \\
\cline{2-2} \cline{4-8}
& SignSGD & & 0.002 & $99.97$ & $93.37 \pm 0.26$ & 19.0 & 12.5 \\
\cline{2-2} \cline{4-8}
& MuonSign & & 0.1 & $99.98$ & $93.31 \pm 0.21$ & 17.7 & 17.2 \\
\cline{2-2} \cline{4-8}
& SGD & & 0.02 & $99.96$ & $93.04 \pm 0.14$ & 21.3 & 12.2 \\
\hline
\end{tabular}
\caption{CIFAR-10: centralized learning on ResNet-18, three seeds. Test accuracy is the mean $\pm$ standard deviation of the last five epochs; train accuracy, ``Ep.\ to $90\%$'' and ``s/epoch'' are means over the same three seeds. $\eta_0$ was selected on a held-out $5$k validation split at the same $75$-epoch horizon the table reports, and the selected rate then retrained on the full $50$k set, as described in \hyperref[app:repro]{Appendix~\ref*{app:repro}}. Every method fits the training set to within $0.05$ points of the others, so the column separates nothing on its own; the differences in the test column are not differences in how far training got.}
\label{tab:cifar_central}
\end{table*}

\paragraph{Selection horizon.} Selection and reporting share the $75$-epoch horizon, so each $\eta_0$ above is the argmax of a validation sweep conducted at the length the table reports. An earlier protocol that selected at $15$ epochs and reported at $75$ was abandoned because the two horizons selected different rates for two of the methods. The selection is not sharply peaked for the sign-after methods: over the rates within a factor of five of its own optimum, SignMuon's validation accuracy varies by $0.13$ points, against $0.65$ for Muon and $3.00$ for SignSGD (\hyperref[fig:cifar_lr]{Figure~\ref*{fig:cifar_lr}}).

\paragraph{Weight decay.} The primary table is unregularized, the setting the theorems analyse and the one \citet{mishra2026signmuon}'s sweep selects. Repeating the top three at the same $\eta_0$ with decoupled decay $5\times10^{-4}$ (seed $0$, decay not re-tuned) displaces each by at most $0.34$ points, of the order of the three-seed standard deviation of the undecayed runs, and leaves the three within $0.23$ of one another: Muon $94.46\%$ against $94.12\%$ undecayed at that seed, EF21-SignMuon $94.45\%$ against $94.42\%$, SignMuon $94.23\%$ against $94.43\%$. The ordering within that interval does change, but at a single seed and over so narrow a range it is not a measured effect; what the ablation establishes is that decay introduces no separation where the primary table shows none. Decay is applied decoupled, $\mathbf{X} \mathrel{*}= 1 - \eta\lambda$, so the \textsc{lmo} sees the true gradient; the coupled convention would only rotate the direction, since every step here is scale-invariant.

% \begin{figure*}[!t]
%     \centering
%     \begin{subfigure}[b]{0.48\textwidth}
%         \centering
%         \includegraphics[width=\textwidth]{images/centralized_images/train_loss_resnet.png} 
%         \caption{Train Loss}
%     \end{subfigure}
%     \hfill
%     \begin{subfigure}[b]{0.48\textwidth}
%         \centering
%         \includegraphics[width=\textwidth]{images/centralized_images/train_acc_resnet.png} 
%         \caption{Train Accuracy}
%     \end{subfigure}
%     \caption{Centralized setting. Comparison of convergence and classification accuracy achieved by different optimization algorithms during CIFAR-10 classifier training on the training set on ResNet-18.}
%     \label{fig:cifar_results_train}
% \end{figure*}

% \newpage
\subsection{Communication accounting}\label{app:commacct}
\hyperref[tab:exp_3]{Table~\ref*{tab:exp_3}} quotes mean bits per parameter per round. This section states precisely what is counted, since the headline ``$32\times$'' of the sign-compression literature is an idealization eroded by three separate effects, only one of which is customarily acknowledged.

Write $P_{\mathrm{mat}}$ for the number of matrix parameters, $P_{\mathrm{aux}}$ for the auxiliary group (biases, BatchNorm affine parameters, the classifier head), $L$ for the number of matrix layers, and $P = P_{\mathrm{mat}}+P_{\mathrm{aux}}$. On CNN2, $P_{\mathrm{mat}} = 762{,}560$, $P_{\mathrm{aux}} = 2{,}146$ and $L=3$. All figures below are \emph{per client per round}: the uplink is what one client sends, the downlink what the server sends to one client.

\paragraph{(i) The uplink alphabet.} The randomized sign of \hyperref[sec:theory]{Section~\ref*{sec:theory}} renders every transmitted symbol a genuine bit, so the uplink costs one bit per parameter with no entropy coding required to realize it.

\paragraph{(ii) The auxiliary group is never compressed.} It travels at full precision in both directions for every method, so a ``one-bit'' channel actually costs
\[
\frac{1\cdot P_{\mathrm{mat}} + 32\,P_{\mathrm{aux}}}{P} = 1.087\ \text{bits},
\]
a $29.4\times$ reduction rather than $32\times$. On CNN2 the group is $0.28\%$ of the parameters; on an architecture with a large embedding or head it would dominate this table, which is why the quantity is computed per model rather than quoted.

\paragraph{(iii) Error feedback carries one scale per layer.} Both EF21 channels transmit the pair $(\mathbf{s}_t,\alpha_t)$ with one full-precision $\alpha_t$ per matrix layer: on the uplink for every EF21 method, and again on the downlink for EF21-MuonSign. That is $32L$ bits, adding $32L/P \approx 1.3\times10^{-4}$ bits per parameter here. That lies four decimal places in, and it is reported for completeness rather than because it alters a conclusion: it is the difference between ``one bit per parameter'' and ``one bit per parameter, plus a constant''.

\paragraph{(iv) Which methods compress the downlink.} The criterion is not whether the method applies a compressor but whether the object the server must distribute is already $\pm1$-valued. Three cases qualify: the majority vote $\mathbf{s}_t^{\mathrm{agg}}$ itself (SignMuon, SignSGD), a signed \textsc{lmo} output (MuonSign), and a primal error-feedback residual (EF21-MuonSign). In the first of these the server broadcasts the vote rather than the model and each client applies the step to its local copy; the copies start from a common $\mathbf{X}_0$ and receive identical updates, so they never drift. The remaining three methods must distribute a dense server-side quantity, namely $\operatorname{polar}(\cdot)$ of the aggregate for MuonUSign and EF21-MuonUSign and a scaled average of signs for EF21-SignMuon, and therefore transmit it at full precision.

\begin{table}[!t]
\centering\small
\setlength{\tabcolsep}{4pt}
\begin{tabular}{@{}lcccc@{}}
\toprule
\textbf{Method} & \textbf{Up (bits)} & \textbf{Down (bits)} & \textbf{Up} & \textbf{Down} \\
\midrule
Muon, MuonServer     & $32$     & $32$     & $1.0\times$  & $1.0\times$ \\
SGD, Adam            & $32$     & $32$     & $1.0\times$  & $1.0\times$ \\
\midrule
SignMuon             & $1.0870$ & $1.0870$ & $29.4\times$ & $29.4\times$ \\
SignSGD              & $1.0870$ & $1.0870$ & $29.4\times$ & $29.4\times$ \\
MuonSign             & $1.0870$ & $1.0870$ & $29.4\times$ & $29.4\times$ \\
EF21-MuonSign        & $1.0871$ & $1.0871$ & $29.4\times$ & $29.4\times$ \\
\midrule
MuonUSign            & $1.0870$ & $32$     & $29.4\times$ & $1.0\times$ \\
EF21-SignMuon        & $1.0871$ & $32$     & $29.4\times$ & $1.0\times$ \\
EF21-MuonUSign       & $1.0871$ & $32$     & $29.4\times$ & $1.0\times$ \\
\bottomrule
\end{tabular}
\caption{Communication per client per round on CNN2, under the randomized-zero convention. Bits are per model parameter; reductions are against a $32$-bit baseline. The fourth decimal separates the error-feedback methods from the rest: it is the per-layer scale of point~(iii). Groups: uncompressed references, one bit in both directions, one bit on the uplink only.}
\label{tab:commacct}
\end{table}

\hyperref[tab:commacct]{Table~\ref*{tab:commacct}} collects the four effects into the per-round cost of each method, and is the source of the Up and Down columns of \hyperref[tab:exp_3]{Table~\ref*{tab:exp_3}}. It is computed by \texttt{federated.algorithms.communication\_bits} from the alphabet and the measured zero rate of the run it describes, so a run made under the legacy ternary convention reports its own higher figure rather than the idealized one. Every run behind \hyperref[tab:exp_3]{Table~\ref*{tab:exp_3}} was made under the randomized convention, which makes these figures realized rather than idealized. Two diagnostics record what the convention had to absorb, both counted \emph{before} the randomized mapping and therefore feeding no accounting. Exact zeros do occur on the uplink, at up to $3.7\%$ of coordinates for MuonSign and $0.5\%$ for SignMuon, and at none at all for the three error-feedback methods, whose compressed quantity is a residual rather than a direction. And the majority vote tied in no coordinate of any evaluated round of any reported run, as at $N=11$ and $\pm1$ client messages it cannot.

\subsection{Federated training results}\label{app:results_fed:3}
\begin{figure*}[!t]
    \centering
    \includegraphics[width=\textwidth]{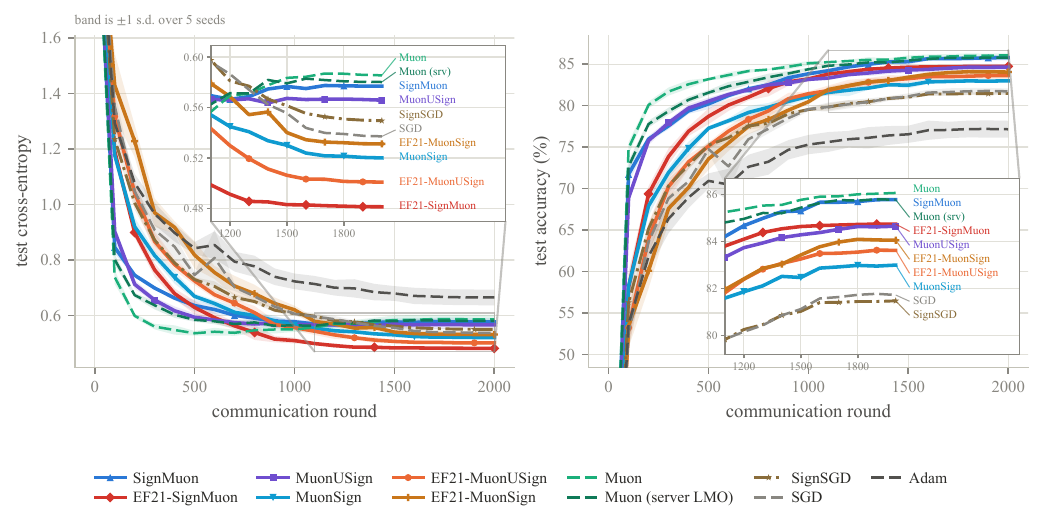}
    \caption{CIFAR-10 federated learning on CNN2 at $N=11$ clients and batch $192$ per client, every method evaluated on the exact server model $\mathbf{X}_t$. Both panels are clipped to exclude round $0$, the untrained model, which is identical for every method. Solid lines are the one-bit methods, dashed the uncompressed references and dash-dot SignSGD, so that hue is not the only channel separating eleven curves. The insets magnify the final $45\%$ of rounds, where every number in \hyperref[tab:exp_3]{Table~\ref*{tab:exp_3}} is read, and name each curve at its own end; Adam finishes clear of the other ten and so falls outside the magnified window. Eleven methods are drawn, one more than \hyperref[tab:exp_3]{Table~\ref*{tab:exp_3}}, the extra being the server-side-LMO Muon control.}
    \label{fig:exp_3}
\end{figure*}

\hyperref[fig:exp_3]{Figure~\ref*{fig:exp_3}} (together with \hyperref[tab:exp_3]{Table~\ref*{tab:exp_3}} in the main text) reports the comparison of optimizers at $N=11$ clients. The figure additionally shows Muon with a server-side LMO, which is not a communication-efficient method but isolates whether moving the oracle off the clients carries a penalty on its own; it does not ($85.74 \pm 0.11\%$ against Muon's $85.98 \pm 0.26\%$), so the gaps in \hyperref[tab:exp_3]{Table~\ref*{tab:exp_3}} are attributable to compression and sign placement rather than to where the oracle runs.

On the sign-after placement the cost of error feedback depends on the setting, and federation is where it is largest: EF21-SignMuon lies $1.00$ points below SignMuon here, against $0.29$ points centrally (\hyperref[tab:cifar_main]{Table~\ref*{tab:cifar_main}}), while on nanoGPT the two are indistinguishable (\hyperref[tab:nanogpt]{Table~\ref*{tab:nanogpt}}).

Two features of the figure do not appear in the table. The threshold column does not order the methods as final accuracy does: MuonUSign crosses $80\%$ in $500$ rounds, ahead of SignMuon's $540$ and EF21-SignMuon's $640$, and finishes below both, so which method leads depends on where the round budget is cut. And the left panel orders them differently again: the lowest test cross-entropy of the eleven is EF21-SignMuon's, $0.482 \pm 0.005$, and the highest bar Adam's is Muon's, $0.586 \pm 0.016$, the reverse of how those two stand on accuracy. Cross-entropy and accuracy are not obliged to agree, and it is accuracy the comparison is about; the loss panel is drawn so that the disagreement is on the record rather than suppressed by the choice of metric.

\subsection{Language-modelling details}\label{app:nanogpt}

\begin{figure*}[!t]
\centering
\includegraphics[width=\textwidth]{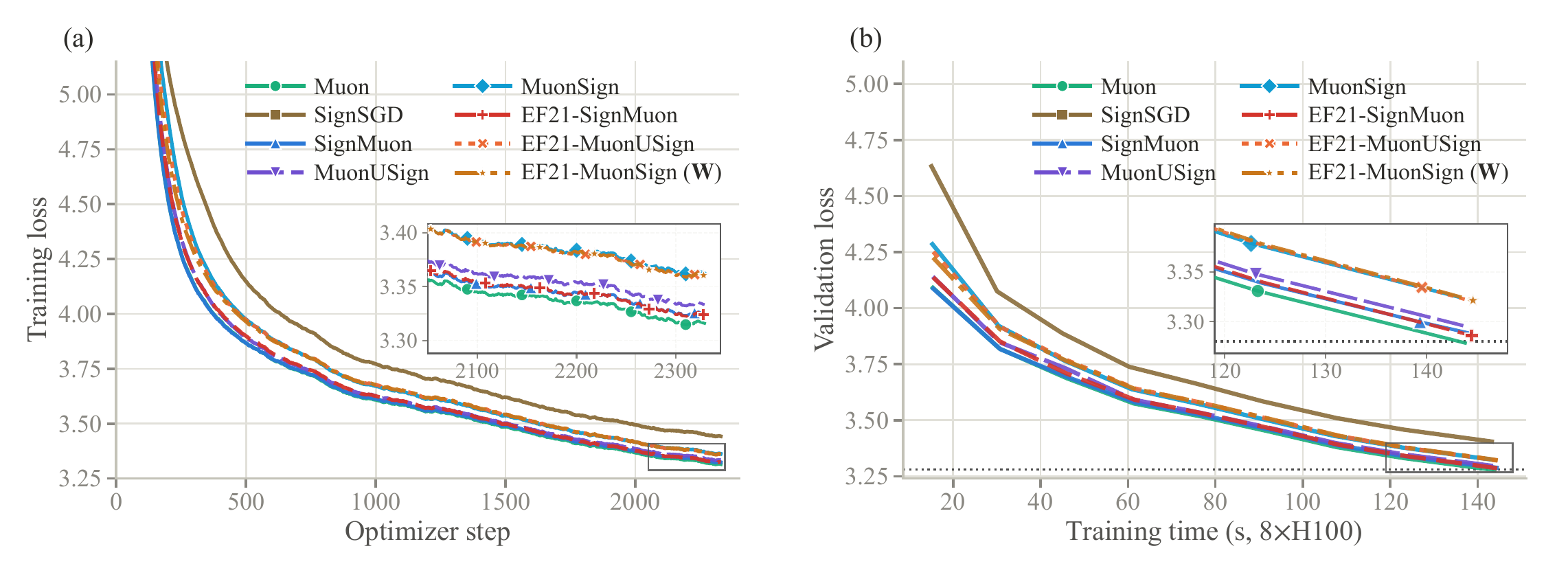}
\caption{NanoGPT, supporting curves. (a)~Training loss against optimizer step (EMA-smoothed; the logged quantity is a per-rank sum over $32{,}768$ tokens, divided out here). (b)~Validation loss against the speedrun clock, which excludes validation and compilation. Insets magnify the boxed tails, as in \hyperref[fig:nanogpt]{Figure~\ref*{fig:nanogpt}}. Marker positions are staggered in both figures so that curves lying within a line width of one another can still be followed individually. The ordering is the same on both axes: no method incurs a measurable wall-clock premium for its compressor or its error-feedback buffers. Note in (a) that EF21-MuonSign's training loss, taken at $\mathbf{W}$, exhibits no degradation throughout, which is what localizes its validation gap to the tracking of $\mathbf{X}$ rather than to training.}
\label{fig:nanogpt_appendix}
\end{figure*}

\paragraph{Setup.} Upstream modded-nanoGPT record \#40 (2025-10-04), the last record before NorMuon and hence the last whose hidden-matrix optimizer is a clean, separable momentum\,$\to$\,LMO\,$\to$\,step Muon, so our variants inject exactly at the LMO. Model: $12$ layers, model dimension $768$, $6$ heads of dimension $128$, vocabulary $50{,}257$; hidden matrices are $768\times3072$ (the merged $\mathbf{Q}/\mathbf{K}/\mathbf{V}/\mathbf{O}$ weight $\mathtt{qkvo\_w}$ is used as four $768\times768$ blocks, and both the LMO and the compressor scale are applied per block; the two MLP matrices are the up-projection $\mathtt{c\_fc}$ and the zero-initialized output projection $\mathtt{c\_proj}$, which maps the $3072$-dimensional hidden activation back to the model dimension). Data: FineWeb10B, the $10$B-token sample of FineWeb \citep{penedo2024fineweb} that the speedrun repository prepares and tokenizes; $262{,}144$ tokens per step, $2330$ steps ($=611$M tokens), validation on the fixed $10{,}485{,}760$-token split. Hardware: one rented $8\times$H100 SXM node ($80$\,GB per GPU; dual Xeon Platinum host, $224$ vCPU, $2$\,TB RAM, PCIe~5.0 $\times16$, NVMe scratch), driver $595.71.05$, running PyTorch $2.10.0{+}\mathrm{cu}128$ under Python $3.12.3$ in a virtual environment of its own rather than the container's torch, since the prebuilt Flash-Attention-3 kernel the record fetches exists for no CUDA-13 build; one process per GPU. Gradients are averaged by \texttt{reduce\_scatter} so the owning rank runs the centralized update and \texttt{all\_gather} returns the parameter, i.e.\ the compression is a property of the update rule, as in the centralized algorithms we analyze.

\paragraph{Hyperparameters.} Matrix/gate optimizer: $\eta_0 = 0.06$ (\textsc{lmo} family) or $0.03$ (\textsc{sign} family), per-layer scaled by the unit-gain rule (\hyperref[app:lrscale]{Appendix~\ref*{app:lrscale}}); Nesterov momentum $\mu = 0.95$, warmed up linearly from $0.85$ over the first $300$ steps and cooled back to $0.85$ over the last $50$; weight decay $0$ (the record's own value); $\eta$ constant then linearly cooled to $0.1\eta_0$ over the final $45\%$ of the $2290$ scheduled iterations, with the $40$-step extension held at that floor; LMO by $5$ Polar-Express iterations. Auxiliary parameters (embeddings, scalars, head) use the record's distributed Adam unchanged: $\eta = 0.008$, $\beta = (0.65, 0.95)$, $\varepsilon = 10^{-8}$, no weight decay, per-parameter multipliers $75$ on embeddings and $5$ on scalars, stepped every other iteration. Nothing above was tuned by us: the \textsc{lmo} family runs at the record's own $\eta_0$, and every value outside the matrix optimizer is the record's. Wall-clock varies by at most $1.1\%$ across all eight methods ($61.4$--$62.1$\,ms/step), all of them some $2\%$ above the record's own $60.4$\,ms/step, which every method pays equally: our port replaces its Triton kernels and batched sharded transport with a pure-torch per-parameter equivalent.

\begin{figure}[!t]
\centering
\includegraphics[width=0.65\textwidth]{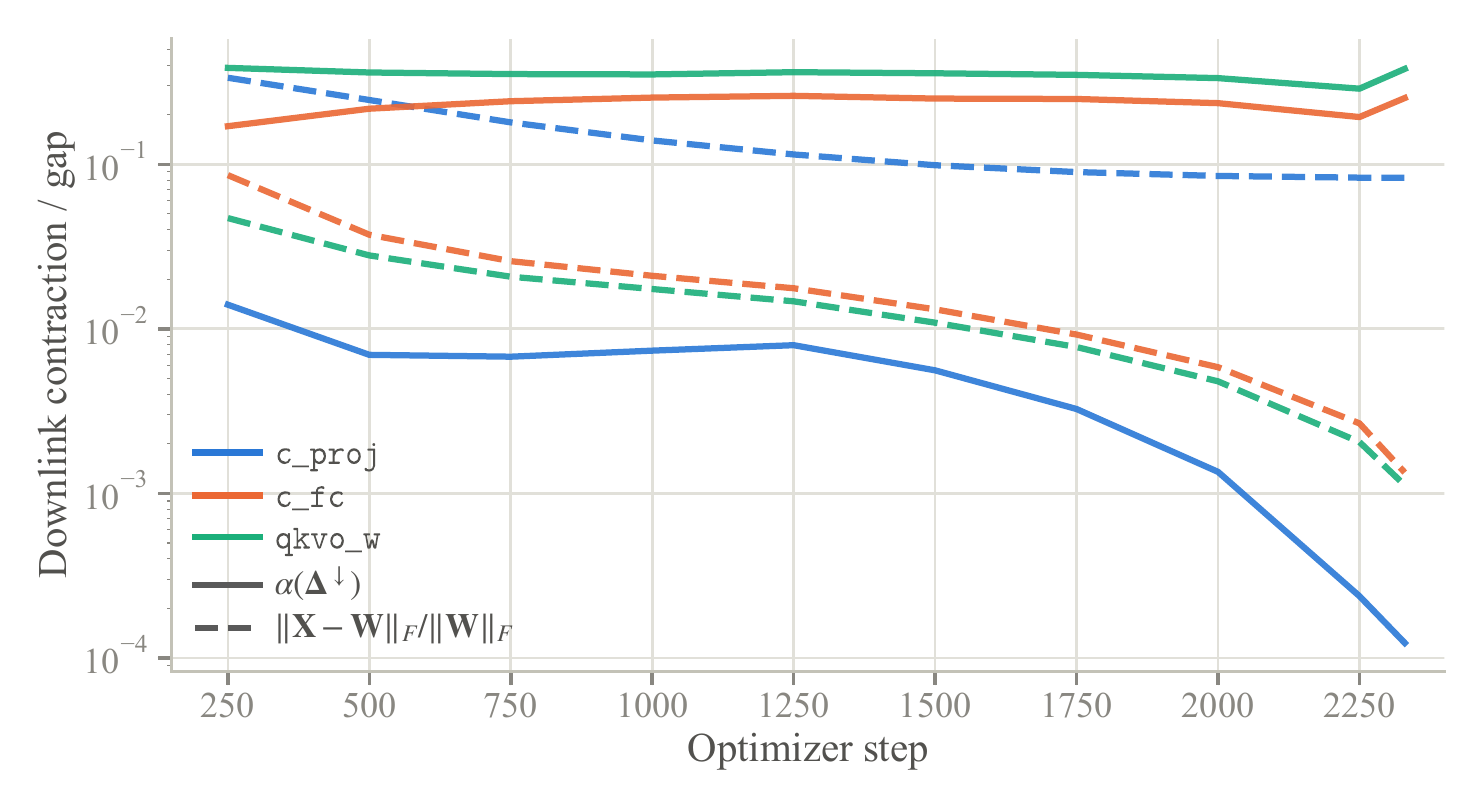}
\caption{The downlink measurement for EF21-MuonSign, per layer type: the contraction $\alpha(\boldsymbol{\Delta}^{\downarrow})$ the scaled sign actually achieves (solid) and the resulting server/broadcast gap $\|\mathbf{X}_t-\mathbf{W}_t\|_F/\|\mathbf{W}_t\|_F$ (dashed). On the zero-initialized $\mathtt{c\_proj}$ the contraction sits one to two orders of magnitude below the other two from step~$250$ on, then collapses over the last $500$ steps to $1.2\times10^{-4}$, three orders below them, while the gap stops decreasing near $10^{-1}$; the other two contract at $\Theta(10^{-1})$ throughout and close their gaps by $42$ and $64\times$ over the run, ending near $10^{-3}$. This is \hyperref[rem:dimension]{Remark~\ref*{rem:dimension}} rendered as a measurement.}
\label{fig:nanogpt_diag}
\end{figure}

\paragraph{Supporting curves.} \hyperref[fig:nanogpt_appendix]{Figure~\ref*{fig:nanogpt_appendix}} gives the two views \hyperref[fig:nanogpt]{Figure~\ref*{fig:nanogpt}} omits: training loss against optimizer step, and validation loss against the speedrun clock. The second is what supports the claim of equal wall-clock in \hyperref[sec:exp_nanogpt]{Section~\ref*{sec:exp_nanogpt}}; the first shows that EF21-MuonSign trains normally at $\mathbf{W}$, which is what places its validation gap in the tracking of $\mathbf{X}$ rather than in training. \hyperref[fig:nanogpt_diag]{Figure~\ref*{fig:nanogpt_diag}} then measures that tracking directly, per layer type.

\paragraph{Compressor diagnostics.} \hyperref[tab:nanogpt_diag]{Table~\ref*{tab:nanogpt_diag}} reports, per layer type at the final step, the contraction each scaled sign achieves, $\alpha(\boldsymbol{\Delta}) = \|\boldsymbol{\Delta}\|_1^2/(d\|\boldsymbol{\Delta}\|_F^2)$, and the relative estimator lag. Three observations merit separate comment. \emph{(i)}~Every uplink is well-contractive, $\alpha \in [0.32, 0.64]$ against the isotropic $2/\pi = 0.637$, and EF21-SignMuon's is uniformly the best, the entries of its orthogonal target being the most evenly spread; yet the uplink \emph{lag} is large ($0.61$--$0.82$), so that the estimator remains far from its target at every step and the methods nonetheless train well, the LMO being scale-invariant and requiring only the direction. \emph{(ii)}~The single anomaly in the table is the $\mathtt{c\_proj}$ downlink, $\alpha = 1.2\times10^{-4}$, a factor $4.8\times10^{3}$ below its uplink on the very same layer. Since both compressors are the same operator, the difference is a property of the residual, not of the compressor: the uplink residual is refreshed by an exogenous stochastic gradient each round, the downlink residual is generated by the compressor's own recursion (\hyperref[rem:dimension]{Remark~\ref*{rem:dimension}}). \emph{(iii)}~The resulting validation loss at the server model $\mathbf{X}$ falls to $4.20$ by step~$750$, rises to $5.52$ by step~$1500$, and then holds there ($5.51$--$5.58$ to the end of the run): a persistent offset above the broadcast model, not a divergence.

\paragraph{Where the two models part.} The $\mathtt{c\_proj}$ anomaly of item~\emph{(ii)} is what separates EF21-MuonSign's two models: on that layer the gap $\|\mathbf{X}_t-\mathbf{W}_t\|_F/\|\mathbf{W}_t\|_F$ stops decreasing at $8\%$ against ${<}1\%$ for every other layer and both gates, and it is the only layer whose gap fails to close, which is what localizes the ${\approx}2.2$-nat offset of item~\emph{(iii)} to it. The reason is the initialization. A layer built from zero receives maximally correlated updates, so its downlink residual concentrates, and a compressor that moves every coordinate by $\operatorname{mean}|\boldsymbol{\Delta}^{\downarrow}|$ cannot catch the coordinates driven hardest. This is the mechanism of \hyperref[rem:dimension]{Remark~\ref*{rem:dimension}} rather than a tuning failure. Lowering $\eta_0$ does not repair it, because the admissible step size would have to shrink by a further $\sqrt{r}$ in the layer rank (\hyperref[rem:normmismatch]{Remark~\ref*{rem:normmismatch}}); only a downlink compressor contractive in the \emph{spectral} norm would.

\begin{table}[!t]
\centering\small
\setlength{\tabcolsep}{3.5pt}
\begin{tabular}{@{}lcccc@{}}
\toprule
& \multicolumn{2}{c}{\textbf{uplink}} & \multicolumn{2}{c}{\textbf{downlink}} \\
\cmidrule(lr){2-3}\cmidrule(lr){4-5}
\textbf{Layer type} & $\alpha$ & lag & $\alpha$ & gap \\
\midrule
\multicolumn{5}{@{}l}{\emph{EF21-SignMuon}} \\
\quad $\mathtt{qkvo\_w}$ ($\times10$)   & 0.640 & 0.66 & -- & -- \\
\quad $\mathtt{c\_fc}$ ($\times11$)     & 0.633 & 0.62 & -- & -- \\
\quad $\mathtt{c\_proj}$ ($\times11$)   & 0.634 & 0.62 & -- & -- \\
\midrule
\multicolumn{5}{@{}l}{\emph{EF21-MuonUSign}} \\
\quad $\mathtt{qkvo\_w}$   & 0.374 & 0.81 & -- & -- \\
\quad $\mathtt{c\_fc}$     & 0.323 & 0.82 & -- & -- \\
\quad $\mathtt{c\_proj}$   & 0.597 & 0.62 & -- & -- \\
\midrule
\multicolumn{5}{@{}l}{\emph{EF21-MuonSign}} \\
\quad $\mathtt{qkvo\_w}$   & 0.406 & 0.79 & 0.380 & 0.0011 \\
\quad $\mathtt{c\_fc}$     & 0.355 & 0.78 & 0.253 & 0.0013 \\
\quad $\mathtt{c\_proj}$   & 0.596 & 0.61 & $\mathbf{1.2{\times}10^{-4}}$ & $\mathbf{0.083}$ \\
\bottomrule
\end{tabular}
\caption{Compressor diagnostics at step $2330$, medians over the identical layers of each type. $\alpha$ is the contraction the scaled sign achieves on that round's residual ($2/\pi$ for an isotropic residual, $1/d$ in the worst case); ``lag'' is $\|\text{target}-\text{estimator}\|_F/\|\text{target}\|_F$; ``gap'' is $\|\mathbf{X}_t-\mathbf{W}_t\|_F/\|\mathbf{W}_t\|_F$. Downlink columns are empty for the two methods that broadcast exactly. Gates are omitted (they are $6\times12$ and $1\times12$).}
\label{tab:nanogpt_diag}
\end{table}

\subsection{Per-Layer Step Sizes: the Unit-Gain Rule}\label{app:lrscale}

This appendix derives the heuristic \eqref{eq:unit_gain_main} stated in the main text, selects its one free exponent by measurement, and delimits its scope. SGD and Adam have no norm-fixed step, so the rule does not apply to them; both are run at one global rate throughout.

Our counterexamples, like the analysis of \citet{mishra2026signmuon}, concern a single matrix, where one scalar step size suffices. A network has layers of very different shapes, and the methods of this paper produce step matrices from \emph{two} families whose norms scale differently with shape, so a single global $\eta$ cannot be simultaneously correct for both families and across layers. This is not a lacuna we are obliged to tolerate: the layer-wise LMO framework to which our convergence result reduces \citep{riabinin2025gluon,gruntkowska2025error} already carries per-layer norms $\|\cdot\|_{(\ell)}$, smoothness constants $L_\ell$ and radii $\eta_\ell$; it is only the \emph{experiments} that have hitherto fixed that radius at a constant. What follows instantiates it. The specific rule is this paper's own; the criterion behind it is borrowed (it is the average-case form of the spectral scaling condition of the maximal-update literature \citep{yang2021tuning,yang2023spectral,large2024modular}), and its \textsc{lmo} branch reproduces the aspect factor Muon already uses in practice \citep{jordan2024muon}, which is the external check we rely on.

\paragraph{The two families.} For a parameter reshaped to $\mathbf{X}\in\mathbb{R}^{m\times n}$ ($m$ the output dimension, $n$ the input dimension, $r=\operatorname{rank}=\min(m,n)$ generically), the step matrix $\mathbf{P}$ that a method applies belongs to one of two families:
\begin{equation}
\label{eq:two_families}
\mathbf{P}=\mathbf{U}\mathbf{V}^\top\ \ (\textsc{lmo}),
\qquad
\mathbf{P}\in\{\pm1\}^{m\times n}\ \ (\textsc{sign}).
\end{equation}
The \textsc{lmo} family comprises Muon, MuonUSign, EF21-MuonUSign, EF21-MuonSign and EF21-SignMuon; the \textsc{sign} family comprises SignMuon, MuonSign and SignSGD.
Both have \emph{exactly} known Frobenius norms: $\|\mathbf{U}\mathbf{V}^\top\|_F^2=\operatorname{tr}(\mathbf{V}\mathbf{U}^\top\mathbf{U}\mathbf{V}^\top)=\operatorname{tr}(\mathbf{V}^\top\mathbf{V})=r$, so $\|\mathbf{P}\|_F=\sqrt{r}$; and a $\pm1$ matrix has $\|\mathbf{P}\|_F=\sqrt{mn}$. (EF21-SignMuon steps along the error-feedback estimator $\mathbf{d}_t^{\mathrm{est}}$ of $\operatorname{polar}(\tilde{\mathbf{M}}_t)$ rather than along the oracle output itself, so for it $\|\mathbf{P}\|_F=\sqrt{r}$ holds only in the limit; we assign it to the \textsc{lmo} family on that basis.)

\paragraph{The criterion.} Define the \emph{RMS gain} of $\mathbf{A}\in\mathbb{R}^{m\times n}$, i.e.\ how much it amplifies a generic input in root-mean-square terms, with $\operatorname{rms}(\mathbf{v}):=\|\mathbf{v}\|/\sqrt{\dim}$:
\begin{equation}
\label{eq:gain_def}
\gamma(\mathbf{A})^2:=\frac{\mathbb{E}_{\mathbf{u}}\!\left[\operatorname{rms}(\mathbf{A}\mathbf{u})^2\right]}{\mathbb{E}_{\mathbf{u}}\!\left[\operatorname{rms}(\mathbf{u})^2\right]},\qquad \mathbf{u}\sim\mathcal{N}(\mathbf{0},\mathbf{I}_n).
\end{equation}
Since $\mathbb{E}\|\mathbf{A}\mathbf{u}\|^2=\operatorname{tr}(\mathbf{A}^\top\mathbf{A})=\|\mathbf{A}\|_F^2$ and $\mathbb{E}\|\mathbf{u}\|^2=n$, \eqref{eq:gain_def} evaluates in closed form:
\begin{equation}
\label{eq:gain_identity}
\gamma(\mathbf{A})=\frac{\|\mathbf{A}\|_F}{\sqrt{m}} .
\end{equation}
The single modelling assumption is that $\mathbf{u}$ is isotropic \emph{and independent of} $\mathbf{A}$, to which we return below. Controlling a layer update's RMS-to-RMS effect is precisely the desideratum of the spectral scaling condition \citep{yang2023spectral} and of the modular norm \citep{large2024modular}; \eqref{eq:gain_def} is its average-case (isotropic-input) version, and when the independence fails, that is, for aligned inputs, where the operator norm governs the gain instead, one recovers the $\mu$P value $a=1$ discussed below.

Every standard initialization with variance $\propto1/n$ has \emph{shape-independent} gain. For He normal ($\sigma^2=2/n$), $\|\mathbf{X}\|_F=\sigma\sqrt{mn}=\sqrt{2m}$, so $\gamma(\mathbf{X})=\sqrt{2}$; for PyTorch's default Kaiming-uniform convolution ($\sigma^2=1/(3n)$), $\gamma(\mathbf{X})=1/\sqrt{3}$. Either way a constant. Requiring the update's gain to be a fixed fraction of the weight's is therefore simply the requirement that \emph{the per-step gain be the same on every layer}, and by \eqref{eq:gain_identity} that is one formula:
\begin{equation}
\label{eq:unit_gain}
\boxed{\ \eta_\ell=\eta_0\,\lambda_\ell,\qquad \lambda_\ell=\frac{\sqrt{m}}{\|\mathbf{P}\|_F}\ }
\end{equation}
which gives $\gamma(\eta_\ell\mathbf{P})=\eta_0$ exactly, for every shape and both families, so that $\eta_0$ \emph{is} the per-step RMS gain. One caveat: the rule is derived from $\|\mathbf{U}\mathbf{V}^\top\|_F=\sqrt{r}$, which holds for the exact oracle. Five Newton--Schulz steps leave the singular values of the returned matrix in a band around $1$ rather than at $1$, so its Frobenius norm falls below $\sqrt{r}$, by $5$--$22\%$ on our layer shapes; for the \textsc{lmo}-terminated methods $\eta_0$ is therefore the per-step RMS gain of the exact step, realized up to that shape-dependent factor. Substituting the two Frobenius norms of \eqref{eq:two_families},
\begin{equation}
\label{eq:unit_gain_closed}
\lambda_\ell^{\textsc{lmo}}=\sqrt{\max\!\left(1,\tfrac{m}{n}\right)},
\qquad
\lambda_\ell^{\textsc{sign}}=\frac{1}{\sqrt{n}} .
\end{equation}

\paragraph{Justification of \eqref{eq:unit_gain_closed}.} The first expression is exactly the aspect-ratio factor $\sqrt{\max(1,m/n)}$ present in the reference Muon implementation \citep{jordan2024muon}, which was introduced as a practical heuristic. It is also the point at which the rule parts company with the geometric alternative: taking the layer norm to be RMS$\to$RMS rather than spectral (\hyperref[app:prelim]{Appendix~\ref*{app:prelim}}) prescribes $\sqrt{m/n}$, which coincides with unit gain for $m\ge n$ and falls below it for $m<n$, where the two disagree and the implemented factor is the unit-gain one. The disagreement is the informative case: on a wide layer the RMS$\to$RMS ball shrinks the step in proportion to $\sqrt{m/n}$, whereas the gain \eqref{eq:gain_identity} of $\mathbf{U}\mathbf{V}^\top$ is already $\sqrt{\min(m,n)/m}=1$ there and needs no correction. Nor does the geometric route reach the \textsc{sign} family at all: a $\pm1$ matrix is the oracle output of no norm (\hyperref[th:1]{Theorems~\ref*{th:1}}--\ref{th:3}), so there is no unit ball whose radius could set its scale, while \eqref{eq:unit_gain} applies to it unchanged. The unit-gain criterion \emph{derives} the aspect factor, and also explains why Muon's step size is known to transfer across widths: its step has $\gamma=\eta$ independently of $n$, so no correction in the input dimension is needed. The second expression is the counterpart the sign family has never been given. Only $\eta_0$ is tuned, and it is now a shape-free quantity; the shape dependence is determined a priori.

Two consequences follow. First, $\lambda_\ell$ is a deterministic function of the layer shape, known to server and clients alike, so per-layer step sizes require \emph{no} communication and leave the one-bit-per-parameter budget intact. Second, on the CIFAR ResNet-18 of our experiments $\lambda_\ell^{\textsc{sign}}$ spans a factor of $13$ across layers (from $1/\sqrt{27}$ at the first convolution to $1/\sqrt{4608}$ in the last stage), so a single global rate is necessarily a compromise: roughly correct for the middle of the network, several-fold too large at the first convolution and too small at the last stage. The \textsc{lmo} family is exempt from this, which is one reason full-precision Muon is easier to tune than its sign-compressed variants.

\paragraph{Selecting the exponent.} Writing $\lambda_\ell^{\textsc{sign}}=n^{-a}$, the unit-gain rule is $a=\tfrac12$. Identity \eqref{eq:gain_identity} assumes the input independent of $\mathbf{A}$, the right model for a \emph{single} step; if the \emph{accumulated} update $\sum_t\eta\mathbf{P}_t$ aligns with the activations, its gain is $\Theta(\eta n)$ rather than $\Theta(\eta\sqrt{n})$, giving $a=1$, the $\mu$P rule $\eta\propto1/n$ for sign-like updates \citep{yang2021tuning}. A direct measurement decides between the regimes: we track the realized gain $\|\mathbf{X}_t-\mathbf{X}_0\|_F/\sqrt{m}$ at a \emph{constant} step size (under a decaying schedule the accumulation saturates and the fit reports the schedule) and fit its growth exponent $h$ in $t$, which is $\tfrac12$ for incoherent accumulation and $1$ for aligned. Over $20$ epochs the fit returns $\hat h = 0.513$ for Muon, $0.515$ for SignMuon, $0.490$ for SignSGD and $0.561$ for MuonSign, each at $R^2 \ge 0.999$. Muon is the control: \eqref{eq:unit_gain_closed} and $\mu$P prescribe the \textsc{lmo} family the identical multiplier, so its exponent is the diagnostic's reading when the rule is not in question, and the sign methods match it. All four values lie near $\tfrac12$ and none near $1$, so the accumulation is incoherent on this network and we adopt $a=\tfrac12$ for both families in every network experiment. The exponent fixes the shape dependence of $\lambda_\ell$ and the transfer of $\eta_0$ across widths; it makes no claim about which $a$ maximizes accuracy at one fixed width, where the choice is largely absorbed into $\eta_0$.

\paragraph{The placement of weight decay is not arbitrary.} The same scale invariance that makes $\lambda_\ell$ necessary also dictates \emph{where} an $\ell_2$ penalty may be applied: $\operatorname{sign}(c\mathbf{M})=\operatorname{sign}(\mathbf{M})$ and $\operatorname{polar}(c\mathbf{M})=\operatorname{polar}(\mathbf{M})$ for all $c>0$. Folding the decay into the gradient, $\widetilde{\mathbf{G}}_t=\mathbf{G}_t+\lambda_{\mathrm{wd}}\mathbf{X}_t$, the convention of \citet{mishra2026signmuon} and of most sign-method implementations, therefore supplies no contraction: the step length $\eta_t\lambda_\ell\|\mathbf{P}\|_F$ is fixed by \eqref{eq:two_families}, and the decay term can only rotate the direction. That rotation is governed by $\rho_t=\lambda_{\mathrm{wd}}\|\mathbf{X}_t\|_F/\|\mathbf{G}_t\|_F$, which drifts from negligible to $\Theta(1)$ as $\|\mathbf{G}_t\|_F$ falls over training and which depends on each method's own momentum scale, so one nominal $\lambda_{\mathrm{wd}}$ is a different perturbation for each method. Decoupled decay, $\mathbf{X}_{t+1}=(1-\eta_t\lambda_{\mathrm{wd}})\mathbf{X}_t-\eta_t\lambda_\ell\mathbf{P}_t$, is by contrast commensurate with the update under the unit-gain rule: its displacement has gain $\eta_t\lambda_{\mathrm{wd}}\gamma(\mathbf{X}_t)$ against the step's $\eta_t$, a ratio free of $\eta_t$, of the layer shape and of the method. We therefore decouple, and use the coupled form only in an ablation. This explains an observation of \citet{mishra2026signmuon}: sweeping $\lambda_{\mathrm{wd}}\in\{0,0.1,0.2\}$ coupled over $330$ CIFAR-10 ResNet-50 configurations, every Muon and Sign-Muon entry in their top ten (the only entries with a decay sweep) selects $\lambda_{\mathrm{wd}}=0$; at the two nonzero values the decay term dominates the gradient in $\widetilde{\mathbf{G}}_t$ for most of training, so the transmitted sign approaches $\operatorname{sign}(\mathbf{X}_t)$ and the sweep rejects this placement of the penalty, not regularization as such.

\paragraph{Scope of the analysis.} Our theorems are stated for unregularized $f$, so we report unregularized runs as the primary comparison and weight decay as an ablation; the reference nanoGPT configuration we build on also uses $\lambda_{\mathrm{wd}}=0$ for every parameter group. Two remarks delimit the gap. First, the \emph{coupled} form is covered verbatim: it is nothing other than running the same method on $f_{\lambda}=f+\tfrac{\lambda_{\mathrm{wd}}}{2}\|\mathbf{X}\|_F^2$, so every rate carries over with $L_i\mapsto L_i+\lambda_{\mathrm{wd}}r_i$. The rank factor is not slack in the bound: our smoothness is measured in the nuclear norm against a spectral-norm displacement (\hyperref[as:2]{Assumption~\ref*{as:2}}), and $\|\lambda_{\mathrm{wd}}\mathbf{Z}\|_*\le\lambda_{\mathrm{wd}}r_i\|\mathbf{Z}\|_{2\to2}$ is tight at $\mathbf{Z}=\mathbf{I}_{r_i}$; $L\mapsto L+\lambda_{\mathrm{wd}}$ would be the Euclidean statement, and these rates are not Euclidean. The paradox is that this is precisely the variant which does not regularize. Second, the \emph{decoupled} form is not covered by our rates, yet it furnishes something the analysis requires. Since $\|\mathbf{P}_t\|_F$ is a known constant and $\lambda_\ell\|\mathbf{P}_t\|_F=\sqrt{m}$ identically under \eqref{eq:unit_gain}, the triangle inequality gives $\|\mathbf{X}_{t+1}\|_F\le(1-\eta_t\lambda_{\mathrm{wd}})\|\mathbf{X}_t\|_F+\eta_t\sqrt{m}$, and hence, for any $\eta_t\lambda_{\mathrm{wd}}\le1$,
\begin{equation}
\gamma(\mathbf{X}_t)\;\le\;\max\!\left\{\gamma(\mathbf{X}_0),\;\lambda_{\mathrm{wd}}^{-1}\right\}
\quad\text{for all }t,
\label{eq:wd_ball}
\end{equation}
a bound on the layer's gain that is uniform in $t$ and \emph{independent of the layer shape}. Norm-constrained updates are what render this possible: for SGD the step length is data-dependent and no such a priori bound exists. Since layer-wise \textsc{lmo} analyses assume smoothness on a bounded region, \eqref{eq:wd_ball} is the statement that decoupled decay supplies that region. We do not claim a convergence rate for the decoupled variant itself.

\paragraph{Sensitivity to the rule.} Equation~\eqref{eq:unit_gain_main} is a heuristic, so we state its scope precisely. Each candidate rule, one global rate ($\lambda_\ell=1$), unit gain ($\lambda_\ell=n^{-1/2}$) or $\mu$P ($\lambda_\ell=n^{-1}$), shifts the selected $\eta_0$ by roughly the multiplier it prescribes; what would matter is a change in the \emph{ordering} of the methods. The exposure is bounded twice over: the \textsc{lmo} family cannot move, unit gain and $\mu$P prescribing it the identical multiplier, and the three \textsc{sign} methods are tuned and reported under one rule, so a wrong exponent rescales them alike. Neither consideration is a measurement, so \hyperref[tab:rule_ablation]{Table~\ref*{tab:rule_ablation}} re-tunes the three \textsc{sign} methods from scratch under each rule on federated CNN2, whose three matrix parameters span a factor of $7.8$ in $\lambda_\ell^{\textsc{sign}}$, and runs each selected rate at the reporting horizon.

Two things follow. The selected $\eta_0$ moves by roughly the multiplier the rule prescribes, which is the rule working and not a defect: measured against one global rate, unit gain raises SignMuon's rate by a factor of $50$ and $\mu$P by $10^3$, against the $\sqrt{n}\in[8.7,67.9]$ and $n\in[75,4608]$ that the three layer shapes prescribe. The ordering, meanwhile, does not move: SignMuon, then MuonSign, then SignSGD under every rule, the first ahead of the last by $4.3$ points under the global rate, $4.3$ under unit gain and $4.0$ under $\mu$P. Within a method the rules agree to within $0.2$ points for SignMuon and $0.1$ for SignSGD, at or below the seed spread; MuonSign is the one case where they separate at all, unit gain standing $0.57$ above $\mu$P and $0.68$ above the global rate, about two seed spreads, and in the direction that favours the rule we adopted. The sign-family ordering of \hyperref[tab:exp_3]{Table~\ref*{tab:exp_3}} therefore does not rest on the exponent, which the exposure argument above could only bound rather than establish.

\begin{table}[!t]
\centering\small
\setlength{\tabcolsep}{4pt}
\begin{tabular}{@{}llccr@{}}
\toprule
\textbf{Method} & \textbf{Rule} & $\boldsymbol{\lambda_\ell}$ & $\boldsymbol{\eta_0}$ & \textbf{Test acc (\%)} \\
\midrule
SignMuon & global            & $1$        & $0.002$  & $85.68 \pm 0.17$ \\
SignMuon & \textbf{unit gain} & $n^{-1/2}$ & $0.1$    & $85.72 \pm 0.24$ \\
SignMuon & $\mu$P            & $n^{-1}$   & $2$      & $85.52 \pm 0.01$ \\
\midrule
MuonSign & global            & $1$        & $0.001$  & $82.26 \pm 0.40$ \\
MuonSign & \textbf{unit gain} & $n^{-1/2}$ & $0.02$   & $82.94 \pm 0.19$ \\
MuonSign & $\mu$P            & $n^{-1}$   & $0.5$    & $82.37 \pm 0.27$ \\
\midrule
SignSGD  & global            & $1$        & $0.0005$ & $81.37 \pm 0.21$ \\
SignSGD  & \textbf{unit gain} & $n^{-1/2}$ & $0.01$   & $81.44 \pm 0.15$ \\
SignSGD  & $\mu$P            & $n^{-1}$   & $0.5$    & $81.47 \pm 0.25$ \\
\bottomrule
\end{tabular}
\caption{The per-layer rule ablation, on the federated CNN2 of \hyperref[tab:exp_3]{Table~\ref*{tab:exp_3}}. Each (method, rule) pair is re-tuned from scratch on the five-point lattice of \hyperref[app:repro]{Appendix~\ref*{app:repro}} and then run at $2000$ rounds: three seeds under the two alternatives, and under unit gain the five seeds \hyperref[tab:exp_3]{Table~\ref*{tab:exp_3}} reports.}
\label{tab:rule_ablation}
\end{table}

\paragraph{Comparison with concurrent work.} \citet{mishra2026signmuon} analyse the normalized update $\mathbf{D}_t=\bar{\mathbf{S}}_t/\sqrt{mn}$, justified by $\|\mathbf{D}_t\|_{\mathrm{op}}\le\|\mathbf{D}_t\|_F=1$ under their spectral-norm smoothness assumption, and remark that updating with $\bar{\mathbf{S}}_t$ directly is equivalent after absorbing $\sqrt{mn}$ into $\eta_t$. That equivalence holds for a single matrix but not across layers of differing shape, and their algorithm applies no shape factor, so their experiments use a single global rate as well. The substitution is also loose in a shape-dependent way: $\bar{\mathbf{S}}_t$ has rank at most $\min(m,n)$, so $\|\bar{\mathbf{S}}_t\|_F\le\sqrt{\min(m,n)}\,\|\bar{\mathbf{S}}_t\|_{\mathrm{op}}$ and the substitution of $\sqrt{mn}$ for the operator norm is loose by up to $\sqrt{\min(m,n)}$. That bound itself grows with depth, from $\sqrt{27}\approx5.2$ at the first convolution of a ResNet-18 to $\sqrt{512}\approx22.6$ in the last stage, so the spectral radius the analysis assigns to the step varies across the network instead of remaining uniform. \hyperref[tab:lrscale]{Table~\ref*{tab:lrscale}} sets the rules side by side.

\begin{table}[!t]
\centering\footnotesize
\setlength{\tabcolsep}{3pt}
\begin{tabular}{@{}llll@{}}
\toprule
\textbf{Rule} & \textbf{\textsc{lmo}} & \textbf{\textsc{sign}} & \textbf{Equalizes} \\
\midrule
global ($a\!=\!0$) & $1$ & $1$ & nothing \\
RMS$\to$RMS ball & $\sqrt{\tfrac{m}{n}}$ & --- & \textsc{lmo} trust region \\
Muon default & $\sqrt{\max(1,\tfrac{m}{n})}$ & $1$ & \textsc{lmo} gain \\
\textbf{unit gain} & $\sqrt{\max(1,\tfrac{m}{n})}$ & $n^{-1/2}$ & per-step gain, both \\
$\mu$P ($a\!=\!1$) & $\sqrt{\max(1,\tfrac{m}{n})}$ & $n^{-1}$ & aligned accumulation \\
Mishra et al. & $r^{-1/2}$ & $(mn)^{-1/2}$ & $\|\mathbf{P}\|_F$ \\
\bottomrule
\end{tabular}
\caption{Per-layer step-size multipliers $\lambda_\ell$. Only $\eta_0$ is tuned; $\lambda_\ell$ is fixed a priori by the shape. The last four rules differ only in the \textsc{sign} family, and the \textsc{lmo} column of the unit-gain rule coincides with the factor already employed in practice \citep{jordan2024muon}, which is our principal evidence that \eqref{eq:gain_def} is the correct criterion. The second row is the multiplier implied by taking the layer geometry to be RMS$\to$RMS rather than spectral: it departs from the other four for $m<n$ and, being a property of a unit ball, is undefined for steps that are the oracle of no ball. The last row is the normalization of \citet{mishra2026signmuon}, which equalizes $\|\mathbf{P}\|_F$ rather than the gain.}
\label{tab:lrscale}
\end{table}

\subsection{Algorithms}\label{app:alg}

% Federated protocol details moved here from the main text (Theory,
% "Federated Learning") to save main-part space.
\paragraph{Federated protocol.}
At the start of round $t$ each client $j$ holds the global model $\mathbf{X}_{t-1}$ and evaluates one stochastic gradient $\mathbf{G}_t^{(j)} = \nabla f_j(\mathbf{X}_{t-1};\xi_t^{(j)})$ at it; clients take no local parameter steps, so one round is one server step and no client-drift term arises. (The released runs accumulate three mini-batches of $64$ at fixed weights to save activation memory; the BatchNorm statistics being frozen (\hyperref[app:repro]{Appendix~\ref*{app:repro}}), the loss is separable across samples and the average is a gradient at batch $192$, except where a client's shard ends in a shorter mini-batch.) Each client updates its own momentum buffer, applies the LMO (\hyperref[alg:muon_lmo]{Algorithm~\ref*{alg:muon_lmo}}), and transmits the elementwise sign,
\begin{equation*}
\mathbf{M}_t^{(j)} = \mu \mathbf{M}_{t-1}^{(j)} + (1-\mu)\,\mathbf{G}_t^{(j)}, \qquad
\mathbf{D}_t^{(j)} = -A\bigl(\mathbf{M}_t^{(j)}\bigr), \qquad
\mathbf{s}_t^{(j)} = \operatorname{sign}\bigl(\mathbf{D}_t^{(j)}\bigr),
\end{equation*}
with the exponential-moving-average momentum of \eqref{eq:signa_general}, matching \hyperref[alg:fed_workerlmo]{Algorithms~\ref*{alg:fed_workerlmo}}--\ref{alg:fed_serverlmo}; the uplink is one bit per matrix parameter. The server aggregates by majority vote, $\mathbf{s}_t^{\mathrm{agg}}=\operatorname{sign}\bigl(\sum_{j=1}^N \mathbf{s}_t^{(j)}\bigr)$, which is $\pm1$ in each component: client messages are $\pm1$-valued by the convention of \hyperref[sec:theory]{Section~\ref*{sec:theory}}, so at an odd client count the vote cannot tie, and at an even count a tie is broken by a fair coin. Momentum having been applied at the clients, the server steps directly, $\mathbf{X}_{t} = \mathbf{X}_{t-1} - \eta_t \mathbf{s}_t^{\mathrm{agg}}$, and the vote rather than the model travels back down the link: every client applies the same $\pm1$-valued update to its local copy of the model, so client and server models remain identical and the downlink carries one bit per parameter as well (\hyperref[app:commacct]{Appendix~\ref*{app:commacct}}). The final classification layer is exempt from the rule and trained with AdamW.

\paragraph{Federated error feedback.}
The repair for the biased sign compressor is EF21 \citep{richtarik2021ef21} in the LMO form of \citet{gruntkowska2025error}: it changes the uplink message, and for EF21-MuonSign the downlink as well. Client $j$ compresses the residual $\Delta_t^{(j)}$ between its estimator and the quantity it would otherwise send, the polar factor $\mathbf{D}_t^{(j)}$ for EF21-SignMuon, whose oracle runs on the client, and the momentum $\tilde{\mathbf{M}}_t^{(j)}$ for EF21-MuonUSign and EF21-MuonSign, whose oracle runs on the server, and transmits the pair $(\operatorname{sign}(\Delta_t^{(j)}),\alpha_t^{(j)})$ with $\alpha_t^{(j)}=\operatorname{mean}|\Delta_t^{(j)}|$. The server accumulates these into a global estimator ($\mathbf{d}_t$ in \hyperref[alg:fed_workerlmo]{Algorithm~\ref*{alg:fed_workerlmo}}, $\mathbf{g}_t$ in \hyperref[alg:fed_serverlmo]{Algorithm~\ref*{alg:fed_serverlmo}}) and either steps along it or applies one LMO to it. The extra scalar is one full-precision number per matrix layer per round, so the uplink stays at ${\approx}1$ bit per parameter; the estimator itself is dense, so the downlink carries a full-precision model unless a second error-feedback loop compresses it, as EF21-MuonSign's does.

\begin{algorithm}[!t]
\algnarrow
\caption{MuonLMO}
\label{alg:muon_lmo}
\begin{algorithmic}[1]
\Statex \textbf{Input}: tensor $\mathbf{Y}$; Newton-Schulz coefficients $a=3.4445,\,b=4.7750,\,c=2.0315$; iteration count $\mathrm{ns\_steps}$
\Statex \textbf{Output}: polar factor $\mathbf{D}\approx\mathbf{U}\mathbf{V}^\top$ of $\mathbf{Y}$
\State $\mathbf{Y} \gets \text{ReshapeTo2D}(\mathbf{Y})$ \Comment{Flatten tensor to matrix $m \times n$}
\State $\mathbf{Y} \gets \mathbf{Y}/\|\mathbf{Y}\|_F$ \Comment{Normalize (optional)}
\For{$k = 1$ to $\mathrm{ns\_steps}$} \Comment{5th-order Newton-Schulz orthogonalization}
    \State $\mathbf{A} \gets \mathbf{Y}\mathbf{Y}^\top$
    \State $\mathbf{Y} \gets a\,\mathbf{Y} - b\,\mathbf{A}\mathbf{Y} + c\,\mathbf{A}^2\,\mathbf{Y}$
\EndFor
\State $\mathbf{D} \gets \text{ReshapeToOriginal}(\mathbf{Y})$ \Comment{Reshape back to the original tensor shape}
\State \Return $\mathbf{D}$
\end{algorithmic}
\end{algorithm}

\begin{algorithm}
\algnarrow
    \caption{SignMuon}
    \label{central_alg}
    \textbf{Input}: Initial model $\mathbf{X}_0$, momentum coefficient $\mu$, learning rate $\eta_t$ \\
    \textbf{Output}: Updated model $\mathbf{X}$
    \begin{algorithmic}[1]
    \State $\mathbf{M}_0 \gets 0$
    \For{$t = 1$ to $T$}
        \State $\mathbf{G}_t \gets \nabla f(\mathbf{X}_{t-1};\xi_t)$ \Comment{Stochastic gradient}
        \State $\mathbf{M}_t \gets \mu \mathbf{M}_{t-1} + (1-\mu)\,\mathbf{G}_t$ \Comment{Momentum accumulation}
        \State $\tilde{\mathbf{M}}_t =
        \begin{cases}
            \mathbf{M}_t, & \text{(default)}, \\
            (1-\mu)\,\mathbf{G}_t + \mu \mathbf{M}_t, & \text{(Nesterov)}
        \end{cases}$
        \State $\mathbf{D}_t \gets \mathrm{MuonLMO}(\tilde{\mathbf{M}}_t)$
        \State $\mathbf{s}_t^{\uparrow} \gets \operatorname{sign}(\mathbf{D}_t)$ \Comment{Uplink sign compression}
        \State $\mathbf{X}_{t} \gets \mathbf{X}_{t-1} - \eta_t \mathbf{s}_t^{\uparrow}$ \Comment{Update parameters}
    \EndFor
    \end{algorithmic}
\end{algorithm}

\begin{algorithm}
\algnarrow
    \caption{EF21-SignMuon}
    \label{ef21_signmuon}
    \textbf{Input}: Initial model $\mathbf{X}_0$, momentum coefficient $\mu$, learning rate $\eta_t$ \\
    \textbf{Output}: Updated model $\mathbf{X}$
    \begin{algorithmic}[1]
    \State $\mathbf{M}_0 \gets 0,\quad \mathbf{d}_0^{\mathrm{est}} \gets 0$
    \For{$t = 1$ to $T$}
        \State $\mathbf{G}_t \gets \nabla f(\mathbf{X}_{t-1};\xi_t)$ \Comment{Stochastic gradient}
        \State $\mathbf{M}_t \gets \mu \mathbf{M}_{t-1} + (1-\mu)\,\mathbf{G}_t$ \Comment{Momentum accumulation (EMA)}
        \State $\tilde{\mathbf{M}}_t =
        \begin{cases}
            \mathbf{M}_t, & \text{(default)}, \\
            (1-\mu)\,\mathbf{G}_t + \mu \mathbf{M}_t, & \text{(Nesterov)}
        \end{cases}$
        \State $\mathbf{D}_t \gets \mathrm{MuonLMO}(\tilde{\mathbf{M}}_t)$
        \State $\Delta_t^{\uparrow} \gets \mathbf{D}_t - \mathbf{d}_{t-1}^{\mathrm{est}}$ \Comment{Uplink residual (polar factor)}
        \State $\alpha_t^{\uparrow} \gets \operatorname{mean}(|\Delta_t^{\uparrow}|)$
        \State $\mathbf{d}_{t}^{\mathrm{est}} \gets \mathbf{d}_{t-1}^{\mathrm{est}} + \alpha_t^{\uparrow}\operatorname{sign}(\Delta_t^{\uparrow})$ \Comment{Uplink EF21}
        \State $\mathbf{X}_{t} \gets \mathbf{X}_{t-1} - \eta_t \mathbf{d}_{t}^{\mathrm{est}}$ \Comment{Update parameters}
    \EndFor
    \end{algorithmic}
\end{algorithm}

\begin{algorithm}[!h]
\algnarrow
    \caption{MuonUSign}
    \label{alg:muon_usign}
    \textbf{Input}: Initial model $\mathbf{X}_0$, momentum coefficient $\mu$, learning rate $\eta_t$ \\
    \textbf{Output}: Updated model $\mathbf{X}$
    \begin{algorithmic}[1]
    \State $\mathbf{M}_0 \gets 0$
    \For{$t = 1$ to $T$}
        \State $\mathbf{G}_t \gets \nabla f(\mathbf{X}_{t-1};\xi_t)$ \Comment{Stochastic gradient}
        \State $\mathbf{M}_t \gets \mu \mathbf{M}_{t-1} + (1-\mu)\,\mathbf{G}_t$ \Comment{Momentum accumulation}
        \State $\tilde{\mathbf{M}}_t =
        \begin{cases}
            \mathbf{M}_t, & \text{(default)}, \\
            (1-\mu)\,\mathbf{G}_t + \mu \mathbf{M}_t, & \text{(Nesterov)}
        \end{cases}$
        \State $\mathbf{s}_t^{\uparrow} \gets \operatorname{sign}(\tilde{\mathbf{M}}_t)$ \Comment{Uplink sign compression}
        \State $\mathbf{D}_t \gets \mathrm{MuonLMO}(\mathbf{s}_t^{\uparrow})$
        \State $\mathbf{X}_{t} \gets \mathbf{X}_{t-1} - \eta_t \mathbf{D}_t$ \Comment{Update parameters}
    \EndFor
    \end{algorithmic}
\end{algorithm}
    
\begin{algorithm}[!h]
\algnarrow
    \caption{MuonSign}
    \label{alg:muon_sign}
    \textbf{Input}: Initial model $\mathbf{X}_0$, momentum coefficient $\mu$, learning rate $\eta_t$ \\
    \textbf{Output}: Updated model $\mathbf{X}$
    \begin{algorithmic}[1]
    \State $\mathbf{M}_0 \gets 0$
    \For{$t = 1$ to $T$}
        \State $\mathbf{G}_t \gets \nabla f(\mathbf{X}_{t-1};\xi_t)$ \Comment{Stochastic gradient}
        \State $\mathbf{M}_t \gets \mu \mathbf{M}_{t-1} + (1-\mu)\,\mathbf{G}_t$ \Comment{Momentum accumulation}
        \State $\tilde{\mathbf{M}}_t =
        \begin{cases}
            \mathbf{M}_t, & \text{(default)}, \\
            (1-\mu)\,\mathbf{G}_t + \mu \mathbf{M}_t, & \text{(Nesterov)}
        \end{cases}$
        \State $\mathbf{s}_t^{\uparrow} \gets \operatorname{sign}(\tilde{\mathbf{M}}_t)$ \Comment{Uplink sign compression}
        \State $\mathbf{D}_t \gets \mathrm{MuonLMO}(\mathbf{s}_t^{\uparrow})$
        \State $\mathbf{s}_t^{\downarrow} \gets \operatorname{sign}(\mathbf{D}_t)$ \Comment{Downlink sign compression}
        \State $\mathbf{X}_{t} \gets \mathbf{X}_{t-1} - \eta_t \mathbf{s}_t^{\downarrow}$ \Comment{Update parameters}
    \EndFor
    \end{algorithmic}
\end{algorithm}
    
\begin{algorithm}[!t]
\algnarrow
    \caption{EF21-MuonUSign}
    \label{central_alg_ef}
    \textbf{Input}: Initial model $\mathbf{X}_0$, momentum coefficient $\mu$, learning rate $\eta_t$ \\
    \textbf{Output}: Updated model $\mathbf{X}$
    \begin{algorithmic}[1]
    \State $\mathbf{M}_0 \gets 0,\quad \mathbf{g}_0^{\mathrm{est}} \gets 0$
    \For{$t = 1$ to $T$}
        \State $\mathbf{G}_t \gets \nabla f(\mathbf{X}_{t-1};\xi_t)$ \Comment{Stochastic gradient}
        \State $\mathbf{M}_t \gets \mu\mathbf{M}_{t-1} + (1-\mu)\,\mathbf{G}_t$ \Comment{Momentum accumulation}
        \State $\tilde{\mathbf{M}}_t =
        \begin{cases}
            \mathbf{M}_t, & \text{(default)}, \\
            (1-\mu)\,\mathbf{G}_t + \mu\mathbf{M}_t, & \text{(Nesterov)}
        \end{cases}$
        \State $\Delta_t^{\uparrow} \gets \tilde{\mathbf{M}}_t - \mathbf{g}_{t-1}^{\mathrm{est}}$ \Comment{Uplink residual}
        \State $\alpha_t^{\uparrow} \gets \operatorname{mean}(|\Delta_t^{\uparrow}|)$
        \State $\mathbf{g}_{t}^{\mathrm{est}} \gets \mathbf{g}_{t-1}^{\mathrm{est}} + \alpha_t^{\uparrow}\operatorname{sign}(\Delta_t^{\uparrow})$ \Comment{Uplink EF21}
        \State $\mathbf{D}_t \gets \mathrm{MuonLMO}(\mathbf{g}_{t}^{\mathrm{est}})$
        \State $\mathbf{X}_{t} \gets \mathbf{X}_{t-1} - \eta_t\mathbf{D}_t$ \Comment{Update parameters}
    \EndFor
    \end{algorithmic}
\end{algorithm}

\begin{algorithm}[!t]
\algnarrow
    \caption{EF21-MuonSign}
    \label{central_alg_ud}
    \textbf{Input}: Initial model $\mathbf{X}_0=\mathbf{W}_0$, momentum coefficient $\mu$, learning rate $\eta_t$ \\
    \textbf{Output}: Updated model $\mathbf{X}$
    \begin{algorithmic}[1]
    \State $\mathbf{M}_0 \gets 0,\quad \mathbf{g}_0^{\mathrm{est}} \gets 0$
    \For{$t = 1$ to $T$}
        \State $\mathbf{G}_t \gets \nabla f(\mathbf{W}_{t-1};\xi_t)$ \Comment{At the broadcast model}
        \State $\mathbf{M}_t \gets \mu\mathbf{M}_{t-1} + (1-\mu)\,\mathbf{G}_t$ \Comment{Momentum accumulation}
        \State $\tilde{\mathbf{M}}_t =
        \begin{cases}
            \mathbf{M}_t, & \text{(default)}, \\
            (1-\mu)\,\mathbf{G}_t + \mu\mathbf{M}_t, & \text{(Nesterov)}
        \end{cases}$
        \State $\Delta_t^{\uparrow} \gets \tilde{\mathbf{M}}_t - \mathbf{g}_{t-1}^{\mathrm{est}}$ \Comment{Uplink residual}
        \State $\alpha_t^{\uparrow} \gets \operatorname{mean}(|\Delta_t^{\uparrow}|)$
        \State $\mathbf{g}_{t}^{\mathrm{est}} \gets \mathbf{g}_{t-1}^{\mathrm{est}} + \alpha_t^{\uparrow}\operatorname{sign}(\Delta_t^{\uparrow})$ \Comment{Uplink EF21}
        \State $\mathbf{D}_t \gets \mathrm{MuonLMO}(\mathbf{g}_{t}^{\mathrm{est}})$
        \State $\mathbf{X}_{t} \gets \mathbf{X}_{t-1} - \eta_t\mathbf{D}_t$ \Comment{Update parameters (server)}
        \State $\Delta_t^{\downarrow} \gets \mathbf{X}_{t} - \mathbf{W}_{t-1}$ \Comment{Downlink residual}
        \State $\alpha_t^{\downarrow} \gets \operatorname{mean}(|\Delta_t^{\downarrow}|)$
        \State $\mathbf{W}_{t} \gets \mathbf{W}_{t-1} + \alpha_t^{\downarrow}\operatorname{sign}(\Delta_t^{\downarrow})$ \Comment{Downlink EF21-P}
    \EndFor
    \end{algorithmic}
\end{algorithm}

All six federated methods are instances of just two templates, separated by \emph{where the Muon LMO is evaluated}. When the sign acts \emph{after} the LMO (the SignMuon family), each client must orthogonalize locally, so the LMO runs on the worker and the client transmits a compressed \emph{direction} (\hyperref[alg:fed_workerlmo]{Algorithm~\ref*{alg:fed_workerlmo}}). When the sign acts \emph{before} the LMO (the MuonUSign/MuonSign family), the client transmits a compressed \emph{gradient}, and the server reconstructs it and applies a single LMO (\hyperref[alg:fed_serverlmo]{Algorithm~\ref*{alg:fed_serverlmo}}). Within each template, a method is fixed by its uplink compressor $\mathcal{C}^{\uparrow}\in\{\text{sign},\text{EF21}\}$ and downlink compressor $\mathcal{C}^{\downarrow}\in\{\text{exact},\text{sign},\text{EF21-P}\}$; \hyperref[tab:fed_master]{Table~\ref*{tab:fed_master}} lists the six instantiations.

The two uplinks aggregate differently, and each aggregation is forced. The EF21 uplink \emph{averages} the decompressed messages, $\mathbf{g}_t=\mathbf{g}_{t-1}+\tfrac1N\sum_j\alpha_t^{(j)}\operatorname{sign}(\Delta_t^{(j)})$, as \citet[Algorithm~3]{gruntkowska2025error} prescribe and as the reduction of \hyperref[app:hyp]{Appendix~\ref*{app:hyp}} requires; replacing that average by a vote would leave the framework and forfeit \hyperref[thm:conv]{Theorem~\ref*{thm:conv}}. The plain sign uplink instead takes a majority vote, $\operatorname{sign}(\sum_j\mathbf{s}_t^{(j)})$, \emph{before} the server LMO. Voting is what keeps the oracle's argument a $\pm1$ matrix, so that the server-side method is exactly the centralized MuonUSign, $\operatorname{polar}(\operatorname{sign}(\cdot))$ of \eqref{eq:three_placements}, evaluated at the aggregated sign; averaging would feed $\operatorname{polar}$ an argument valued in $\{-1,-1+2/N,\dots,1\}$ and define a different method, one that agrees with MuonUSign only at $N=1$. The choice also matches the sign-compression literature it inherits from \citep{bernstein2019signsgdmajority}. It carries no consequence for the downlink of this family, the polar factor being dense either way. As in the centralized setting, both templates are applied per matrix parameter, while vector parameters and the final classification layer are optimized with AdamW.

\begin{table}[!t]
    \centering\small
    \setlength{\tabcolsep}{5pt}
    \begin{tabular}{@{}lccc@{}}
    \toprule
    \textbf{Method} & \textbf{LMO} & \textbf{Uplink $\mathcal{C}^{\uparrow}$} & \textbf{Downlink $\mathcal{C}^{\downarrow}$} \\
    \midrule
    SignMuon        & worker & sign\,/\,MV & exact  \\
    EF21-SignMuon   & worker & EF21        & exact  \\
    \addlinespace
    MuonUSign       & server & sign\,/\,MV & exact  \\
    MuonSign        & server & sign\,/\,MV & sign   \\
    EF21-MuonUSign  & server & EF21        & exact  \\
    EF21-MuonSign   & server & EF21        & EF21-P \\
    \bottomrule
    \end{tabular}
    \caption{The six federated methods as instantiations of the two templates: \hyperref[alg:fed_workerlmo]{Algorithm~\ref*{alg:fed_workerlmo}} (worker-side LMO; rows~1--2) and \hyperref[alg:fed_serverlmo]{Algorithm~\ref*{alg:fed_serverlmo}} (server-side LMO; rows~3--6). Each method is fixed by the LMO location and the uplink/downlink compressors. ``MV'': majority vote; ``EF21-P'': primal (model-side) error feedback; ``exact'': the server applies no downlink compressor. That is not the same as a full-precision downlink: SignMuon's server has nothing to compress because the object it distributes, the majority vote, is already $\pm1$-valued, so its downlink is one bit per parameter all the same (\hyperref[app:commacct]{Appendix~\ref*{app:commacct}}). The three methods that do broadcast a dense model are MuonUSign, EF21-SignMuon and EF21-MuonUSign.}
    \label{tab:fed_master}
\end{table}

\begin{algorithm}[!t]
\algnarrow
    \caption{Federated SignMuon / EF21-SignMuon (worker-side LMO)}
    \label{alg:fed_workerlmo}
    \textbf{Input}: initial model $\mathbf{X}_0$; clients $N$; rounds $T$; learning rate $\eta_t$; momentum $\mu$; uplink compressor $\mathcal{C}^{\uparrow}\in\{\text{sign},\text{EF21}\}$ (\hyperref[tab:fed_master]{Table~\ref*{tab:fed_master}}) \\
    \textbf{Output}: global model $\mathbf{X}_T$
    \begin{algorithmic}[1]
    \State $\mathbf{M}_0^{(j)} \gets 0$, $\mathbf{d}_0^{(j)} \gets 0$ for all $j$; \quad $\mathbf{d}_0 \gets 0$
    \State \textbf{broadcast} $\mathbf{X}_0$ \emph{once}; every client keeps a local copy, refreshed below from the downlink message alone
    \For{$t = 1$ to $T$}
        \For{$j = 1$ to $N$ \textbf{in parallel}} \Comment{client $j$, holding $\mathbf{X}_{t-1}$}
            \State $\mathbf{G}_t^{(j)} \gets \nabla f_j(\mathbf{X}_{t-1};\,\xi_{t}^{(j)})$
            \State $\mathbf{M}_t^{(j)} \gets \mu \mathbf{M}_{t-1}^{(j)} + (1-\mu)\,\mathbf{G}_t^{(j)}$
            \State $\tilde{\mathbf{M}}_t^{(j)} \gets \mathbf{M}_t^{(j)}$ \Comment{or $(1-\mu)\mathbf{G}_t^{(j)}+\mu\mathbf{M}_t^{(j)}$ (Nesterov)}
            \State $\mathbf{D}_t^{(j)} \gets \mathrm{MuonLMO}\bigl(\tilde{\mathbf{M}}_t^{(j)}\bigr)$ \Comment{LMO on the client}
            \If{$\mathcal{C}^{\uparrow} = \text{EF21}$} \Comment{EF21-SignMuon}
                \State $\Delta_t^{(j)} \gets \mathbf{D}_t^{(j)} - \mathbf{d}_{t-1}^{(j)}$;\quad $\alpha_t^{(j)} \gets \operatorname{mean}(|\Delta_t^{(j)}|)$
                \State $\mathbf{d}_{t}^{(j)} \gets \mathbf{d}_{t-1}^{(j)} + \alpha_t^{(j)}\operatorname{sign}(\Delta_t^{(j)})$
                \State \textbf{send} $\bigl(\operatorname{sign}(\Delta_t^{(j)}),\, \alpha_t^{(j)}\bigr)$
            \Else \Comment{SignMuon}
                \State \textbf{send} $\mathbf{s}_t^{(j)} \gets \operatorname{sign}(\mathbf{D}_t^{(j)})$
            \EndIf
        \EndFor
        \Statex \quad\emph{on the server:}
        \If{$\mathcal{C}^{\uparrow} = \text{EF21}$} \Comment{EF21-SignMuon: $\mathbf{d}_t$ is dense, 32 bits down}
            \State $\mathbf{d}_t \gets \mathbf{d}_{t-1} + \tfrac{1}{N}\sum_{j} \alpha_t^{(j)}\operatorname{sign}(\Delta_t^{(j)})$ \Comment{aggregate direction}
            \State $\mathbf{X}_{t} \gets \mathbf{X}_{t-1} - \eta_t\,\mathbf{d}_t$;\; \textbf{broadcast} $\mathbf{X}_t$
        \Else \Comment{SignMuon: the vote is $\pm1$, 1 bit down}
            \State $\hat{\mathbf{s}}_t \gets \operatorname{sign}\bigl(\sum_{j} \mathbf{s}_t^{(j)}\bigr)$ \Comment{majority vote}
            \State $\mathbf{X}_{t} \gets \mathbf{X}_{t-1} - \eta_t\,\hat{\mathbf{s}}_t$;\; \textbf{broadcast} $\hat{\mathbf{s}}_t$ \Comment{clients apply the same step}
        \EndIf
    \EndFor
    \end{algorithmic}
\end{algorithm}

\begin{algorithm}[!t]
\algnarrow
    \caption{Federated MuonUSign / MuonSign / EF21-MuonUSign / EF21-MuonSign (server-side LMO)}
    \label{alg:fed_serverlmo}
    \textbf{Input}: initial model $\mathbf{X}_0$; clients $N$; rounds $T$; learning rate $\eta_t$; momentum $\mu$; compressors $\mathcal{C}^{\uparrow}\in\{\text{sign},\text{EF21}\}$, $\mathcal{C}^{\downarrow}\in\{\text{exact},\text{sign},\text{EF21-P}\}$ (\hyperref[tab:fed_master]{Table~\ref*{tab:fed_master}}) \\
    \textbf{Output}: global model $\mathbf{X}_T$
    \begin{algorithmic}[1]
    \State $\mathbf{M}_0^{(j)} \gets 0$, $\mathbf{g}_0^{(j)} \gets 0$ for all $j$; \quad $\mathbf{g}_0 \gets 0$
    \State \textbf{broadcast} $\mathbf{X}_0$ \emph{once}; every client holds the broadcast model $\mathbf{W}_0 \gets \mathbf{X}_0$, refreshed below from the downlink message alone ($\mathbf{W}_t=\mathbf{X}_t$ unless $\mathcal{C}^{\downarrow}=\text{EF21-P}$)
    \For{$t = 1$ to $T$}
        \For{$j = 1$ to $N$ \textbf{in parallel}} \Comment{client $j$, holding $\mathbf{W}_{t-1}$}
            \State $\mathbf{G}_t^{(j)} \gets \nabla f_j(\mathbf{W}_{t-1};\,\xi_{t}^{(j)})$
            \State $\mathbf{M}_t^{(j)} \gets \mu \mathbf{M}_{t-1}^{(j)} + (1-\mu)\,\mathbf{G}_t^{(j)}$
            \State $\tilde{\mathbf{M}}_t^{(j)} \gets \mathbf{M}_t^{(j)}$ \Comment{or $(1-\mu)\mathbf{G}_t^{(j)}+\mu\mathbf{M}_t^{(j)}$ (Nesterov)}
            \If{$\mathcal{C}^{\uparrow} = \text{EF21}$} \Comment{EF21-MuonUSign / EF21-MuonSign}
                \State $\Delta_t^{(j)} \gets \tilde{\mathbf{M}}_t^{(j)} - \mathbf{g}_{t-1}^{(j)}$;\quad $\alpha_t^{(j)} \gets \operatorname{mean}(|\Delta_t^{(j)}|)$
                \State $\mathbf{g}_{t}^{(j)} \gets \mathbf{g}_{t-1}^{(j)} + \alpha_t^{(j)}\operatorname{sign}(\Delta_t^{(j)})$
                \State \textbf{send} $\bigl(\operatorname{sign}(\Delta_t^{(j)}),\, \alpha_t^{(j)}\bigr)$
            \Else \Comment{MuonUSign / MuonSign}
                \State \textbf{send} $\mathbf{s}_t^{(j)} \gets \operatorname{sign}\bigl(\tilde{\mathbf{M}}_t^{(j)}\bigr)$
            \EndIf
        \EndFor
        \Statex \quad\emph{on the server:}
        \If{$\mathcal{C}^{\uparrow} = \text{EF21}$}
            \State $\mathbf{g}_t \gets \mathbf{g}_{t-1} + \tfrac{1}{N}\sum_{j} \alpha_t^{(j)}\operatorname{sign}(\Delta_t^{(j)})$ \Comment{reconstruct gradient}
        \Else
            \State $\mathbf{g}_t \gets \operatorname{sign}\bigl(\sum_{j} \mathbf{s}_t^{(j)}\bigr)$ \Comment{majority vote}
        \EndIf
        \State $\mathbf{D}_t \gets \mathrm{MuonLMO}(\mathbf{g}_t)$ \Comment{single LMO on the server}
        \If{$\mathcal{C}^{\downarrow} = \text{sign}$} \Comment{MuonSign: 1 bit/param down}
            \State $\mathbf{X}_{t} \gets \mathbf{X}_{t-1} - \eta_t \operatorname{sign}(\mathbf{D}_t)$;\; \textbf{broadcast} $\operatorname{sign}(\mathbf{D}_t)$
            \State $\mathbf{W}_{t} \gets \mathbf{W}_{t-1} - \eta_t\operatorname{sign}(\mathbf{D}_t)$ \Comment{$=\mathbf{X}_t$; one model}
        \ElsIf{$\mathcal{C}^{\downarrow} = \text{EF21-P}$} \Comment{EF21-MuonSign: 1 bit/param down}
            \State $\mathbf{X}_{t} \gets \mathbf{X}_{t-1} - \eta_t \mathbf{D}_t$ \Comment{server model}
            \State $\Delta_t^{\downarrow} \gets \mathbf{X}_{t} - \mathbf{W}_{t-1}$;\quad $\alpha_t^{\downarrow} \gets \operatorname{mean}(|\Delta_t^{\downarrow}|)$
            \State $\mathbf{W}_{t} \gets \mathbf{W}_{t-1} + \alpha_t^{\downarrow}\operatorname{sign}(\Delta_t^{\downarrow})$;\; \textbf{broadcast} $\bigl(\operatorname{sign}(\Delta_t^{\downarrow}),\alpha_t^{\downarrow}\bigr)$ \Comment{clients apply the same refresh}
        \Else \Comment{exact downlink: MuonUSign / EF21-MuonUSign}
            \State $\mathbf{X}_{t} \gets \mathbf{X}_{t-1} - \eta_t \mathbf{D}_t$;\; \textbf{broadcast} $\mathbf{X}_t$;\quad $\mathbf{W}_{t} \gets \mathbf{X}_t$
        \EndIf
    \EndFor
    \end{algorithmic}
\end{algorithm}

\end{document}